\RequirePackage{ifpdf}
\ifpdf 
\documentclass[pdftex]{sigma}
\else
\documentclass{sigma}
\fi

\newcommand{\mh}[1]{\mbox{\huge{$#1$}}}

\usepackage[all]{xy}

\numberwithin{equation}{section}
\newtheorem{thm}{Theorem}[section]
\newtheorem{prop}[thm]{Proposition}
\newtheorem{cor}[thm]{Corollary}
\newtheorem{lem}[thm]{Lemma}

\theoremstyle{definition}

\newtheorem{dfn}[thm]{Definition}
\newtheorem{rmk}[thm]{Remark}
\newtheorem{nt}[thm]{Notation}

\newcommand{\Rep}{\text{\rm Rep}}

\newcommand{\cHil}{\mathcal{H}}
\newcommand{\nphi}{\mathfrak{n}_\varphi}

\newcommand{\ICM}{\text{\rm IC}(M)}
\newcommand{\Dom}{\text{\rm Dom}}
\newcommand{\ICMM}{\text{\rm IC}(M, M_1)}

\newcommand{\Ind}{{\rm Ind}\,}

\newcommand{\IR}{\text{\rm IR}}
\newcommand{\ep}{\varepsilon}
\newcommand{\sgn}{\text{\rm sgn}}

\newcommand{\rFs}[5]{\,_{#1}\varphi_{#2} \left( \genfrac{.}{.}{0pt}{}{#3}{#4}
\ ;#5 \right)}

\begin{document}

\allowdisplaybreaks

\renewcommand{\thefootnote}{$\star$}

\renewcommand{\PaperNumber}{087}

\FirstPageHeading

\ShortArticleName{Spherical Fourier Transforms on Locally Compact Quantum Gelfand Pairs}

\ArticleName{Spherical Fourier Transforms\\ on Locally Compact Quantum Gelfand Pairs\footnote{This paper is a
contribution to the Special Issue ``Relationship of Orthogonal Polynomials and Special Functions with Quantum Groups and Integrable Systems''. The
full collection is available at
\href{http://www.emis.de/journals/SIGMA/OPSF.html}{http://www.emis.de/journals/SIGMA/OPSF.html}}}

\Author{Martijn CASPERS}

\AuthorNameForHeading{M.~Caspers}

\Address{Radboud Universiteit Nijmegen, IMAPP, Heyendaalseweg 135, \\
6525 AJ Nijmegen, The Netherlands}
\Email{\href{mailto:caspers@math.ru.nl}{caspers@math.ru.nl}}
\URLaddress{\url{http://www.math.ru.nl/~caspers/}}

\ArticleDates{Received April 14, 2011, in f\/inal form August 30, 2011;  Published online September 06, 2011}

\Abstract{We study Gelfand pairs for locally compact quantum groups. We give an operator algebraic interpretation and show that the quantum Plancherel transformation restricts to a spherical Plancherel transformation. As an example, we turn the quantum group analogue of the normaliser of $SU(1,1)$ in $SL(2,\mathbb{C}$) together with its diagonal subgroup into a pair for which every irreducible corepresentation admits at most two vectors that are invariant with respect to the quantum subgroup. Using a $\mathbb{Z}_2$-grading, we obtain product formulae for little $q$-Jacobi functions.}


\Classification{16T99; 43A90}

\section{Introduction}

In the classical setting of locally compact groups, a Gelfand pair consists of a locally compact group $G$, together with a compact subgroup $K$ such that the convolution algebra of bi-$K$-invariant $L^1$-functions on $G$ is commutative. See \cite{Dij} or \cite{Far} for a comprehensive introduction. Gelfand pairs give rise to spherical functions and a spherical Fourier transform which decomposes bi-$K$-invariant functions on $G$ as an integral of spherical functions, see \cite[Theorem 6.4.5]{Dij} or \cite[Th\'eor\`eme IV.2]{Far}.

For many examples, this decomposition is made precise \cite{Dij}. The examples include the group of motions of the plane together with its diagonal subgroup and the pair $(SO_0(1,n), SO(n))$, where $SO_0(1,n)$ is the connected component of the identity of $SO(1,n)$. In particular the spherical functions are determined and one can derive product formulae for these type of functions.

Since the introduction of quantum groups, Gelfand pairs were studied in a quantum context, see for example \cite{Flo, PodVai, VaiGel, VaiHyp} and also the references given there. These papers  consider  pairs of quantum groups that are both compact. For such pairs it suf\/f\/ices to stay with a purely \mbox{(Hopf-)}al\-gebraic approach. Under the assumption that every irreducible unitary corepresentation admits only one matrix element that is invariant under both the left and right action of the subgroup, these quantum groups are called (quantum) Gelfand pairs. Classically, this is equiva\-lent to the commutativity assumption on the convolution algebra of bi-$K$-invariant elements. If the matrix coef\/f\/icients form a commutative algebra one speaks of a strict (quantum) Gelfand pair. In the group setting every Gelfand pair is automatically strict and as such strictness is a~purely non-commutative phenomenon.

For quantum groups, many deformations of classical Gelfand pairs do indeed form a quantum Gelfand pair that moreover is strict. As a compact example, $(SU_q(n), U_q(n-1))$ forms a strict Gelfand pair \cite{VaiHyp}. In a separate paper \cite{VaiE2} Vainerman introduces the quantum group of motions of the plane, together with the circle as a subgroup as an example of a Gelfand pair of which the larger quantum group is non-compact.  As a result a product formula for the Hahn--Exton $q$-Bessel functions, also known as $_1\varphi_1$ $q$-Bessel functions, is obtained \cite[ Corollary, p.~324]{VaiE2}, see also \cite[Corollary~6.4]{Koe}. However, a comprehensive general framework of quantum Gelfand pairs in the non-compact operator algebraic setting was unavailable at that time.

At the turn of the millennium, locally compact (l.c.) quantum groups have been put in an operator algebraic setting by
Kustermans and Vaes in their papers \cite{KusVae, KusVaeII},
see also \cite{KusLec, Tim, DaeleLcqg}. The def\/initions give a C$^\ast$-algebraic and a von Neumann algebraic interpretation of locally compact quantum groups. Many aspects of abstract harmonic analysis have found a suitable interpretation in this von Neumann algebraic framework.  In particular, Desmedt proved in his thesis \cite{Des} that there is an analogue of the Plancherel theorem, which gives a decomposition of the left regular corepresentation of a l.c.\ quantum group.

From this perspective, it is a natural question if the study of Gelfand pairs can be continued in the l.c.\ operator algebraic setting. In this paper we give this interpretation. Motivated by Desmedt's proof of the quantum Plancherel theorem, we def\/ine the necessary structures to obtain a classical Plancherel--Godement theorem \cite[Theorem 6.4.5]{Dij} or \cite[Th\'eor\`eme IV.2]{Far}. For this the operator algebraic interpretation of Gelfand pairs is essential.

We keep the setting a bit more general than one would expect. For a classical Gelfand pair of groups, one can prove that the larger group is unimodular from the commutativity assumption on bi-$K$-invariant elements. Here we will study pairs of quantum groups for which the smaller quantum group is compact and we {\it assume} that the larger group is unimodular. We will not impose the classically stronger commutativity assumption. The reason for this is that we would like to study $SU_q(1,1)_{\rm ext}$ together with its diagonal subgroup. However, the natural analogue of the commutativity assumption would exclude this example.

We mention that it is known that the notion of a quantum subgroup is in a sense too restrictive. Using Koornwinder's twisted primitive elements, it is possible to def\/ine double coset spaces associated with $SU_q(2)$ and get so called $(\sigma, \tau)$-spherical elements, see \cite{KoeAsk} for this particular example. See also \cite{KoeSto} for a similar study of $SU_q(1,1)$ on an algebraic level. The subgroup setting then corresponds to the limiting case $\sigma, \tau \rightarrow \infty$. In the present paper we do not incorporate such a general setting.

Motivated by the Hopf-algebraic framework, we introduce the non-compact analogues of bi-$K$-invariant functions and its dual \cite{Flo, VaiHyp} and equip these with weights. We will do this in a von Neumann algebraic manner and for the dual structure also in a C$^\ast$-algebraic manner. We prove that the C$^\ast$-algebraic weight lifts to the von Neumann algebraic weight. Moreover, we establish a~spherical analogue of a theorem by Kustermans~\cite{Kus} which establishes a correspondence between representations of the (universal) C$^\ast$-algebraic dual quantum group and corepresentations of the quantum group itself.  Eventually, this structure culminates in a quantum Plancherel--Godement theorem, as an application of \cite[Theorem~3.4.5]{Des}. This illustrates the advantage of an operator algebraic interpretation above the Hopf algebraic approach. In particular, we get a~spherical $L^2$-Fourier transform, or spherical Plancherel transformation, and we show in principle that this is a restriction of the non-spherical Plancherel transformation.

As an example, we treat the f\/irst example of a quantum Gelfand pair involving a $q$-deforma\-tion of $SU(n,1)$. Namely, we treat the quantum analogue of the normalizer of $SU(1,1)$ in $SL(2, \mathbb{C})$, which we denote by $SU_q(1,1)_{\rm ext}$, see \cite{KoeKus} and \cite{GroKoeKus}. We identify the circle as its diagonal subgroup and study the spherical properties of this pair.
We  see that the classical commutativity assumption on the convolution algebra is too restrictive to capture $SU_q(1,1)_{\rm ext}$ with its diagonal subgroup. Nevertheless, the pair exhibits properties reminiscent of classical Gelfand pairs. In particular, we  see how one can  derive product formulae using gradings on this quantum group and its dual.

We mention that the $q$-deformation of $SU(1,1)$ was f\/irst established on the operator algebraic level in~\cite{KoeKus}. The construction heavily relies on $q$-analysis. More recently, de Commer~\cite{ComSU} was able to obtain $SU_q(1,1)_{\rm ext}$ using Galois co-objects. So far, the higher dimensional $q$-deformations of~$SU(n,1)$ remain undef\/ined on the von Neumann algebraic level.

\vspace{-2mm}

\subsection*{Structure of the paper}

In Section~\ref{SectHomogeneousSpaces} we study the homogeneous space of left and right invariant elements. We also give their dual spaces. Main goal here is the introduction of the von Neumann algebras $N$ and $\hat{N}$ as well as the C$^\ast$-algebras $\hat{N}_c$ and $\hat{N}_u$. These are homogeneous counterparts of the von Neumann algebra of a quantum group and its dual as well as the underlying reduced and universal dual C$^\ast$-algebra.

 In Section \ref{SectWeights}, we study the natural weights on these homogeneous spaces. Since the weight on the larger quantum group is generally not a state, the analysis is much more intricate compared to the compact Hopf-algebraic approach. We prove that the C$^\ast$-algebraic weights def\/ined here lift to the von Neumann algebraic (dual) weight. This result is a major ingredient for the quantum Plancherel--Godement theorem, since it allows us to apply Desmedt's auxiliary theorem \cite[Theorem~3.4.5]{Des}.
 Next, we introduce the necessary terminology of corepresentations that admit a vector that is invariant under the action of a subgroup. This is worked out in Section \ref{SectSphericalCoreps}.

In Section \ref{SectCorrespondence} we elaborate on a spherical version of Kustermans' result \cite{Kus}: there is a 1-1 correspondence between representations of the universal dual of a quantum group and corepresentations of the quantum group itself.  This will form the essential bridge between  \cite[Theorem~3.4.5]{Des} and the quantum Plancherel--Godement  Theorem \ref{ThmPG}.
Eventually, Section~\ref{SectPlancherelGodement} combines the results of Sections~\ref{SectHomogeneousSpaces}--\ref{SectCorrespondence} to prove a quantum version of the Plancherel--Godement theorem.

In Section~\ref{SectExample} we work out the example of $SU_q(1,1)_{\rm ext}$  together with its diagonal subgroup. We determine all the objects def\/ined in Sections~\ref{SectHomogeneousSpaces}--\ref{SectPlancherelGodement}. As an application of the theory we f\/ind product formulae for little $q$-Jacobi-functions that appear as matrix coef\/f\/icients of irreducible corepresentations.

\vspace{-1mm}

\section{Preliminaries and notation}\label{SectPreliminaries}

We brief\/ly recall the def\/inition and essential results from the theory of locally compact quantum groups. The results can be found in \cite{KusVae, KusVaeII} and \cite{Kus}. For an introduction we refer to~\cite{Tim} and~\cite{KusLec}. For the theory of weights on von Neumann algebras we refer to \cite{TakII}.

We use the following notational conventions. If $\pi$ and $\rho$ are linear maps, we write $\pi\rho$ for the composition $\pi \circ \rho$. $\iota$ denotes the identity homomorphism. The symbol $\otimes$ will be used for either the tensor product of two elements, of linear maps, the von Neumann algebraic tensor product or the tensor product of representations. We use the leg-numbering notation for operators. For example, if $W \in B(\cHil) \otimes B(\cHil)$, we write $W_{23}$ for $1 \otimes W$ and $W_{13} = (\Sigma \otimes 1) W_{23}(\Sigma \otimes 1)$, where $\Sigma: \cHil \otimes \cHil \rightarrow \cHil \otimes \cHil$ is the f\/lip.  For a linear map $A$, we denote $\Dom(A)$ for its domain.

Let $B$ be a Banach $\ast$-algebra. With a representation of $B$, we mean a $\ast$-homomorphism from $B$ to the bounded operators on a Hilbert space, which is referred to as the representation space. If $B$ is a C$^\ast$-algebra, we write $\Rep(B)$ for the equivalence classes of representations of $B$ and $\IR(B)$ for the equivalence classes of irreducible representations of $B$. With equivalence, we mean unitary equivalence. With slight abuse of notation we sometimes write $\pi \in \Rep(B)$ (or $\pi \in \IR(B)$) to mean that $\pi$ is an (irreducible) representation of $B$ instead of looking at its class.

 If $M$ is a von Neumann algebra and $\omega \in M_\ast$, we denote $\bar{\omega}$ for the normal functional def\/ined by $\bar{\omega}(x) = \overline{\omega(x^\ast)}$. For $\omega \in M_\ast$ and $x,y \in M$, we denote $x\omega$, $\omega y$, $x \omega y$ for the functionals def\/ined by $(x\omega)(z) = \omega(zx)$, $(\omega y)(z) = \omega(yz)$, $(x\omega y)(z) = \omega(yzx)$, where $z \in M$.

 Let $\phi$ be a weight on  $M$. Let $\mathfrak{n}_\phi = \{ x
\in M \mid \phi(x^\ast x) < \infty \}$, $\mathfrak{m}_\phi =
\mathfrak{n}_\phi^\ast \mathfrak{n}_\phi$. Let $\sigma^\phi$ be the modular automorphism group of $\phi$. We denote $\mathcal{T}_\phi$ for the Tomita
algebra def\/ined by
\[
 \mathcal{T}_\phi = \left\{ x \in M \mid x \textrm{ is analytic w.r.t. } \sigma^\phi \textrm{ and } \forall \,z \in \mathbb{C}: \sigma_z^\phi(x) \in \mathfrak{n}_\phi \cap \mathfrak{n}_\phi^\ast \right\}.
\]
For $x,y \in \mathcal{T}_\phi$, we write $x \phi y$ for the normal functional determined by $(x \phi y)(z) = \phi(yzx)$, $z \in M$.

\subsection*{Quantum groups}

 We use the Kustermans--Vaes def\/inition of a locally compact quantum group \cite{KusVae,KusVaeII}, see also \cite{KusLec,Tim,DaeleLcqg}.
\begin{dfn}
A locally compact quantum group $(M, \Delta)$ consists of the
following data:
\begin{enumerate}\itemsep=0pt
 \item A von Neumann algebra $M$;
\item A unital, normal $\ast$-homomorphism $\Delta: M \rightarrow M \otimes M$ satisfying the coassociativity relation $(\Delta \otimes \iota) \Delta = (\iota \otimes \Delta) \Delta$;
\item  Two  normal, semi-f\/inite, faithful weights
$\varphi$, $\psi$ on $M$ so that
\begin{alignat*}{4}
& \varphi \left( (\omega \otimes \iota )\Delta(x)\right)=
\varphi(x)\omega(1), \qquad && \forall \, \omega \in M^+_*,\  \forall\,
x\in \mathfrak{m}^+_\varphi
\qquad && \text{(left invariance);}& \\
& \psi\left( (\iota \otimes\omega)\Delta(x)\right)=
\psi(x)\omega(1), \qquad && \forall \, \omega \in M^+_*,\ \forall\,
x\in \mathfrak{m}^+_\psi \qquad && \text{(right invariance)}. &
\end{alignat*}
$\varphi$ is the left Haar weight and $\psi$ the right Haar
weight.
\end{enumerate}
\end{dfn}
\noindent Note that we suppress the Haar weights in the notation. A locally compact quantum group $(M, \Delta)$ is called compact if $\varphi$ and $\psi$ are states. $(M, \Delta)$ is called unimodular if $\varphi = \psi$. Compact quantum groups are unimodular.

 The triple $(\cHil, \pi,
\Lambda)$ denotes the GNS-construction with respect to the left
Haar weight $\varphi$. We may assume that $M$ acts on the
GNS-space $\cHil$. We use $J$ and $\nabla$ for the modular conjugation and modular operator of $\varphi$ and $\sigma$ for the modular automorphism group of $\varphi$. Recall that there is a constant $\nu \in \mathbb{R}^+$ called the scaling constant such that $\psi \sigma_t = \nu^{-t} \psi$. By applying \cite{VaeRad}, we see that there is a positive, self-adjoint operator $\delta$, called the modular element, such that $\psi = \varphi_\delta$, i.e.\ formally $\psi(\cdot) = \varphi(\delta^{\frac{1}{2}} \cdot \delta^{\frac{1}{2}} )$. For compact quantum groups, the scaling constant and the modular element are trivial.

\subsection*{Multiplicative unitary}
There exists a unique unitary operator $W\in
B(\cHil \otimes \cHil)$ def\/ined by
\[
W^\ast \left( \Lambda (a)\otimes  \Lambda(b) \right) = \left(
\Lambda \otimes \Lambda \right) \left( \Delta (b)(a\otimes
1)\right), \qquad a,b \in \nphi.
\]
$W$ is known as the multiplicative unitary. It satisf\/ies the
pentagon equation $W_{12}W_{13}W_{23}=W_{23}W_{12}$ in $B(\cHil
\otimes \cHil \otimes \cHil)$. Here we use the usual leg numbering notation. Moreover, $W$~implements the comultiplication, i.e.\ $\Delta(x) = W^\ast (1 \otimes x) W$ and $W \in M \otimes B(\cHil)$.

\subsection*{The unbounded antipode}

To $(M, \Delta)$ one can associate an unbounded map called the antipode $S: \Dom(S) \subseteq M \rightarrow M$. It can be def\/ined as the $\sigma$-strong-$\ast$ closure of the map $(\iota \otimes \omega)(W) \mapsto (\iota \otimes \omega)(W^\ast)$, where $\omega \in B(\cHil)_\ast$.
 One can prove that there exists a unique $\ast$-anti-automorphism $R: M \rightarrow M$ and a unique strongly continuous one-parameter group of $\ast$-automorphisms $\tau: \mathbb{R} \rightarrow {\rm Aut}(M)$ such that
\[
S = R \tau_{-i/2}, \qquad R^2 = \iota, \qquad \tau_t R = R \tau_t, \qquad \varphi \tau_t = \nu^{-t} \varphi, \qquad \forall\, t \in \mathbb{R}.
\]
$R$ is called the unitary antipode and $\tau$ is called the scaling group. Moreover,
\begin{gather}\label{EqnRCom}
\Delta R = \Sigma_{M,M} (R \otimes R) \Delta, \qquad \psi = \varphi R,
\end{gather}
where $\Sigma_{M,M}: M \otimes M \rightarrow M \otimes M$ is the f\/lip. Using the relative invariance property of the left Haar weight with respect to the scaling group, we def\/ine $P$ to be the positive operator on~$\cHil$ such that $P^{it} \Lambda(x) = \nu^{\frac{t}{2}} \Lambda(\tau_t(x))$, $t \in \mathbb{R}$, $x \in \nphi$.
 We use the notation
\[
    M_\ast^{\sharp} = \big\{ \omega \in M_\ast \mid \exists\, \theta \in M_\ast: (\theta \otimes \iota)(W) = (\omega \otimes \iota)(W)^\ast \big\}.
\]
For $\omega \in M_\ast^\sharp$, $\omega^\ast$ is def\/ined by $(\omega^\ast \otimes \iota)(W) = (\omega \otimes \iota)(W)^\ast$. In that case $\omega^\ast(x) = \overline{\omega}(S(x))$, $x \in \Dom(S)$. For $\omega \in M_\ast^\sharp$, we set $\Vert \omega \Vert_\ast = \max \{\Vert \omega \Vert, \Vert \omega^\ast \Vert \}$. $M_\ast^\sharp$ becomes a Banach-$\ast$-algebra with this norm.

\subsection*{The dual quantum group}
In \cite{KusVae, KusVaeII}, it is
proved that there exists a dual locally compact quantum group
$(\hat{M},\hat{\Delta})$, so that
$(\hat{\hat{M}},\hat{\hat{\Delta}}) = (M,\Delta)$. The dual left
and right Haar weight are denoted by $\hat{\varphi}$ and
$\hat{\psi}$. Similarly, all other dual objects will be denoted by
a hat, i.e.\ $\hat{\nabla}, \hat{J}, \hat{\delta}, \hat{\sigma}_t, \hat{W}, \ldots$.  $P$ is self-dual, i.e.\ $\hat{P} = P$. By construction,
\[\hat{M} = \overline{\big\{ (\omega
\otimes \iota)(W) \mid \omega \in M_\ast
\big\}}^{\sigma\text{-strong-}\ast}.
\] Furthermore, $\hat{W} =
\Sigma W^\ast \Sigma$, where $\Sigma$ denotes the f\/lip on $\cHil
\otimes \cHil$. This implies that $W \in M \otimes \hat{M}$ and
\[
M = \overline{\big\{ (\iota \otimes \omega)(W) \mid \omega \in
\hat{M}_\ast \big\}}^{\sigma\text{-strong-}\ast}.
\]
 For $\omega \in
M_\ast$, we use the standard notation $\lambda(\omega) = (\omega
\otimes \iota) (W) \in \hat{M}$. Then $\lambda: M_\ast^\sharp \rightarrow \hat{M}$ is a~representation.

 We denote $\mathcal{I}$ for the set of
$\omega \in M_\ast$, such that $\Lambda(x) \mapsto \omega(x^\ast)$,
$x \in \nphi$, extends to a bounded functional on $\cHil$. By the
Riesz theorem, for every $\omega \in \mathcal{I}$, there is a
unique vector denoted by $\xi(\omega) \in \cHil$ such that
$\omega(x^\ast) = \langle \xi(\omega)  , \Lambda(x) \rangle$, $x \in
\nphi$. The dual left Haar weight $\hat{\varphi}$ is def\/ined to be
the unique normal, semi-f\/inite, faithful weight on $\hat{M}$, with
GNS-construction $(\cHil, \iota, \hat{\Lambda})$ such that
$\lambda(\mathcal{I})$ is a $\sigma$-strong-$\ast$/norm core for
$\hat{\Lambda}$ and $\hat{\Lambda}(\lambda(\omega)) = \xi(\omega)$,
$\omega \in \mathcal{I}$. By def\/inition, this means that $\{ (\lambda(\omega), \xi(\omega)) \mid \omega \in \mathcal{I} \}$ is dense in the graph of $\hat{\Lambda}$ with respect to the product of the $\sigma$-strong-$\ast$ topology on $M$ and the norm topology on $\cHil$.

\subsection*{Corepresentations}

A (unitary) corepresentation is a unitary operator $U \in M \otimes B(\cHil_U)$ that satisf\/ies the relation $(\Delta \otimes \iota) (U) = U_{13}U_{23}$. In this paper all corepresentations are assumed to be unitary. It follows from the pentagon equation that the multiplicative unitary $W$ is a corepresentation on the GNS-space $\cHil$. Two corepresentations $U_1$, $U_2$ are equivalent if there is a unitary $T: \cHil_{U_1} \rightarrow \cHil_{U_2}$ such that $T (\omega \otimes \iota)(U_1) v = (\omega \otimes \iota)(U_2) T v$ for every $\omega \in M_\ast$, $v \in \cHil_{U_1}$. We denote $\ICM$ for the set of equivalence classes of all unitary corepresentations. For any corepresentation $U \in M \otimes B(\cHil_U)$ and $\omega \in B(\cHil_U)_\ast$, $(\iota \otimes \omega)(U) \in \Dom(S)$ and
\[
S (\iota \otimes \omega)(U) = (\iota \otimes \omega)(U^\ast) .
\]

\subsection*{Reduced quantum groups}
To every locally compact quantum group $(M, \Delta)$ one can associate a reduced C$^\ast$-algebraic quantum group. We def\/ine $M_c$ to be the norm closure of $\{ (\iota \otimes \omega)(W) \mid \omega \in \hat{M}_\ast \}$, which is a C$^\ast$-algebra. We restrict $\Delta$, $\varphi$ and $\psi$ to $M_c$ and denote the respective restrictions by $\Delta_c$, $\varphi_c$, $\psi_c$. In fact, $\Delta_c$ should be considered as a map into the multiplier algebra of the minimal tensor product of $M_c$ with itself. The GNS-construction representation $(\cHil, \Lambda, \pi)$ then restricts to a~GNS-representation of the C$^\ast$-algebraic weight $\varphi_c$, which is denoted by $(\cHil, \Lambda_c, \pi_c)$.
It is proven in \cite{KusVae} that $(M_c, \Delta_c)$ forms a reduced C$^\ast$-algebraic quantum group.
Similarly, one def\/ines the reduced dual  C$^\ast$-algebraic quantum group $(\hat{M}_c, \hat{\Delta}_c)$. The associated objects are denoted with a~hat.

\subsection*{Universal quantum groups}
Universal quantum groups were introduced by Kustermans \cite{Kus}.
 For $\omega \in M_\ast^\sharp$, we def\/ine
\[
    \Vert \omega \Vert_u = \sup \{ \Vert \pi(\omega) \Vert \mid \pi \textrm{ a }   \textrm{representation of } M_\ast^\sharp   \}.
\]
Recall that with a representation, we mean a $\ast$-homomorphism to the bounded operators on a~Hilbert space.
Note that this def\/ines a norm since the representation $\lambda$ is injective.
 Let $\hat{M}_u$ be the completion of $M_\ast^\sharp$ with respect to $\Vert \cdot \Vert_u$. We let $\lambda_u: M_\ast^\sharp \rightarrow \hat{M}_u$ denote the canonical embedding. $\hat{M}_u$ carries the following universal property: if $\pi$ is a~representation of $M_\ast^\sharp$, then there is a~unique representation $\rho$ of $\hat{M}_u$ such that $\pi = \rho \lambda_u$. In particular, from the representation $\lambda$ we get a surjective map $\hat{\vartheta}: \hat{M}_u \rightarrow \hat{M}_c$. We def\/ine a universal weight on $\hat{M}_u$ by setting $\hat{\varphi}_u = \hat{\varphi}_c \hat{\vartheta}$, and $\hat{\psi}_u = \hat{\psi}_c \hat{\vartheta}$. The GNS-representation of $\hat{\varphi}_u$ is given by $(\cHil, \hat{\Lambda}_u = \hat{\Lambda}_c \hat{\vartheta}, \hat{\pi}_u = \hat{\pi}_c \hat{\vartheta})$. If $U \in M \otimes B(\cHil_U)$ is a corepresentation of $M$, then the map $M_\ast^\sharp \rightarrow B(\cHil_U): \omega \mapsto (\omega \otimes \iota)(U)$ determines a representation of $\hat{M}_u$, which we denote by $\pi_U$. In fact, it is shown in~\cite{Kus} that this establishes a~1-1 correspondence between corepresentations of $M$ and non-degenerate representations of~$\hat{M}_u$.

For completeness, we mention that $\hat{M}_u$ can be equipped with a comultiplication $\hat{\Delta}_u$, which is a map from $\hat{M}_u$ to the multiplier algebra of the minimal tensor product of $\hat{M}_u$ with itself, such that $(\hat{M}_u, \hat{\Delta}_u)$ is a universal C$^\ast$-algebraic quantum group in the sense of~\cite{Kus}.
Similarly, one def\/ines $M_u, \Delta_u, \vartheta, \varphi_u, \psi_u, \ldots$.

\section{Spherical Fourier transforms}\label{SectHomogeneousSpaces}

Let $(G, K)$ be a pair consisting of a locally compact group $G$ and a compact subgroup $K$. If $L^1(K \backslash G \slash K)$, the bi-$K$-invariant $L^1$-functions on $G$ equipped with the convolution product is a~commutative algebra, then $(G,K)$ is a called a Gelfand pair. Examples of Gelfand pairs can be found in \cite[Chapter~7]{Dij} and \cite{Far}.

A notion of Gelfand pairs for compact quantum groups was introduced by Koornwinder~\cite{Koo}. We brief\/ly recall the def\/inition. Consider two unital Hopf-algebras $H$, $H_1$. Denote $\Delta$ for the comultiplication of $H$, denote $\varphi_1$ for the Haar functional on $H_1$. Suppose that there exists a~surjective morphism $\pi: H \rightarrow H_1$, so that $H_1$ is identif\/ied as a quantum subgroup of $H$. Now consider the left and right coactions
$
\Delta^l = (\pi \otimes \iota) \Delta,  \Delta^r = (\iota \otimes \pi) \Delta.
$
 Def\/ine
\[
H_1 \backslash H   = \big\{ h \in H \mid \Delta^l(h) = 1 \otimes h\big\}, \qquad  H \slash H_1   = \big\{ h \in H \mid \Delta^r(h) = 1 \otimes h\big\}
\]
 and set $H_1 \backslash H  \slash H_1 = (H_1 \backslash H)  \cap (H  \slash H_1)$. Classically, $H_1 \backslash H  \slash H_1$ corresponds to the algebra of bi-$K$-invariant elements. Set
\[
\tilde{\Delta}  = (\iota \otimes \varphi_1 \pi \otimes \iota)(\iota \otimes \Delta) \Delta.
\]
Now, the following def\/inition characterizes a quantum Gelfand pair. In fact, there are more equivalent def\/initions. We state the one which is closest to the theory we develop in the present section.
\begin{dfn}
Let $(H, H_1)$ be as above. $(H, H_1)$ is called a Gelfand pair if $\tilde{\Delta}$ is cocommutative, i.e.\ $\tilde{\Delta} = \Sigma_{H, H} \tilde{\Delta}$. The pair is called a strict Gelfand pair if moreover $H_1 \backslash H  \slash H_1$ is commutative. Here $\Sigma_{H, H}$ denotes the f\/lip.
\end{dfn}

A pair of compact groups $(G, K)$ is a Gelfand pair if and only if the Hopf-algebra of matrix coef\/f\/icients of unitary f\/inite dimensional representations form a quantum Gelfand pair (which is automatically strict).  Many deformations of classical Gelfand pairs form   strict Gelfand pairs in the Hopf-algebraic setting. Examples  can be found in for instance \cite{Flo, KoeAsk, NouYamMim, VaiGel} and \cite{VaiHyp}.

The aim of this section is to give a general framework of Gelfand pairs in the locally compact quantum group setting as introduced by Kustermans and Vaes \cite{KusVae, KusVaeII}. This puts the earlier studies as mentioned in the introduction in a non-compact, von Neumann algebraic setting.

One of the main motivations for the operator algebraic approach is that we can def\/ine sphe\-ri\-cal Fourier transforms. In particular, we show that  the structure   presented here allows us to prove a decomposition theorem analogous to the classical Plancherel--Godement theorem  \cite[Th\'eor\`eme~IV.2]{Far} or \cite[Section~6]{Dij}. The proof is an application of Desmedt's auxiliary result \cite[Theorem~3.4.5]{Des}.

As explained in the introduction, we do not assume a natural quantum analogue of the classical commutativity assumption on the convolution algebra of bi-invariant functions. Instead, we assume unimodularity of the larger quantum group, which is classically a weaker assumption. This allows us to cover the example of $SU_q(1,1)_{\rm ext}$, see Section~\ref{SectExample}.

\begin{nt}
Throughout Sections \ref{SectHomogeneousSpaces}--\ref{SectPlancherelGodement}, we f\/ix a locally compact quantum group $(M, \Delta)$ together with a closed quantum subgroup $(M_1,\Delta_1)$ which we assume to be the compact. Recall \cite[Def\/inition~2.9]{VaeVaiLow} that this means that we have a surjective $\ast$-homomorphism $\pi: M_u \rightarrow (M_1)_u$ on the level of universal C$^\ast$-algebras and the induced dual $\ast$-homomorphism $\hat{\pi}: (\hat{M}_1)_u \rightarrow \hat{M}_u$ lifts to a map on the level of von Neumann algebras $\hat{\pi}: \hat{M}_1 \rightarrow \hat{M}$, which with slight abuse of notation is denoted by $\hat{\pi}$ again. When we encounter $\hat{\pi}$ in this paper, we always mean the von Neumann algebraic map.

Note $\pi$ and $\hat{\pi}$ are in principal also used for the GNS-representations of $M$ and $\hat{M}$. However, we omit the maps most of the time, since $M$ and $\hat{M}$ are identif\/ied with their GNS-representations. In that case we explicitly need the GNS-representations, we will mention this.

 We mention that from a certain point, see Notation~\ref{NtUnimodular}, we will assume that $(M, \Delta)$ is a~unimodular quantum group.

We use $\Sigma_{M_1, M}: M_1 \otimes M \rightarrow M \otimes M_1$ to denote the f\/lip. The objects associated with $(M_1, \Delta_1)$ will be equipped with a subscript 1, i.e.\ $S_1, R_1, \tau_1, \nu_1, \varphi_1, \ldots$.
\end{nt}

\subsection*{Homogeneous spaces}
Due to \cite[Proposition 3.1]{VaeVaiExt}, there are canonical right coactions of $(M_1, \Delta_1)$ on $M$, denoted by $\beta, \gamma: M \rightarrow M \otimes M_1$, which are normal $\ast$-homomorphisms uniquely determined by
\begin{gather}\label{EqnBeta}
 (\beta \otimes \iota)(W)  = W_{13} \left((\iota \otimes \hat{\pi})(W_1)\right)_{23}, \qquad
 (\gamma \otimes \iota)(W) =  \left((R_1 \otimes \hat{\pi})(W_1)\right)_{23} W_{13}.
\end{gather}
We have the relation $\gamma = (R \otimes \iota)  \beta  R$. The map $\beta$ corresponds classically to right translation, whereas $\gamma$ corresponds to left translation.

\begin{rmk}
 In \cite[Section 3]{VaeVaiExt}, the roles of $(M_1, \Delta_1)$ and $(\hat{M}_1, \hat{\Delta}_1)$ are interchanged. Using our conventions for the roles of $(M_1, \Delta_1)$ and $(\hat{M}_1, \hat{\Delta}_1)$, recall the left action $\mu$ of $(M_1, \Delta_1)$ on~$M$ and the left action $\theta$ of $(M_1, \Delta^{\rm op})$ on $M$ from \cite[Proposition~3.1]{VaeVaiExt}. By this proposition, $\beta$ equals the right coaction $\Sigma_{M_1, M}   \theta$. $\gamma$ equals the right coaction $\Sigma_{M_1, M} (R_1 \otimes \iota) \mu$.
\end{rmk}

\begin{lem}\label{LemInterchangeLR}
As maps $M \rightarrow M \otimes M \otimes M_1$, we have an equality
\begin{gather}
 (\iota \otimes \Sigma_{M_1, M}) (\beta \otimes \iota)  \Delta  = (\iota \otimes \iota \otimes R_1)(\iota \otimes \gamma) \Delta.\label{EqnGammaBetaRelI}
\end{gather}
\end{lem}
\begin{proof}
For the left hand side, using the pentagon equation and \cite[Proposition 3.1]{VaeVaiExt},
 \[
  (\beta \otimes \iota \otimes \iota)  (\Delta\otimes \iota) W = (\beta \otimes \iota \otimes \iota)  W_{13} W_{23} =
W_{14} ((\iota \otimes \hat{\pi})(W_1))_{24} W_{34}.
 \]
For the right hand side, using again \cite[Proposition 3.1]{VaeVaiExt},
\begin{gather*}
 (\iota \otimes \iota \otimes R_1 \otimes \iota)(\iota \otimes \gamma \otimes \iota) (\Delta \otimes \iota) W =
 (\iota \otimes \iota \otimes R_1 \otimes \iota)(\iota \otimes \gamma \otimes \iota) W_{13} W_{23} \\
\qquad {}=   (\iota \otimes \iota \otimes R_1 \otimes \iota) W_{14} ((R_1 \otimes \iota)(\iota \otimes \hat{\pi})(W_1))_{34} W_{24} =
W_{14} ((\iota \otimes \hat{\pi})(W_1))_{34} W_{24}.
\end{gather*}
The lemma follows by the fact that the elements $\{ (\iota \otimes \omega)(W) \mid \omega \in \hat{M}_\ast\}$ are $\sigma$-strong-$\ast$ dense in~$M$.
\end{proof}

\begin{dfn}
We denote $M^\beta$ for the f\/ixed point algebra $\left\{ x \in M \mid \beta(x) = x \otimes 1 \right\}$. Similarly, $M^\gamma$ denotes the f\/ixed point algebra of $\gamma$. By def\/inition of $\gamma$ we f\/ind $M^\gamma = R(M^\beta)$.  We def\/ine
\[
N = M^\beta \cap M^\gamma.
\]
\end{dfn}
Note that $M^\beta$, $M^\gamma$ and $N$ are von Neumann algebras. Furthermore, $\Delta( M^\beta) \subseteq M \otimes M^\beta$. Also, by $M^\gamma = R(M^\beta)$, (\ref{EqnRCom}) and~(\ref{EqnGammaBetaRelI}) it follows that $\Delta(M^\gamma) \subseteq M^\gamma \otimes M$.

We recall from \cite{VaeUni} that we have normal, faithful operator valued weights,
\[
 T_\beta: \  M^+ \rightarrow (M^{\beta})^{+}: \ x \mapsto (\iota \otimes \varphi_1) \beta(x);\qquad
 T_\gamma: \ M^+ \rightarrow (M^{\gamma})^{+}: \ x \mapsto (\iota \otimes \varphi_1) \gamma(x).
\]
Since $(M_1, \Delta_1)$ is compact, $T_\beta$ and $T_\gamma$ are f\/inite. We extend the domains of $T_\beta$ and $T_\gamma$ to $M$ in the usual way. We denote the extensions again by $T_\beta$ and $T_\gamma$. The composition of~$T_\beta$ and~$T_\gamma$ forms a well-def\/ined map on $M$. Note that $T_\beta(x^\ast) = T_\beta(x)^\ast$ and $T_\gamma(x^\ast) = T_\gamma(x)^\ast$, where $x \in M$.

\begin{rmk}
The spaces $M^\gamma$, $M^\beta$ where already introduced in \cite{VaeImp} as homogeneous spaces. They also fall within the def\/inition of a homogeneous space as introduced  by Kasprzak \cite[Remark~3.3]{Kas}. Moreover, we stretch that $T_\beta$ and $T_\gamma$ are conditional expectation values, which properties have been studied in the related papers~\cite{SalmiSkalski} and~\cite{Tomatsu}.
\end{rmk}

\begin{lem}\label{LemTProperties} $T_\gamma: M \rightarrow M^\gamma$ and $T_\beta: M \rightarrow M^\beta$ satisfy the following properties:
\begin{enumerate}\itemsep=0pt
\item\label{LemTI} $T_\beta T_\gamma = T_\gamma T_\beta $;
\item\label{LemTII} $\Delta  T_\beta = (\iota \otimes T_\beta) \Delta$ and  $\Delta T_\gamma = (T_\gamma \otimes \iota)  \Delta$;
\item\label{LemTIII} $(\iota \otimes T_\gamma)  \Delta = (T_\beta \otimes \iota)  \Delta$;
\item\label{LemTIV}  $T_\gamma S \subseteq S T_\beta$ and  $T_\beta S \subseteq S T_\gamma$.
\end{enumerate}
\end{lem}

\begin{proof} (\ref{LemTI}) This follows from the fact $(\iota \otimes \Sigma_{M_1, M_1} ) (\gamma \otimes \iota) \beta = (\beta \otimes \iota) \gamma$, which can be established as in the proof of Lemma \ref{LemInterchangeLR}.

(\ref{LemTII}) We prove that $\Delta T_\gamma = (T_\gamma \otimes \iota)  \Delta$, the other equation can be proved similarly using $\beta = (R \otimes \iota)\gamma R$ and (\ref{EqnRCom}). We f\/ind:
\begin{gather*}  (\Delta T_\gamma \otimes \iota)(R \otimes \iota)(W)  =
 (\iota \otimes \iota \otimes \varphi_1 \otimes \iota) (\Delta \otimes \iota \otimes \iota) (R \otimes \iota \otimes \iota)(\beta \otimes \iota)(W) \\
\phantom{(\Delta T_\gamma \otimes \iota)(R \otimes \iota)(W) } =
   (\Sigma_{M,M} \otimes \iota) (R \otimes R \otimes \iota)(\Delta \otimes \iota) (W (1\otimes \hat{\pi}((\varphi_1 \otimes \iota)(W_1)))) \\
\phantom{(\Delta T_\gamma \otimes \iota)(R \otimes \iota)(W) } =
(R \otimes R \otimes \iota) W_{23} W_{13}  (1\otimes 1\otimes \hat{\pi}((\varphi_1 \otimes \iota)(W_1)))  \\
\phantom{(\Delta T_\gamma \otimes \iota)(R \otimes \iota)(W) } =
 (\iota \otimes \varphi_1 \otimes \iota \otimes \iota)(R \otimes \iota \otimes \iota \otimes \iota) (\beta \otimes R \otimes \iota) W_{23}W_{13} \\
\phantom{(\Delta T_\gamma \otimes \iota)(R \otimes \iota)(W) } =
 (T_\gamma \otimes \iota \otimes \iota) (\Sigma_{M,M} \otimes \iota) (R \otimes R \otimes \iota)W_{13} W_{23} \\
\phantom{(\Delta T_\gamma \otimes \iota)(R \otimes \iota)(W) } =
 (T_\gamma \otimes \iota \otimes \iota)(\Delta \otimes \iota)(R \otimes \iota)(W).
\end{gather*}
Now, the equation follows by taking slices of the second leg of $W$.

(\ref{LemTIII}) Since $(M_1, \Delta_1)$ is compact, it is unimodular and hence $\varphi_1 = \varphi_1  R_1$. Using Lemma \ref{LemInterchangeLR} we f\/ind:
\begin{gather*}
 (\iota \otimes T_\gamma)   \Delta = (\iota \otimes \iota \otimes \varphi_1) (\iota \otimes \gamma) \Delta = (\iota \otimes \iota \otimes \varphi_1   R_1) (\iota \otimes \gamma) \Delta \\
 \phantom{(\iota \otimes T_\gamma)   \Delta}{} =
 (\iota \otimes \iota \otimes \varphi_1)\Sigma_{23} (\beta \otimes \iota) \Delta =  (T_\beta \otimes \iota)   \Delta.
\end{gather*}

(\ref{LemTIV})
It follows from \cite[Proposition 2.1, Corollary 2.2]{KusVaeII} that $(\tau_t \otimes \iota)(W) = (\iota \otimes \hat{\tau}_{-t})(W), t \in \mathbb{R}$. Therefore,
\[
(\beta \tau_t \otimes \iota)(W) = (\beta \otimes \hat{\tau}_{-t})(W) = (\iota \otimes \iota \otimes \hat{\tau}_{-t})(W_{13}((\iota \otimes \hat{\pi})(W_1))_{23}).
\]
By \cite[Proposition 5.45]{KusVae}, we know that $\hat{\pi}(\hat{\tau}_1)_t = \hat{\tau}_t \hat{\pi}$. Here $(\hat{\tau}_1)_t$ denotes the scaling group of $(\hat{M}_1, \hat{\Delta}_1)$. Continuing the equation, we f\/ind
\begin{gather*}
  (\beta \tau_t \otimes \iota)(W)  = (\iota \otimes \hat{\tau}_{-t})(W)_{13} ((\iota \otimes   \hat{\pi}(\hat{\tau}_1)_{-t})(W_1))_{23} \\
\phantom{(\beta \tau_t \otimes \iota)(W)}{}
=  (\tau_t \otimes \iota)(W)_{13} ((\tau_1)_{t}\otimes \hat{\pi})(W_1))_{23} =
(\tau_t \otimes (\tau_1)_{t} \otimes \iota) (\beta \otimes \iota)(W).
\end{gather*}
We have $\varphi_1 (\tau_1)_t =  \varphi_1$, since $(M_1, \Delta_1)$ is compact. Hence, $T_\beta \tau_t = \tau_t T_\beta$. With a similar computation, involving the relation $(R \otimes \hat{R})(W) = W^\ast$, see \cite[Proposition 2.1, Corollary 2.2]{KusVaeII}, we f\/ind $T_\beta R = R T_\gamma$. The proposition follows, since $S = R \tau_{-i/2}$.
\end{proof}

\begin{dfn}
For $x \in N$, we def\/ine
\begin{gather}\label{EqnShiftedCom}
\Delta^\natural(x) = (\iota \otimes T_\gamma) \Delta(x) = (T_\beta \otimes \iota) \Delta(x),
\end{gather}
 see Lemma \ref{LemTProperties} (\ref{LemTIII}). This is the von Neumann algebraic version of \cite[equation~(4)]{VaiHyp}.
\end{dfn}
  Recall that $\Delta(N) \subseteq M^\gamma \otimes M^\beta$, so that $\Delta^\natural(N) \subseteq N \otimes N$. Moreover, using (\ref{LemTI})--(\ref{LemTIII}) of the previous lemma, $\Delta^\natural$ is coassociative, i.e.
\begin{gather}
   (\iota \otimes \Delta^\natural) \Delta^\natural  = (\iota \otimes \iota \otimes T_\gamma)(\iota \otimes \Delta)(\iota \otimes T_\gamma)\Delta
 = (\iota \otimes T_\beta   T_\gamma \otimes \iota) (\iota \otimes \Delta) \Delta \nonumber\\
\phantom{(\iota \otimes \Delta^\natural) \Delta^\natural }{}
=   (\iota \otimes T_\gamma   T_\beta \otimes \iota) (\Delta \otimes \iota) \Delta
 =  ( \Delta^\natural \otimes \iota) \Delta^\natural.\label{EqnCoas}
\end{gather}
Note that $\Delta^\natural$ is unital, but generally not multiplicative.

\begin{dfn}
 We def\/ine a convolution product $\ast^\natural$ on $N_\ast$,
\[
\omega_1 \ast^\natural \omega_2 = (\omega_1 \otimes \omega_2) \Delta^\natural, \qquad  \omega_1, \omega_2 \in N_\ast.
\]
\end{dfn}

This convolution product is associative, since by (\ref{EqnCoas}) $\Delta^\natural$ is coassociative.

\begin{dfn}
If $\omega \in N_\ast$, then we def\/ine
\[
\tilde{\omega} = \omega T_\gamma  T_\beta = \omega T_\beta T_\gamma \in M_\ast.
\]
 For $\omega \in M_\ast$, we put $\tilde{\omega} = (\omega \vert_N)^\sim$.
\end{dfn}

\begin{rmk}
Note that $\Delta^\natural$ is the von Neumann algebraic version of $\tilde{\Delta}$ \cite{VaiHyp}, which was used to def\/ine hypergroup structures. See also the remarks at the beginning of this section. Here, we will not focus on hypergroups for two reasons. First of all, we stay mostly at the measurable von Neumann algebraic setting, which does not allow one to directly def\/ine the generalized shift operators~\cite{VaiHyp}. Moreover, we will not assume that $N$ is Abelian, i.e.\ what is called a strict Gelfand pair in \cite{VaiHyp}.
\end{rmk}

\begin{prop}\label{PropInvariance}
 The map $N_\ast \ni \omega \mapsto \tilde{\omega}$ defines a bijective, norm preserving correspondence between $N_\ast$ and the   functionals $\theta \in M_\ast$ that satisfy the invariance properties:
\begin{enumerate}\itemsep=0pt
 \item\label{LemInvarianceI} $(\theta \otimes \iota) \beta(x) = \theta(x) 1_{M_1}$ for all $x \in M$;
\item\label{LemInvarianceII} $(\theta \otimes \iota) \gamma(x) = \theta(x) 1_{M_1}$ for all $x \in M$.
\end{enumerate}
Moreover, for $\omega_1, \omega_2 \in N_\ast$ we have $ \left( \omega_1 \ast^\natural \omega_2 \right)^\sim = \tilde{\omega}_1 \ast \tilde{\omega}_2$.
\end{prop}
\begin{proof}
Using right invariance of $\varphi_1$, for $x \in M$,
\begin{gather*}
(T_\beta \otimes \iota) \beta(x) = (\iota \otimes \varphi_1 \otimes \iota)(\beta \otimes \iota) \beta(x) = (\iota \otimes \varphi_1 \otimes \iota)(\iota \otimes \Delta_1) \beta(x) \\
\phantom{(T_\beta \otimes \iota) \beta(x)}{}
= (\iota \otimes \varphi_1) \beta(x) \otimes 1_{M_1} = T_\beta(x)\otimes 1_{M_1}.
 \end{gather*}
Similarly, the right invariance of $\varphi_1$ implies $(T_\gamma \otimes \iota) \gamma(x) = T_\gamma(x)\otimes 1_{M_1}$, $x\in M$. Using (\ref{LemTI}) of Lemma \ref{LemTProperties} one easily verif\/ies that for $\omega \in N_\ast$, $\tilde{\omega}$ satisf\/ies the invariance properties~(\ref{LemInvarianceI}) and~(\ref{LemInvarianceII}).

Let $\omega \in N_\ast$. For $x \in N$ we have $\omega(x) = \tilde{\omega}(x)$, so that $N_\ast \ni \omega \mapsto \tilde{\omega}$ is injective. If $\theta \in M_\ast$ satisf\/ies (\ref{LemInvarianceI}), then for $x \in M$,
\[
 \theta   T_\beta (x) = \theta (\iota \otimes \varphi_1) \beta(x) =  \varphi_1 (\theta \otimes \iota) \beta(x) = \theta(x).
\]
Similarly, if $\theta \in M_\ast$ satisf\/ies~(\ref{LemInvarianceII}), then $\theta   T_\gamma(x) = \theta(x)$. We f\/ind that $\theta = ( \theta \vert_N )^\sim$ if $\theta$ satisf\/ies (\ref{LemInvarianceI}) and (\ref{LemInvarianceII}). So  $N_\ast \ni \omega \mapsto \tilde{\omega}$ ranges over the normal functionals on $M$ that satisfy the invariance properties (\ref{LemInvarianceI}) and (\ref{LemInvarianceII}).

Using the left invariance of $\varphi_1$, it is a straightforward check that $T_\gamma T_\gamma = T_\gamma$. Then, using (\ref{LemTI})--(\ref{LemTIII}) of Lemma \ref{LemTProperties}, we f\/ind:
\begin{gather*}
(\omega_1 \ast^\natural \omega_2)^\sim  = (\omega_1 \otimes \omega_2) (\iota \otimes T_\gamma) \Delta T_\gamma T_\beta =
(\omega_1 \otimes \omega_2) (T_\gamma \otimes T_\gamma T_\beta) \Delta \\
\phantom{(\omega_1 \ast^\natural \omega_2)^\sim}{}
=
(\omega_1 \otimes \omega_2) (T_\gamma \otimes T_\gamma^2 T_\beta) \Delta =
(\omega_1 \otimes \omega_2) (T_\gamma T_\beta \otimes T_\gamma T_\beta) \Delta
= \tilde{\omega_1} \ast \tilde{\omega_2}.\tag*{\qed}
\end{gather*}
\renewcommand{\qed}{}
\end{proof}

\begin{prop}\label{PropSharpTilde}
 For $\omega \in M_\ast^\sharp$, we have $\tilde{\omega} \in M_\ast^\sharp$ and $(\tilde{\omega})^\ast = (\omega^\ast)^{\sim}$.
\end{prop}
\begin{proof}
Using (\ref{LemTIV}) of Lemma \ref{LemTProperties}, for $x \in \Dom(S)$, we f\/ind
\[
 \overline{\tilde{\omega}}(S(x)) = \overline{\omega(T_\beta T_\gamma (S(x)^\ast))} = \overline{\omega((T_\beta T_\gamma (S(x)))^\ast)} = \overline{ \omega(S(T_\beta T_\gamma(x))^\ast)} = (\omega^\ast)^\sim(x).
\]
So $(\omega^\ast)^\sim \in M_\ast$ has the property $((\omega^\ast)^\sim \otimes \iota)(W) = (\tilde{\omega}  \otimes \iota)(W)^\ast$. This proves that $\tilde{\omega} \in M_\ast^\sharp$ and $\tilde{\omega}^\ast = (\omega^\ast)^\sim$.
\end{proof}

The following proposition is proved in \cite[Proposition 3.1]{VaeVaiExt}.

\begin{prop}\label{PropTGamma}
If $x \in \nphi$, then $T_\gamma (x) \in \nphi$ and $\Lambda(T_\gamma(x)) = \hat{\pi}((\varphi_1 \otimes \iota)(W_1)) \Lambda(x)$.
\end{prop}

Proposition \ref{PropTGamma} def\/ines an orthogonal projection $\hat{\pi}((\varphi_1 \otimes \iota)(W_1))$ for which we simply write
\[
P_\gamma = \hat{\pi}((\varphi_1 \otimes \iota)(W_1)).
\]
 Classically, it corresponds to projecting onto the space of functions that are left invariant with respect to the compact subgroup. Note that $P_\gamma \in \hat{M}$ and
\[
N = T_\beta T_\gamma(M) = \left\{ (\iota \otimes \omega_{P_\gamma v, P_\gamma w})(W) \mid v,w \in \cHil \right\}.
\]

We need a similar result as Proposition \ref{PropTGamma} for $T_\beta$. For this we need unimodularity of $(M, \Delta)$. Classically, if $G$ is a group with compact subgroup $K$ such that $(G, K)$ forms a Gelfand pair, one can prove that $G$ is unimodular, see \cite[Proposition~I.1]{Far}. The natural def\/inition of a quantum Gelfand pair would be to require that $\Delta^\natural$ is cocommutative. However, we like to stretch the def\/inition of a Gelfand pair a bit to handle the example of $SU_q(1,1)_{{\rm ext}}$. The following essential result, see Proposition \ref{PropTBeta}, is the motivation of assuming the (classically) weaker condition that $(M, \Delta)$ is unimodular, see Notation \ref{NtUnimodular}.

First, we need the following lemma. Note that for $a, b \in \mathcal{T}_\varphi$, we have $a \varphi b \in M_\ast$ and hence $\mathcal{T}_{\varphi} \varphi \mathcal{T}_{\varphi}$ is a subset of $M_\ast$. Recall that for $\omega \in \mathcal{I}$, $\xi(\omega) \in \cHil$ is def\/ined using the Riesz theorem as the unique vector such that $\langle  \xi(\omega) , \Lambda(x)\rangle = \omega(x^\ast)$, $x \in \nphi$. By \cite[Lemma~8.5]{KusVae}, $\mathcal{T}_\varphi \varphi \mathcal{T}_\varphi$ is included in~$\mathcal{I}$.

\begin{lem}\label{LemTechLem}
Let $\omega \in \mathcal{I}$. There exists a net $(\omega_j)_j$ in $\mathcal{T}_{\varphi} \varphi \mathcal{T}_{\varphi}$ such that $\omega_j \rightarrow \omega$ in norm and $\xi(\omega_j) \rightarrow \xi(\omega)$ in norm.
\end{lem}

\begin{proof}
 For $\omega \in \mathcal{I}$ we def\/ine the norm $\Vert \omega \Vert_{\mathcal{I}} = \max \{ \Vert \omega \Vert, \Vert \xi(\omega) \Vert \}$. We have to prove that $\mathcal{T}_{\varphi} \varphi \mathcal{T}_{\varphi}$ is dense in $\mathcal{I}$ with respect to this norm. This is exactly what is obtained in the proof of  \cite[Proposition 3.4]{Cas}. Indeed, let $ L$ and $k$ be as in \cite{Cas}. As indicated in the introduction of \cite{Cas},  $\mathcal{T}_{\varphi}^2 \subseteq L$ and for $a,b \in \mathcal{T}_{\varphi}$, $k(ab) = \sigma_i(b) \varphi a$, see also \cite[Corollary 2.15]{Cas}. So $k(\mathcal{T}_{\varphi}^2 ) = \mathcal{T}_{\varphi} \varphi \mathcal{T}_{\varphi}$. The proposition yields that $k(\mathcal{T}_{\varphi}^2 )$ is dense in $\mathcal{I}$.
\end{proof}

\begin{prop}\label{PropTBeta}
Suppose that $(M, \Delta)$ is unimodular. For $x \in \nphi$, we have $T_\beta(x) \in \nphi$. The map $\Lambda(x) \mapsto \Lambda(T_\beta(x))$ is bounded and it extends continuously to the projection $\hat{J} P_\gamma \hat{J}$.
\end{prop}
\begin{proof}
By Proposition \ref{PropTGamma}, we see that $\Lambda(T_\gamma(x^\ast)) = \Lambda(T_\gamma(x)^\ast)$, $x \in \nphi \cap \nphi^\ast$. Denote $T$ for the closure of the map $\Lambda(x) \mapsto \Lambda(x^\ast)$, $x \in \nphi \cap \nphi^\ast$. We see that $P_\gamma \cHil$ is an invariant subspace for~$T$. Since $T = J \nabla^{1/2}$, we f\/ind that $\nabla^{it}$, $t \in \mathbb{R}$ commutes with $P_\gamma$.

Recall \cite[Lemma 8.8, Proposition 8.9]{KusVae} that $\hat{\nabla}^{it} =  P^{it} J \delta^{it} J$ and by Pontrjagin duality  $\nabla^{it} =  \hat{P}^{it} \hat{J} \hat{\delta}^{it} \hat{J}$. Using $\delta = 1$ and the self-duality $\hat{P} = P$, we see that $\hat{\nabla}^{it} =  \nabla^{it} \hat{J} \hat{\delta}^{-it} \hat{J}$. Since $\hat{\delta}$ is af\/f\/iliated with $\hat{M}$ and using the previous paragraph, this shows that $P_\gamma \hat{\nabla}^{it} = \hat{\nabla}^{it} P_\gamma$. Hence $\hat{\sigma}_t (P_\gamma) = P_\gamma$.

Now, let $a, b \in \mathcal{T}_{\hat{\varphi}}$ and put $\omega = a\hat{ \varphi} b \in \hat{M}_\ast$. Then,
\[
T_\beta\left( (\iota \otimes \omega)(W^\ast) \right) = (\iota \otimes \omega)\left( (1 \otimes P_\gamma) W^\ast \right).
\]
Since $P_\gamma \in \hat{M}$ is invariant under $\hat{\sigma}_t$, it follows from \cite[Chapter VIII, Lemma 2.4 (ii) and Lemma 2.5]{TakII} that $\omega P_\gamma \in \hat{\mathcal{I}}$ and
\[
\hat{\xi}(\omega P_\gamma ) = \hat{\Lambda}(a \hat{\sigma}_{-i}(b) \hat{\sigma}_{-i}(P_\gamma))  = \hat{J} P_\gamma \hat{J} \hat{\Lambda}(a \hat{\sigma}_{-i}(b) ).
\]

Now, let $\omega \in \hat{\mathcal{I}}$. We prove the proposition for $x = \hat{\lambda}(\omega) \in \nphi$. Let $(\omega_j)_{j \in J}$ be a net in $\mathcal{T}_{\hat{\varphi}} \hat{\varphi} \mathcal{T}_{\hat{\varphi}}$ that converges to $\omega$ and such that $\hat{\xi}(\omega_j)$ converges to $\hat{\xi}(\omega)$, c.f.\ Lemma \ref{LemTechLem}. Then, $T_\beta(\hat{\lambda}(\omega_j)) \in \nphi$ and  $T_\beta(\hat{\lambda}(\omega_j)) \rightarrow T_\beta(x)$ in the $\sigma$-weak topology. Furthermore, $\Lambda(T_\beta(\hat{\lambda}(\omega_j))) = \hat{J} P_\gamma \hat{J}\hat{\xi}(\omega_j)$ is norm convergent. Since $\Lambda$ is $\sigma$-weak/weak closed, and $\Dom(\Lambda) = \nphi$, this proves that $T_\beta(x) \in \nphi$ and $\Lambda(T_\beta(x)) =   \hat{J} P_\gamma \hat{J} \Lambda(x)$.

 Since the elements in $\hat{\lambda}(\hat{\mathcal{I}})$ form a $\sigma$-strong-$\ast$/norm core for $\Lambda$, this proves the proposition.
\end{proof}

We will write $P_\beta$ for the projection $\hat{J} P_\gamma \hat{J}$. In particular $P_\beta \in \hat{M}'$. Under the assumption that $(M, \Delta)$ is unimodular, we see that $P_\beta$ projects onto the elements that are right invariant with respect to the closed quantum subgroup $(M_1 ,\Delta_1)$. Since we will need this interpretation of $P_\beta$, i.e.\ Proposition \ref{PropTBeta}, we assume unimodularity from now on.

\begin{nt}\label{NtUnimodular}
From this point we assume that $(M, \Delta)$ is unimodular, i.e.\ $\varphi = \psi$.
\end{nt}

We are ready to def\/ine the dual structures associated with $N$. We def\/ine left and right invariant analogues of the dual von Neumann algebraic quantum group and the universal dual C$^\ast$-algebraic quantum group. These duals can be constructed by means of the multiplicative unitary $W$ associated with $(M, \Delta)$. We def\/ine
\begin{gather*}
N_\ast^\sharp = \left\{\omega \in N_\ast \mid \exists\, \theta \in N_\ast: (\tilde{\omega} \otimes \iota)(W)^\ast = (\tilde{\theta} \otimes \iota)(W) \right\}.
\end{gather*}
For $\omega \in N_\ast^\sharp$, we set $\Vert \omega \Vert_\ast = \max \{ \Vert \omega \Vert, \Vert \omega^\ast \Vert \}$. Then, $N_\ast^\sharp$ becomes a Banach-$\ast$-algebra with respect to this norm. Proposition \ref{PropSharpTilde} shows that $\omega \in N_\ast^\sharp$ if and only if $\tilde{\omega} \in M_\ast^\sharp$. Note that  $N_\ast^\sharp$ is dense in $N_\ast$. Indeed, the restriction map $M_\ast \rightarrow N_\ast: \omega \mapsto \omega\vert_N$  is continuous and the image of the subset $M_\ast^\sharp \subseteq M_\ast$ is contained in $N_\ast^\sharp$ by Proposition \ref{PropSharpTilde}. The inclusion $M_\ast^\sharp  \subseteq M_\ast$ is dense, see  \cite[Lemma 2.5]{KusVaeII}. Hence $N_\ast^\sharp$ is dense in $N_\ast$. Using this together with Propositions \ref{PropInvariance} and~\ref{PropSharpTilde}, we see that
\[
\hat{N} = \overline{\left\{ (\tilde{\omega} \otimes \iota)(W) \mid \omega \in N_\ast \right\}}^{\sigma\text{-strong-}\ast},
\]
is a $\ast$-subalgebra of $\hat{M}$. Since we can conveniently write
$
(\tilde{\omega} \otimes \iota)(W) = P_\gamma (\omega \otimes \iota)(W) P_\gamma
$
by (\ref{EqnBeta}),
 we see that $\hat{N}$ is a von Neumann algebra if considered as acting on $P_\gamma \cHil$, so that $P_\gamma$ is its unit. In particular $\hat{N} = P_\gamma \hat{M} P_\gamma$.

We def\/ine $\hat{N}_c$ to be the norm closure of the set $\left\{ (\tilde{\omega} \otimes \iota)(W) \mid \omega \in N_\ast \right\}$. Then $\hat{N}_c$ is a C$^\ast$-subalgebra of the reduced dual C$^\ast$-algebra $\hat{M}_c$.

 For $\omega \in N_\ast^\sharp$, we def\/ine
\[
    \Vert \omega \Vert_u^\natural = \sup \big\{ \Vert \pi(\omega) \Vert \mid \pi \textrm{ a representation of } N_\ast^\sharp  \big\}.
\]
Note that this def\/ines a norm since the representation $\omega \mapsto (\tilde{\omega} \otimes \iota)(W)$ is injective as follows using the bijective correspondence established in Proposition~\ref{PropInvariance}.
 Let $\hat{N}_u$ be the completion of~$N_\ast^\sharp$ with respect to~$\Vert \cdot \Vert_u^\natural$.

Recall the map $\lambda: M_\ast^\sharp \rightarrow \hat{M}: \omega \mapsto (\omega \otimes \iota)(W)$. We set
\[
\lambda^\natural: \ N_\ast^\sharp \rightarrow \hat{N}: \ \omega \mapsto (\tilde{\omega} \otimes \iota)(W).
\]
  Note that the image of $\lambda^\natural$ is contained in $\hat{N}_c$ and we will use this implicitly. $\lambda_u: M_\ast^\sharp \rightarrow \hat{M}_u$ is the canonical inclusion and similarly
\[
\lambda_u^\natural: \ N_\ast^\sharp \rightarrow \hat{N}_u
\]
 denotes the canonical inclusion.

Recall that $\hat{N}_u$ is a C$^\ast$-algebra with the following universal property: if $\pi$ is a representation of $N_\ast^\sharp$ on a Hilbert space, then there is a unique representation $\rho$ of $\hat{N}_u$ such that $\pi = \rho \lambda_u^\natural$.  By this universal property, the map $N_\ast^\sharp \rightarrow \hat{M}_u: \omega \mapsto \lambda_u(\tilde{\omega}) $ extends to a representation
\[
\iota_u: \ \hat{N}_u \rightarrow \hat{M}_u.
\]
  Similarly, the map $\lambda^\natural: N_\ast^\sharp \rightarrow \hat{N}_c$ gives rise to a surjective map
 \[
\hat{\vartheta}^\natural: \ \hat{N}_u \rightarrow \hat{N}_c.
\]
 In particular, $\hat{\vartheta} \iota_u = \hat{\vartheta}^\natural$, where $\hat{\vartheta}: \hat{M}_u \rightarrow \hat{M}_c$ was the canonical projection induced by the representation $\lambda: M_\ast^\sharp \rightarrow \hat{M}_c$.

\begin{rmk}
Note that we do not claim that $\iota_u: \hat{N}_u \rightarrow \hat{M}_u$ is injective. In fact, this is generally not true, see the comments made in Remark \ref{RmkIotaUIsNotInjective} and the paragraph before this remark.
\end{rmk}

\section{Weights on homogeneous spaces}\label{SectWeights}

We introduce weights on the von Neumann algebras and C$^\ast$-algebras as were introduced in Section~\ref{SectHomogeneousSpaces}. We study their GNS-representations and prove Proposition \ref{PropWAstLift}, which is essential for implementing \cite[Theorem~3.4.5]{Des}.

Recall that the C$^\ast$-algebraic weights $\hat{\varphi}_u$, $\hat{\varphi_c}$ were def\/ined in  Section \ref{SectPreliminaries}.
The weights on the von Neumann algebras $M$, $\hat{M}$ and the C$^\ast$-algebras $\hat{M}_c$, $\hat{M}_u$ restrict to weights on $N$, $\hat{N}$ and~$\hat{N}_c$,~$\hat{N}_u$ by setting
\[
\varphi^\natural = \varphi \vert _N, \qquad \hat{\varphi}^\natural = \hat{\varphi} \vert_{\hat{N}} \qquad \textrm{and} \qquad \hat{\varphi}_c^\natural = \hat{\varphi}_c\vert_{\hat{N}_c},\qquad \hat{\varphi}_u^\natural = \hat{\varphi}_u \iota_u = \hat{\varphi}_c \hat{\vartheta} \iota_u = \hat{\varphi}_c^\natural \hat{\vartheta}^\natural,
\]
 respectively. We prove that $\varphi^\natural$ and $\hat{\varphi}^\natural$ are normal, semi-f\/inite, faithful weights  and $\hat{\varphi}_c^\natural$ and $\hat{\varphi}_u^\natural$ are lower semi-continuous, densely-def\/ined, non-zero weights. Here the assumption made in Notation~\ref{NtUnimodular} becomes essential.

\begin{prop}\label{PropNormalWeightVNA}
$\varphi^\natural$ is a normal, semi-finite, faithful weight on $N$. Its GNS-representation is given by $(P_\gamma P_\beta \cHil, \Lambda \vert_{N \cap \nphi}, \pi\vert_N)$.
\end{prop}
\begin{proof}
Trivially, $\varphi^\natural$ is normal and faithful. Since $\nphi$ is $\sigma$-weakly dense in $M$, $T_\beta$ and $T_\gamma$ are $\sigma$-weakly continuous and $N = T_\gamma T_\beta(M)$, Propositions \ref{PropTGamma} and \ref{PropTBeta} prove that $\varphi^\natural$ is semi-f\/inite. It is straightforward to check that $(P_\gamma P_\beta \cHil, \Lambda \vert_{N \cap \nphi}, \pi\vert_N)$ satisf\/ies all the properties of a~GNS-representation \cite[Section VII.1]{TakII}.
\end{proof}

\begin{prop}\label{PropDualSqInt}
  For $\omega \in \mathcal{I}$, we find $\tilde{\omega}  \in \mathcal{I}$ and $\xi(\tilde{\omega}  ) = P_\beta P_\gamma \xi(\omega)$. The set $\mathcal{I}_N = \{ \omega \in N_\ast \mid \tilde{\omega} \in \mathcal{I} \}$ is dense in $N_\ast$. The set $\{ \xi(\tilde{\omega}) \mid \omega \in \mathcal{I}_N \}$ is a dense subset of $P_\beta P_\gamma \cHil$.
\end{prop}
\begin{proof}
For $x \in \nphi$, using Propositions \ref{PropTGamma} and \ref{PropTBeta},
\[
 \Lambda(x) \mapsto \omega (T_\beta  T_\gamma(x^\ast)) =   \omega ((T_\beta T_\gamma(x))^\ast) = \langle \xi(\omega), \Lambda(T_\beta T_\gamma(x)) \rangle = \langle P_\beta P_\gamma \xi(\omega), \Lambda(x) \rangle.
\]
The f\/irst claim now follows by the def\/initions of $\mathcal{I}$ and $\xi(\cdot)$. Moreover, we f\/ind $\mathcal{I}_N = \{ \omega\vert_N \mid \omega \in \mathcal{I} \}$, so that the second claim follows by \cite[Lemma 8.5]{KusVae}. The last claim also follows from \cite[Lemma 8.5]{KusVae}.
\end{proof}

\begin{prop}\label{PropNSFWeightVNA}
 $\hat{\varphi}^\natural$  is a normal, semi-finite, faithful weight on $\hat{N}$. Its GNS-representation is given by $(P_\gamma P_\beta \cHil, \hat{\Lambda} \vert_{\hat{N} \cap \mathfrak{n}_{\hat{\varphi}}}, \hat{\pi}\vert_{\hat{N}})$.
\end{prop}
\begin{proof}
By Proposition \ref{PropDualSqInt}, $\{ (\tilde{\omega} \otimes \iota)(W) \mid \omega \in \mathcal{I}_N \} \subseteq \hat{N}$ is a $\sigma$-strong-$\ast$ dense subset of $\hat{N}$ contained in $\mathfrak{n}_{\hat{\varphi}}$. This proves that $\hat{\varphi}^\natural$ is semi-f\/inite. Trivially,  $\hat{\varphi}^\natural$ is normal and faithful.

 To prove the claim about the GNS-representation, we only need to show that the image of~$\hat{\Lambda} \vert_{\hat{N} \cap \mathfrak{n}_{\hat{\varphi}}}$ is contained in   $P_\beta P_\gamma \cHil$ and that the inclusion is dense.
For $\omega \in \mathcal{I}_N$, we have $\lambda(\tilde{\omega}) \in \hat{N} \cap \mathfrak{n}_{\hat{\varphi}}$ and $\hat{\Lambda}(\lambda(\tilde{\omega})) = \xi(\tilde{\omega}) \in P_\beta P_\gamma \cHil$ by Proposition \ref{PropDualSqInt}.  Now, let $x \in \hat{N} \cap \mathfrak{n}_{\hat{\varphi}}$. Since the elements $\lambda(\mathcal{I})$ form a $\sigma$-strong-$\ast$/norm core for $\hat{\Lambda}$, we can f\/ind a net $(\omega_j)_{j \in J}$ in $\mathcal{I}$ such that $\lambda(\omega_j) \rightarrow x$ in the $\sigma$-strong-$\ast$ topology and $\xi(\omega_j) \rightarrow \hat{\Lambda}(x)$ in norm.  Consider the net $(\tilde{\omega}_j)_{j \in J}$. We f\/ind $\lambda(\tilde{\omega}_j)  = P_\gamma \lambda(\omega_j) P_\gamma \rightarrow P_\gamma x P_\gamma = x$. And $\xi(\tilde{\omega}_j) = P_\beta P_\gamma \xi(\omega_j)$ is norm convergent to $ P_\beta P_\gamma \hat{\Lambda}(x)$. Since $\hat{\Lambda}$ is $\sigma$-strong-$\ast$/norm closed, it follows that $\hat{\Lambda}(x) =   P_\beta P_\gamma \hat{\Lambda}(x) \in P_\beta P_\gamma \cHil$. Moreover, it follows from Proposition \ref{PropDualSqInt} that the range of  $\hat{\Lambda} \vert_{\hat{N} \cap \mathfrak{n}_{\hat{\varphi}}}$ is dense in $P_\beta P_\gamma \cHil$.
\end{proof}

We refer to \cite[Section 1.1]{KusVae} for the def\/inition of a GNS-representation for a C$^\ast$-algebraic weight.

\begin{prop}\label{PropNSFWeightReduced}
 $\hat{\varphi}_c^\natural$  is a   proper $($i.e.\ densely defined, lower semi-continuous, non-zero$)$ weight on $\hat{N}_c$. Its GNS-representation is given by $(P_\gamma P_\beta \cHil, \hat{\Lambda} \vert_{\hat{N}_c \cap \mathfrak{n}_{\hat{\varphi}_c}}, \hat{\pi}\vert_{\hat{N}_c})$.
\end{prop}
\begin{proof}
By Proposition \ref{PropDualSqInt}, $\{ (\tilde{\omega} \otimes \iota)(W) \mid \omega \in \mathcal{I}_N \} \subseteq \hat{N}$ is a norm dense subset of $\hat{N}_c$ contained in~$\mathfrak{n}_{\hat{\varphi}}$. The lower semi-continuity and non-triviality  follow since  $\hat{\varphi}_c^\natural$ is a restriction of the faithful weight~$\hat{\varphi}_c$. The claim on the GNS-representation follows exactly as in the proof of Proposi\-tion~\ref{PropNSFWeightVNA}.
\end{proof}

The following lemma can be found as \cite[Lemma 4.2]{Kus}

\begin{lem}\label{LemDensity}
 $\mathcal{I} \cap M_\ast^\sharp$ is dense in $M_\ast^\sharp$ with respect to the norm $\Vert \cdot \Vert_\ast$.
\end{lem}

\begin{prop}\label{PropNSFWeightCast}
$\hat{\varphi}_u^\natural $ is a proper $($i.e.\ densely defined, lower semi-continuous, non-zero$)$ weight on $\hat{N}_u$. Its GNS-representation is given by  $(P_\gamma P_\beta \cHil, \hat{\Lambda}_{u}  \iota_u, \hat{\pi}_{u}  \iota_u )$.
\end{prop}
\begin{proof}
Lemma \ref{LemDensity} shows that  $\mathcal{I} \cap M_\ast^\sharp$ is dense in $M_\ast^\sharp$.
Since~$\mathcal{I}_N$ consists of the restrictions to~$N$ of functionals in~$\mathcal{I}$, see Proposition~\ref{PropDualSqInt}, and~$N_\ast^\sharp$ consists of the restrictions to~$N$ of functionals in~$M_\ast^\sharp$, see Proposition \ref{PropSharpTilde}, it follows that $\mathcal{I}_N \cap N_\ast^\sharp$ is dense in $N_\ast^\sharp$. Hence, $\lambda_u^\natural(\mathcal{I}_N \cap N_\ast^\sharp)$ is dense in $\hat{N}_u$. Moreover, for $\omega \in \mathcal{I}_N \cap N_\ast^\sharp$,
\begin{gather*}
 \hat{\varphi}^\natural_u(\lambda_u^\natural(\omega)^\ast\lambda_u^\natural(\omega) ) = \hat{\varphi}_c \hat{\vartheta} \iota_u (\lambda_u^\natural(\omega)^\ast\lambda_u^\natural(\omega) )
\\
\phantom{\hat{\varphi}^\natural_u(\lambda_u^\natural(\omega)^\ast\lambda_u^\natural(\omega) )}{}
=  \hat{\varphi}_c \hat{\vartheta} (\lambda_u(\tilde{\omega})^\ast\lambda_u(\tilde{\omega}) ) = \hat{\varphi}( (\tilde{\omega} \otimes \iota)(W)^\ast (\tilde{\omega} \otimes \iota)(W) )
< \infty .
\end{gather*}
So   $\lambda_u^\natural(\mathcal{I}_N \cap N_\ast^\sharp)$ is contained in $\mathfrak{n}_{\hat{\varphi}_u^\natural}$. Thus, $\hat{\varphi}^\natural_u$ is densely def\/ined.
 That $\hat{\varphi}_u^\natural $ is lower semi-continuous follows from \cite[Def\/inition~1.5]{KusVae}. Take any $\omega \in N_\ast^\sharp$ such that $\tilde{\omega} \not = 0$, which exists since all functionals in $N_\ast$ are given by restrictions of functionals in $M_\ast$, see Proposition~\ref{PropInvariance}. Then, $\hat{\varphi}_u^\natural ( \lambda_u^\natural( \omega^\ast \ast^\natural \omega) ) = \hat{\varphi}_u( \lambda_u (\tilde{\omega}^\ast \ast \tilde{\omega}) ) = \hat{\varphi}_c(\lambda(\tilde{\omega}^\ast \ast \tilde{\omega})) \not = 0 $, where $\hat{\varphi}_c$ is the faithful left invariant weight on the reduced C$^\ast$-algebraic dual $(\hat{M}_c, \hat{\Delta}_c)$.

Finally, we have to prove that $\hat{\Lambda}_u \iota_u$ maps $\hat{N}_u$ densely into $P_\beta P_\gamma \cHil$. The proof is similar to the one of Proposition \ref{PropNSFWeightVNA}, but since the dif\/ference is subtle we state it here. Observe that $\lambda(\mathcal{I} \cap M^\sharp_\ast)$ is a $\sigma$-strong-$\ast$/norm core for $\hat{\Lambda}$ as follows from \cite[Proposition 2.6]{KusVaeII}. So for every $x \in \mathfrak{n}_{\hat{\varphi}^\natural_u}$ we have $(\hat{\vartheta}\iota_u)(x) \in \nphi$ and hence there exists a net $(\omega_j)_{j \in J}$ in $\mathcal{I} \cap M^\sharp_\ast$ such that $\lambda(\omega_j) \rightarrow (\hat{\vartheta}\iota_u)(x)$ in the $\sigma$-strong-$\ast$ topology and $\xi(\omega_j) \rightarrow \Lambda ((\hat{\vartheta}\iota_u) (x)) = \hat{\Lambda}_u(\iota_u(x))$.  Consider $(\tilde{\omega}_j)_{j \in J}$. Then, $\lambda(\tilde{\omega_j}) \rightarrow (\hat{\vartheta}\iota_u)(x)$ and $\xi(\tilde{\omega}_j) = P_\beta P_\gamma \xi(\omega_j) \rightarrow \hat{\Lambda}_u(\iota_u(x)) $. Hence $\hat{\Lambda}_u(\iota_u(x)) \in P_\beta P_\gamma \cHil$. The range of $\hat{\Lambda}_u \iota_u$ is dense in $P_\beta P_\gamma \cHil$ by Proposition \ref{PropDualSqInt}.
\end{proof}

Next, we like to prove that $\hat{\varphi}^\natural$ is essentially the W$^\ast$-lift of $\hat{\varphi}_u^\natural$, see \cite[Def\/inition~1.31]{KusVae}. A~priori this question is ill-def\/ined, since these weights are def\/ined on dif\/ferent von Neumann algebras. Indeed, $\hat{\varphi}^\natural$ is a weight on $\hat{N}$, which by def\/inition acts on $P_\gamma \cHil$, whereas the~W$^\ast$-lift of $\hat{\varphi}_u^\natural$ is a~weight on $(\hat{\pi}_u\iota_u(\hat{N}_u)) ''$.
Since $P_\beta \in \hat{M}'$, we see that $P_\gamma P_\beta \cHil$ is an invariant subspace of~$\hat{N}$. By Proposition~\ref{PropNSFWeightCast}, the   von Neumann algebra  $(\hat{\pi}_u\iota_u(\hat{N}_u)) ''$ equals the restriction of $\hat{N}$ to $P_\beta P_\gamma \cHil$, i.e.\ $(\hat{\pi}_u\iota_u(\hat{N}_u)) '' = \hat{N} P_\beta$ acting on $P_\beta P_\gamma \cHil$.

The point is that  that $\hat{N}$ and $\hat{N} P_\beta $ are in fact isomorphic. This follows in fact from Proposition~\ref{PropNSFWeightVNA}, but we give a dif\/ferent argument here. We claim, more precisely, that the map $\hat{N} \rightarrow \hat{N} P_\beta: x \mapsto x P_\beta$ is an isomorphism.
 Indeed, for any $x$ in the center of $\hat{M}$, $\hat{J} x \hat{J} = x$. Since $P_\beta = \hat{J} P_\gamma \hat{J}$, every projection in the center of $\hat{M}$  majorizes $P_\beta$ if and only if it majorizes $P_\gamma$. Therefore, the central supports of $P_\beta$ and $P_\gamma$ are equal. It follows from \cite[Theorem~10.3.3]{KadRin} that $\hat{N}$ is isomorphic to $P_\beta \hat{N} P_\beta = \hat{N} P_\beta$, where the isomorphism is given by the map $\hat{N} \rightarrow \hat{N} P_\beta: x \mapsto x P_\beta$.

We emphasize that $\hat{N}$ will always be considered as a von Neumann algebra acting on $P_\gamma \cHil$, whereas if we encounter $(\hat{\pi}_u\iota_u(\hat{N}_u)) ''$ we assume that it acts on $P_\beta P_\gamma \cHil$. The described isomorphism $\hat{N} \simeq (\hat{\pi}_u\iota_u(\hat{N}_u)) ''$, makes the following proposition well-def\/ined. A similar argument holds on the reduced level.

\begin{prop}\label{PropWAstLift}
The W$^\ast$-lifts of $\hat{\varphi}_{c}^\natural$ and $\hat{\varphi}_{u}^\natural$   to $\hat{N}$ both equal  $\hat{\varphi}^\natural$.
\end{prop}

\begin{proof}
We prove the proposition for $\hat{\varphi}_{u}^\natural$, the proof for $\hat{\varphi}_{c}^\natural$ is similar.
 Recall from Proposition~\ref{PropNSFWeightVNA} that $(P_\gamma P_\beta \cHil,   \hat{\Lambda}\vert_{\hat{N} \cap \mathfrak{n}_{\hat{\varphi}}} , \hat{\pi}\vert_{\hat{N}})$ gives the GNS-representation of $\hat{\varphi}^\natural$.
Denote the W$^\ast$-lift of~$\hat{\varphi}_{u}^\natural$ by~$\phi$. We denote its GNS-representation by $(\cHil_\phi, \pi_\phi, \Lambda_\phi)$. Recall from Proposition~\ref{PropNSFWeightCast} that the GNS-representation of $\hat{\varphi}_u^\natural$ was given by  $(P_\gamma P_\beta \cHil, \hat{\Lambda}_{u}  \iota_u, \hat{\pi}_{u}  \iota_u )$.

 From \cite[Proposition 1.32]{KusVae}, the elements $\hat{\pi}_{u}  \iota_u(x), x \in \mathfrak{n}_{\hat{\varphi}_{u}^\natural}$, form a $\sigma$-strong-$\ast$/norm core for~$\Lambda_\phi$ and we may take $\cHil_\phi = P_\gamma P_\beta  \cHil$. If we can prove that    the elements $\hat{\pi}_{u}  \iota_u(x), x \in \mathfrak{n}_{\hat{\varphi}_{u}^\natural}$, also form a $\sigma$-strong-$\ast$/norm core for $\hat{\Lambda}\vert_{\hat{N} \cap \mathfrak{n}_{\hat{\varphi}}}$, then $\phi$ and $\hat{\varphi}^\natural$ have identical GNS-representations and hence they are equal. So let $x \in \mathfrak{n}_{\hat{\varphi}} \cap \hat{N}$. By \cite[Proposition~2.6]{KusVaeII}, $\lambda(\mathcal{I} \cap M_\ast^\sharp)$  forms a~$\sigma$-strongly-$\ast$/norm core for $\hat{\Lambda}$. So let $(\omega_j)_{j \in J}$ be a net in $\mathcal{I} \cap M_\ast^\sharp$  such that $\lambda(\omega_j) \rightarrow x$ in the $\sigma$-strong-$\ast$ topology and $\xi(\omega_j) \rightarrow \hat{\Lambda}(x)$ in norm. Then, $\xi(\tilde{\omega}_j) = P_\gamma P_\beta \xi(\omega_j)$ converges in norm. Furthermore,
\[
\lambda(\tilde{\omega}_j) = (\tilde{\omega}_j \otimes \iota)(W) = P_\gamma (\omega_j \otimes \iota)(W) P_\gamma = P_\gamma \lambda(\omega_j) P_\gamma \rightarrow P_\gamma x P_\gamma = x,
\]
where the convergence is in the $\sigma$-strong-$\ast$ topology. Note that $\lambda(\tilde{\omega}_j) = (\hat{\pi}_u \iota_u \lambda_u^\natural)(\omega_j\vert_N)$. All in all, we conclude that $(y_j)_{j \in J} := ((\lambda_u^\natural)(\omega_j\vert_N))_{j \in J}$ is a net in $ \mathfrak{n}_{\hat{\varphi}_{u}^\natural}$ such that $\hat{\pi}_{u}  \iota_u(y_j)$ is $\sigma$-strong-$\ast$ convergent to $x$ and $\hat{\Lambda}\vert_{\hat{N} \cap \mathfrak{n}_{\hat{\varphi}}}(\hat{\pi}_{u}  \iota_u(y_j))$ is convergent.
\end{proof}

\begin{rmk}
In particular, Proposition \ref{PropNSFWeightVNA} implies that $\hat{\varphi}_{c}^\natural$ and $\hat{\varphi}_{u}^\natural$ are approximately KMS-weights, see \cite[Proposition 1.35]{KusVae}.
\end{rmk}

\begin{rmk}
Note that on one hand, we have a map $N_\ast \rightarrow \hat{N}: \omega \mapsto \lambda(\tilde{\omega}) = (\tilde{\omega} \otimes \iota)(W)$. On the other hand, we can def\/ine a map $\hat{N}_\ast \rightarrow N: \omega \mapsto \hat{\lambda}(P_\gamma \omega P_\gamma) = (\iota \otimes \omega)(1\otimes P_\gamma)(W^\ast)(1\otimes P_\gamma)$. These can be considered as the spherical $L^1$-Fourier transforms. Since both $\varphi^\natural$ and $\hat{\varphi}^\natural$ are normal, semi-f\/inite, faithful weights, one can proceed as in \cite{Cas} to obtain a spherical $L^p$-Fourier transform which then is a restriction of the $L^p$-Fourier transform as def\/ined in \cite[Theorem 5.6]{Cas}.
\end{rmk}

\section{Spherical corepresentations}\label{SectSphericalCoreps}
 We introduce the necessary terminology for corepresentations that admit vectors that are invariant under the action of a quantum subgroup. These corepresentations can be considered as spherical corepresentations.

\begin{dfn}\label{DfnMOneInv}
 Let $U \in M \otimes B(\cHil_U)$ be a corepresentation. Then  $v \in \cHil_U$   is called a $M_1$-invariant vector if
\[
 (\omega \otimes \iota)((T_\beta \otimes \iota)(U)) v = (\omega \otimes \iota)(U)v, \qquad \forall\, \omega\in M_\ast.
\]
 We denote
\[
\cHil_U^{M_1} = \left\{v \in \cHil_U \mid v \textrm{ is } M_1 \textrm{-invariant} \right\}.
\]
 Note that $\cHil_U^{M_1}$ is a closed subspace of $\cHil_U$. We denote
$
\ICMM
$
 for the equivalence classes of irreducible corepresentations of $M$ that admit non-trivial $M_1$-invariant vectors. We refer to such corepresentations as spherical corepresentations. If $\{ (\omega \otimes \iota)(U) \cHil_U^{M_1} \mid \omega \in M_\ast \}$ is dense in $\cHil_U$, then we call $U$ homogeneously cyclic. It should be clear that every irreducible corepresentation of $M$ that admits a non-trivial $M_1$-invariant vector is homogeneously cyclic.
 \end{dfn}

\begin{prop}\label{PropMOneInv}
$ v \in P_\gamma \cHil$ if and only if $v$ is $M_1$-invariant for $W$.
\end{prop}
\begin{proof}
Observe that $ (T_\beta \otimes \iota)(W) = W \left(\iota \otimes \hat{\pi}((\varphi_1 \otimes \iota)(W_1))\right) = W (\iota \otimes P_\gamma)$. This yields the only if part. The other implication follows by taking a  net $(\omega_j)_{j \in J}$ in $M_\ast$ such that $\lambda(\omega_j) \rightarrow 1$ in the $\sigma$-strong-$\ast$ topology. Then,
\[
v \leftarrow (\omega_j \otimes \iota)(W) v =  (\omega_j \otimes \iota)(T_\beta \otimes \iota)(W) v = (\omega_j \otimes \iota)( W)  P_\gamma v \rightarrow P_\gamma v,
\]
where the convergence is in norm.
\end{proof}

\begin{lem}
 Let $U \in M \otimes B(\cHil_U)$ be a corepresentation. Let $v \in \cHil_U$ and $\omega \in N_\ast$. The vector $(\tilde{\omega} \otimes \iota)(U) v$ is $M_1$-invariant.
\end{lem}

\begin{proof}
 This follows from the following series of equalities for which we use Lemma \ref{LemTProperties} (\ref{LemTIII}) and $T_\gamma^2 = T_\gamma$. For $\theta \in M_\ast$,
 \begin{gather*}
  (\theta \otimes \iota)((T_\beta \otimes \iota)(U))(\tilde{\omega} \otimes \iota)(U)v = (\theta  T_\beta \otimes \omega  T_\beta T_\gamma \otimes \iota) U_{13} U_{23} v  \\
    \qquad{}
=   (\theta  T_\beta \otimes \omega T_\beta T_\gamma \otimes \iota) (\Delta \otimes \iota)(U) v
=
(\theta  \otimes \omega T_\beta T_\gamma \otimes \iota) (\Delta \otimes \iota)(U) v\\
\qquad{}
= (\theta \otimes \iota)(U) (\tilde{\omega} \otimes \iota)(U)v.\tag*{\qed}
 \end{gather*}
 \renewcommand{\qed}{}
\end{proof}

We mention the following two propositions in order to compare our framework with the setting of classical Gelfand pairs of groups. The proof of the f\/irst one is completely analogous to the proof of \cite[Proposition II.6]{Far} or \cite[Lemma 6.2.3]{Dij}. For Proposition \ref{PropAbDim}, one proves that the  representation $N_\ast^\sharp \rightarrow B(\cHil_U^{M_1}): \omega \rightarrow (\tilde{\omega} \otimes \iota)(U)\vert_{\cHil_U^{M_1}}$ is irreducible. This can be done along the lines of \cite[Th\'eor\`eme III.1]{Far} or \cite[Proposition 6.31]{Dij}.

\begin{prop}
 Let $U \in M \otimes B(\cHil_U)$ be a corepresentation. Suppose that $U$ admits a $M_1$-invariant, cyclic vector $v$. If $\dim(\cHil_U^{M_1}) = 1$, then $U$ is irreducible.
\end{prop}

\begin{prop}\label{PropAbDim}
Assume that $\hat{N}$ is Abelian. Let $U \in M \otimes B(\cHil_U)$ be an irreducible   corepresentation. Then, $\dim(\cHil_U^{M_1}) \leq 1$.
\end{prop}

In particular, suppose that  $\hat{N}$ is Abelian, and that $U \in M \otimes B(\cHil_U)$ is an irreducible corepresentation with $M_1$-invariant unit vector $v$. Then $x = (\iota \otimes \omega_{v,v})(U)$ satisf\/ies
\begin{gather}\label{EqnProdForm}
\Delta^\natural(x) = x \otimes x.
\end{gather}
Hence, $x$ is a character of the convolution algebra $N_\ast^\sharp$. It can be considered as a quantum spherical function or quantum spherical element. The equality (\ref{EqnProdForm}) allows one to derive product formulae as is done in for example \cite{VaiE2, VaiGel,VaiHyp}. Here we keep to a more general setting and do not assume that $\hat{N}$ is Abelian in order to include the example of $SU_q(1,1)_{{\rm ext}}$ in Section~\ref{SectExample}.

\begin{rmk}\label{RmkInvSubs}
Let $U \in M \otimes B(\cHil_U)$ be a corepresentation. It follows in particular that $\cHil_U^{M_1}$ is a closed invariant subspace for the  representation $N_\ast^\sharp \rightarrow B(\cHil_U): \omega \mapsto (\tilde{\omega} \otimes \iota)(U)$. By the universal property of $\hat{N}_u$, we see that this gives rise to a representation of $\hat{N}_u$ on $\cHil_U^{M_1}$. Of course, this representation can be trivial. If $U$ is homogeneously cyclic, then the corresponding representation of $\hat{N}_u$ is non-degenerate.  Indeed, suppose that there exist $w \in \cHil_U^{M_1}$, such that for all $v \in \cHil_U^{M_1}$ and $\omega \in M_\ast^\sharp$, $  \langle (\tilde{\omega} \otimes \iota )(U) v, w \rangle =0 $. Then, using (\ref{LemTIV}) of Lemma \ref{LemTProperties},
\begin{gather*}
0 = \langle (\omega  T_\gamma T_\beta \otimes \iota)(U) v, w \rangle =
\langle (\omega T_\gamma  \otimes \iota)(U) v, w \rangle =
\langle  v, (\omega  T_\gamma  \otimes \iota)(U)^\ast w \rangle  \\
\phantom{0}{} =  \langle  v, (\overline{\omega  T_\gamma}  \otimes \iota)(U^\ast) w \rangle =
\langle  v, (\overline{\omega}  T_\gamma  \otimes \iota)(S \otimes \iota)(U) w \rangle  \\
\phantom{0}{} =  \langle  v, (\overline{\omega} \otimes \iota) (S \otimes \iota) (T_\beta \otimes \iota)(U) w \rangle =
 \langle  v, (\omega^\ast  \otimes \iota)(T_\beta \otimes \iota)(U) w\rangle  \\
\phantom{0}{} =  \langle  v, (\omega^\ast  \otimes \iota)(U) w \rangle =
\langle (\omega  \otimes \iota)(U)  v,  w \rangle.
\end{gather*}
We see that for all $v \in \cHil_U^{M_1}$ and $\omega \in M_\ast^\sharp$,
$  \langle (\omega \otimes \iota )(U) v, w \rangle =0,$
which proves that $w = 0$, since~$U$ is homogeneously cyclic. Hence, every non-zero homogeneously cyclic corepresentation~$U$ of~$M$ gives rise to a  non-degenerate representation of $\hat{N}_u$.
\end{rmk}

\begin{dfn}\label{DfnRep}
Let $U \in M \otimes B(\cHil_U)$ be a homogeneously cyclic corepresentation of $M$ on a~Hilbert space $\cHil_U$. Then, we get a representation $\pi_U^\natural$ of $\hat{N}_u$ determined by
\[
\pi_U^\natural: \ \lambda_u^\natural(\omega) \mapsto (\tilde{\omega} \otimes \iota)(U) \vert_{\cHil_{U}^{M_1}}, \qquad \omega \in N_\ast^\sharp.
\]
 We emphasize that the representation Hilbert space of $\pi_U^\natural$ is $\cHil_U^{M_1}$.
\end{dfn}

\begin{rmk}\label{RmkIr}
Let $U \in M \otimes B(\cHil_U)$ be a homogeneously cyclic corepresentation of $M$. $U$ is irreducible if and only if $\pi_U^\natural$ is irreducible. Indeed, if $U$ is reducible, then clearly $\pi_U^\natural$ is reducible. The only if part follows from a computation similar to the one in Remark \ref{RmkInvSubs}.
\end{rmk}

 Recall that we denote $\pi_U$ for the representation of $\hat{M}_u$ given by $\lambda_u(\omega) \mapsto (\omega \otimes \iota)(U)$. Then, $\pi_U^\natural$ equals the restriction of $\pi_U \iota_u$ to $\cHil_U^{M_1}$.

We will need the following result for Theorem \ref{ThmPG}. We refer to \cite{DixVNA} for the theory of direct integration and the def\/inition of a fundamental sequence.

\begin{prop}\label{PropMeas}
Let $X$ be a measure space, with standard measure $\mu$. Suppose that for every $x \in X$, we have a homogeneously cyclic corepresentation $U_x$ of $M$ on a Hilbert space $\cHil_x$ such that $(\pi_{U_x}^\natural)_{x \in X}$ is a $\mu$-measurable field of representations of $\hat{N}_u$. Suppose that $\hat{M}_u$ is separable. Then, $(\cHil_x)_{x \in X}$ is a $\mu$-measurable field of Hilbert spaces such that $(\cHil_x^{M_1})_{x \in X}$ is a $\mu$-measurable field of subspaces and $(U_x)_{x \in X}$ is a $\mu$-measurable field of corepresentations.
\end{prop}
\begin{proof}
Let $\omega_i \in N_\ast^\sharp$, $i \in \mathbb{N}$ be a such that $\lambda_u(\omega_i)$, $i \in \mathbb{N}$ is dense in $\hat{M}_u$. Since $(\pi_{U_x}^\natural)_{x \in X}$ is $\mu$-measurable, we have a fundamental sequence $(e_x^j)_{x\in X}$, $j \in \mathbb{N}$ for the $\mu$-measurable f\/ield of Hilbert spaces $(\cHil_x^{M_1})_{x \in X}$. For $i, j \in \mathbb{N}$, $x \in X$, def\/ine $f^{i,j}_x \in \cHil_x$ by
\[
f^{i,j}_x = (\omega_i \otimes \iota)(U_x) e^j_x.
\]
We claim that $(f^{i,j}_x)_{x\in X}$ is a fundamental sequence for $(\cHil_x)_{x \in X}$. Indeed, since $U_x$ is homogeneously cyclic, the span of $(\omega_i \otimes \iota)(U_x) e^j_x$, $i, j \in \mathbb{N}$ is dense in $\cHil_x$. Moreover,
\begin{gather}
\langle f^{i,j}_x , f^{i',j'}_x \rangle =   \langle (\omega_i \otimes \iota)(U_x) e^j_x, (\omega_{i'} \otimes \iota)(U_x) e^{j'}_x \rangle =
\langle ((\omega_{i'}^\ast \ast \omega_i) \otimes \iota)(U_x) e^j_x, e^{j'}_x \rangle
\nonumber \\
\phantom{\langle f^{i,j}_x , f^{i',j'}_x \rangle}{}
=
\langle ((\omega_{i'}^\ast \ast \omega_i)^\sim \otimes \iota)(U_x) e^j_x, e^{j'}_x \rangle = \langle \pi_{U_x}^\natural(\lambda_u^\natural((\omega_{i'}^\ast \ast \omega_i)\vert_N))  e^j_x, e^{j'}_x \rangle,\label{EqnMeasField}
\end{gather}
where the third equality follows by a computation similar to the one in Remark~\ref{RmkInvSubs}.
Since $(\pi_{U_x}^\natural)_{x\in X}$ is a $\mu$-measurable f\/ield of representations, we see that (\ref{EqnMeasField}) is a $\mu$-measurable function of~$x$. Hence $(f^{i,j}_x)_{x\in X}$ is a fundamental sequence. Moreover, by a similar computation as (\ref{EqnMeasField}), for any $\omega \in M_\ast^\sharp$,
\begin{gather*}
\langle (\omega \otimes \iota)(U_x) f^{i,j}_x, f^{i',j'}_x \rangle =
\langle ((\omega_{i'}^\ast \ast \omega \ast \omega_i) \otimes \iota)(U_x) e^j_x, e^{j'}_x \rangle\\
\phantom{\langle (\omega \otimes \iota)(U_x) f^{i,j}_x, f^{i',j'}_x \rangle}{} =
 \langle \pi_{U_x}^\natural(\lambda_u^\natural((\omega_{i'}^\ast \ast \omega \ast \omega_i)\vert_N))  e^j_x, e^{j'}_x \rangle,
\end{gather*}
is a $\mu$-measurable function of $x$. For $\omega = \omega_{v, w}$ with $v,w \in \Dom(\hat{\nabla}^{\frac{1}{2}}) \cap \Dom(\hat{\nabla}^{- \frac{1}{2}} )$, we have $\omega \in M_\ast^\sharp$ by \cite[Proposition~1.10]{Des}. Using \cite[Proposition~II.1.10 and~II.2.1]{DixVNA} it is straightforward to prove that $(U_x)_{x \in X}$ is a $\mu$-measurable f\/ield of operators.
\end{proof}

\begin{prop}\label{PropEqRepCorep}
 Let $U_1$ and $U_2$ be homogeneously cyclic corepresentations of $M$ on Hilbert spaces $\cHil_1$ and $\cHil_2$ respectively. Suppose that the  representations of $N_\ast^\sharp$ given by $\pi_i: \omega \mapsto (\tilde{\omega} \otimes \iota)(U_i)\vert_{\cHil_{U_i}^{M_1}}$, $\omega \in N_\ast^\sharp$, $i \in \{ 1,2\}$ are equivalent. Then $U_1$ and $U_2$ are equivalent.
\end{prop}

\begin{proof}
Let $T: \cHil_{U_1}^{M_1} \rightarrow \cHil_{U_2}^{M_1} $ be the unitary intertwiner between $\pi_1$ and $\pi_2$. Let $Q_0$ be the mapping
\begin{gather*}
 \left\{ (\omega \otimes \iota)(U_1) v \mid \omega \in M_\ast^\sharp, v \in \cHil_{U_1}^{M_1} \right\}  \rightarrow
 \left\{ (\omega \otimes \iota)(U_2) w \mid \omega \in M_\ast^\sharp, w \in \cHil_{U_2}^{M_1} \right\}, \\
(\omega \otimes \iota)(U_1) v   \mapsto (\omega \otimes \iota)(U_2) T v.
\end{gather*}
This map is well-def\/ined and isometric. Indeed, for $\omega \in M_\ast^\sharp$ and $v \in \cHil_{U_1}^{M_1}$,
\begin{gather*}
  \Vert (\omega \otimes \iota)(U_2) T v \Vert^2 = \langle (\omega^\ast \ast \omega \otimes \iota)(U_2) T v, T v \rangle =
 \langle( (\omega^\ast \ast \omega )^\sim \otimes \iota)(U_2) T v, T v \rangle \\
 \phantom{\Vert (\omega \otimes \iota)(U_2) T v \Vert^2}{}  =
\langle( (\omega^\ast \ast \omega )^\sim\otimes \iota)(U_1)  v,  v \rangle =
 \langle (\omega^\ast \ast \omega \otimes \iota)(U_1)  v,  v \rangle =
 \Vert (\omega \otimes \iota)(U_1) v \Vert^2,
\end{gather*}
where the second equality follows from a similar calculation as in Remark~\ref{RmkInvSubs}.
Since $U_1$ and $U_2$ are homogeneously cyclic, $Q_0$ is densely def\/ined and has dense range. Let $Q: \cHil_{U_1} \rightarrow \cHil_{U_2}$ be the unitary extension of $Q_0$.  Let $\omega, \omega_1, \omega_2 \in M_\ast^\sharp$ and $v, w \in \cHil_{U_1}^{M_1}$,
\begin{gather*}
  \langle (\omega \otimes \iota)(U_1) (\omega_1 \otimes \iota)(U_1) v, (\omega_2 \otimes \iota)(U_1) w \rangle =
 \langle (\omega_2^\ast \ast \omega \ast \omega_1 \otimes \iota)(U_1)  v, w \rangle  \\
\qquad{}=  \langle ( (\omega_2^\ast \ast \omega \ast \omega_1)^\sim \otimes \iota)(U_1)  v, w \rangle =
 \langle ( (\omega_2^\ast \ast \omega \ast \omega_1)^\sim \otimes \iota)(U_2) T  v,  T w \rangle  \\
\qquad {} = \langle (\omega_2^\ast \ast \omega \ast \omega_1 \otimes \iota)(U_2) T  v,  T w \rangle =
\langle (\omega \otimes \iota)(U_2) (\omega_1 \otimes \iota)(U_2) T v, (\omega_2 \otimes \iota)(U_2) T w \rangle  \\
\qquad {}= \langle (\omega \otimes \iota)(U_2) Q (\omega_1 \otimes \iota)(U_1)  v, Q (\omega_2 \otimes \iota)(U_1)  w \rangle,
\end{gather*}
where the second equality follows again from a similar calculation as in Remark~\ref{RmkInvSubs}.
Since $U_1$ is homogeneously cyclic, this proves that $Q$ intertwines~$U_1$ with~$U_2$.
\end{proof}

Note that the converse of the previous proposition is clear: if $U_1$ and $U_2$ are equivalent corepresentations, then the corresponding representations as considered in Remark \ref{RmkInvSubs} are equivalent.

\begin{rmk}\label{RmkAlsoUniv}
Proposition \ref{PropEqRepCorep} and its converse also hold on the universal level. So let~$U_1$ and~$U_2$ be homogeneously cyclic corepresentations of $M$. $\pi_{U_1}^\natural$ and~$\pi_{U_2}^\natural$
 are equivalent if and only if~$U_1$ and~$U_2$ are equivalent.
\end{rmk}

\section[Representations of $\hat{M}_u$, $\hat{M}_c$, $\hat{N}_u$ and $\hat{N}_c$]{Representations of $\boldsymbol{\hat{M}_u}$, $\boldsymbol{\hat{M}_c}$, $\boldsymbol{\hat{N}_u}$ and $\boldsymbol{\hat{N}_c}$}\label{SectCorrespondence}

In this section, we compare the representations of the C$^\ast$-algebras def\/ined in Sections~\ref{SectPreliminaries} and~\ref{SectHomogeneousSpaces}. Main objective is to prove that the representations of $\hat{N}_c$ `lift' to representations of~$\hat{M}_c$.

Let us give a more elaborate discussion. There are three special types of representations within $\Rep(\hat{N}_u)$.
\begin{enumerate}\itemsep=0pt
\item\label{ItemRepI} As explained in Remark \ref{RmkInvSubs}, the corepresentations of $M$ give rise to representations of~$\hat{N}_u$. Recall \cite{Kus} that the corepresentations of $M$ are in 1-to-1 correspondence with non-degenerate representations of~$\hat{M}_u$. Hence, the representations of $\hat{M}_u$ give rise to representations of $\hat{N}_u$. This correspondence can be described more directly: if~$\pi$ is a representation of $\hat{M}_u$ on a Hilbert space $\cHil_\pi$, then $\pi\iota_u$ is the corresponding representation of~$\hat{N}_u$ on the closure of $((\pi\iota_u)(\hat{N}_u)) \cHil_\pi$.  Note that by Remark~\ref{RmkAlsoUniv}, this assignment  descends to a well-def\/ined, injective map on the equivalence classes of representations,
\[
\Rep(\hat{M}_u) \hookrightarrow \Rep(\hat{N}_u): \ \pi \mapsto \pi \iota_u.
\]
\item\label{ItemRepII} If $\pi$ is a representation of $\hat{N}_c$, then $\pi \hat{\vartheta}^\natural$ is a representation of~$\hat{N}_u$. These representations correspond to the representations that are weakly contained in the GNS-representation of~$\hat{\varphi}^\natural_u$. Indeed this follows, since by Proposition~\ref{PropNSFWeightCast}, this GNS-representation is given by
\[
\hat{\pi}_u \iota_u = \hat{\pi}_c \hat{\vartheta} \iota_u = \hat{\pi}_c \hat{\vartheta}^\natural = \hat{\pi}\vert_{\hat{N}_c} \hat{\vartheta}^\natural.
\]
Hence, every representation   $\pi \hat{\vartheta}^\natural$, with $\pi \in \Rep(\hat{N}_c)$ is weakly contained in the GNS-representation of~$\hat{N}_u$. The other way around, any representation of $\hat{N}_u$ that is weakly contained in the GNS-representation of $\hat{\varphi}^\natural_u$ factors through the canonical projection $\hat{\vartheta}^\natural$.
\item\label{ItemRepIII} If $\pi$ is a representation of $\hat{M}_c$, then $\pi \hat{\vartheta}^\natural = \pi \hat{\vartheta} \iota_u$ is a representation of $\hat{N}_u$.  Here, we  used that $\hat{N}_c \subseteq \hat{M}_c$.
\end{enumerate}

The main results of this section will be the following. We prove that every representation of~$\hat{N}_c$ comes from a representation of~$\hat{M}_c$, i.e.\ the representations of $\hat{N}_u$ obtained in~(\ref{ItemRepII}) and~(\ref{ItemRepIII}) are the same ones.

\begin{thm}\label{ThmCorrespondenceC}
For every non-degenerate representation $\rho \in \Rep(\hat{N}_c)$, there exists a non-degenerate representation $\pi \in \Rep(\hat{M}_c)$ on a Hilbert space $\cHil_\pi$ such that $\rho$ is equivalent to the restriction of $\pi\vert_{\hat{N}_c}$ to the closure of $\pi(\hat{N}_c) \cHil_{\pi}$.
\end{thm}
  \begin{proof}
Let $M_\ast^{\sharp,\beta}$ denote space of functionals $\theta \in M_\ast^\sharp$ such that $\theta T_\beta = \theta$. By a similar argument as in the proof of Proposition \ref{PropSharpTilde}, we see that if $\theta \in M_\ast^\sharp$, then $\theta T_\beta \in M_\ast^{\sharp, \beta}$ and $(\theta T_\beta)^\ast = \theta^\ast T_\gamma$.  In particular, $(N_\ast^\sharp)^\sim \subseteq M_\ast^{\sharp,\beta}$. Note that for $\theta_1, \theta_2 \in  M_\ast^{\sharp, \beta}$, we have that $(\theta_1^\ast \ast \theta_2)^\sim = \theta_1^\ast \ast \theta_2$.

We complete $M_\ast^{\sharp,\beta}$ into a (right) Hilbert $\hat{N}_c$-module. For $\theta, \theta_1, \theta_2 \in M_\ast^{\sharp,\beta}$ and $\omega \in N_\ast^\sharp$, we put
\begin{gather}
\theta \cdot \omega    =  \theta \ast \tilde{\omega}  \in M_\ast^{\sharp,\beta}, \label{EqnRightAction}\\
\langle \theta_1, \theta_2 \rangle_{N_\ast^\sharp}   =  (\theta_1^\ast \ast \theta_2) \vert_N  \in N_\ast^\sharp. \label{EqnInnerProduct}
\end{gather}
The fact that (\ref{EqnRightAction}) is in $M_\ast^{\sharp,\beta}$ follows from Lemma~\ref{LemTProperties}. That~(\ref{EqnInnerProduct}) is in $N_\ast^\sharp$ follows from Proposition~\ref{PropSharpTilde}.
This gives a right $N_\ast^\sharp$-module structure on $M_\ast^{\sharp,\beta}$. We will apply \cite[Lemma~2.16]{RaeWil} to get a Hilbert $\hat{N}_c$-module $X$. Here we consider $N_\ast^\sharp$ is a subalgebra of $\hat{N}_c$  by means of the map~$\lambda^\natural$. To continue, note that conditions~(a) and~(c) of \cite[Def\/inition~2.1]{RaeWil} are indeed satisf\/ied. (b)~follows, since for $\theta_1, \theta_2 \in M_\ast^{\sharp,\beta}$ and $\omega \in N_\ast^\sharp$,
\begin{gather*}
\langle \theta_1, \theta_2 \cdot \omega \rangle_{N_\ast^\sharp}  = (\theta_1^\ast \otimes \theta_2 \otimes \omega) (\iota \otimes \iota \otimes T_\beta T_\gamma) (\Delta \otimes \iota)  \Delta\vert_N, \\
\langle \theta_1, \theta_2 \rangle_{N_\ast^\sharp} \ast^\natural \omega  = (\theta_1^\ast \otimes \theta_2 )  \Delta \vert_N \ast^\natural \omega = (\theta_1^\ast \otimes \theta_2 \otimes \omega) (\iota \otimes \iota \otimes T_\gamma) (\Delta \otimes \iota) \Delta \vert_N.
\end{gather*}
Since $(\iota \otimes T_\beta) \Delta\vert_N = \Delta T_\beta\vert_N = \Delta\vert_N$, these expressions are equal. (d) follows, since for $\theta \in M_\ast^{\sharp,\beta}$,
\[
\lambda^\natural((\theta^\ast \ast \theta )\vert_N ) = \lambda(\theta^\ast \ast \theta) = \lambda(\theta)^\ast \lambda(\theta) \geq 0,
\]
in $\hat{N}_c$.  This def\/ines the right  Hilbert $\hat{N}_c$-module $X$ and we denote its norm by $\Vert\cdot \Vert_X$.

We are able to def\/ine an action of $M_\ast^\sharp$ on $X$ by  means of adjointable operators which extends the convolution product on $M_\ast^\sharp$. Indeed, for $\omega\in M_\ast^\sharp$, and $\theta_1, \theta_2 \in M_\ast^{\sharp,\beta} $,
\begin{gather*}
\langle \omega \ast \theta_1, \theta_2 \rangle_{N_\ast^\sharp} = \theta_1^\ast \ast \omega^\ast \ast \theta_2\vert_{N} =  \langle  \theta_1, \omega^\ast \ast \theta_2 \rangle_{N_\ast^\sharp}.
\end{gather*}
Furthermore, for $\omega \in M_\ast^\sharp$ and $\theta \in M_\ast^{\sharp, \beta}$,
\begin{gather}\label{EqnAdjAction}
\Vert \omega \ast \theta \Vert_X^2 = \Vert \lambda (\theta^\ast \ast \omega^\ast \ast \omega \ast \theta) \Vert \leq \Vert \lambda(\omega) \Vert^2 \Vert \lambda(\theta^\ast \ast \theta) \Vert = \Vert \lambda(\omega) \Vert^2 \Vert \theta\Vert_X^2.
\end{gather}
So the action of $\omega \in M_\ast^\sharp$ is bounded. This allows us to extend the action of $\omega$ to an adjointable operator on $X$. We denote this operator by $L_\omega$.  Moreover, from~(\ref{EqnAdjAction}) we get a representation of~$\hat{M}_c$ on the Hilbert module $X$ by means of adjointable operators. This representation is uniquely determined by $\lambda^\natural(\omega)   \mapsto L_\omega$, $\omega \in N_\ast^\sharp$.

Now, let $\rho$ be a representation of $\hat{N}_c$ on a Hilbert space $\cHil_\rho$. By  \cite[Proposition 2.66]{RaeWil} we get an induced representation of $\hat{M}_c$ on $X  \otimes_{\hat{N}_c} \cHil_\rho$. Let us denote the latter by~$\Ind \rho$. Let $(a_j)_{j \in J}$ be an approximate unit for $\hat{M}_c$. Then, for $\theta \in M_\ast^{\sharp, \beta}, v \in \cHil_\rho$, we see that $(\Ind \rho)(a_j) (\theta \otimes v) \rightarrow (\theta \otimes v)$ in $X$. So $\Ind \rho$ is non-degenerate.

 Note that the completion of $  (N_\ast^{\sharp})^{\sim}   \otimes_{\hat{N}_c}  \cHil_\rho$ is a closed subspace of the Hilbert space $X  \otimes_{\hat{N}_c} \cHil_\rho$ that is isomorphic to $\cHil_\rho$ via the unitary  extension $T$ of
\begin{gather}\label{EqnUnitaryIntertwiner}
( N_\ast^\sharp) ^\sim  \otimes_{\hat{N}_c} \cHil_\rho \rightarrow \cHil_\rho: \ \tilde{\omega} \otimes v \mapsto \rho(\lambda ^\natural(\omega)) v.
\end{gather}
The map (\ref{EqnUnitaryIntertwiner}) extends unitarily since $\rho$ is non-degenerate. For $\omega, \theta \in N_\ast^\sharp$, $v \in \cHil_\rho$,
\[
(\Ind \rho)(\lambda^\natural(\omega))  (\tilde{\theta} \otimes v) = \tilde{\omega} \ast \tilde{\theta} \otimes v = (\omega \ast^\natural \theta)^\sim \otimes v,
\]
 so that $(N_\ast^\sharp)^\sim  \otimes \cHil_\rho$ is an invariant subspace for $ \Ind \rho $. We denote its closure by $Y$. Moreover, for $\omega, \theta \in N_\ast^\sharp$ and $v \in \cHil_\rho$,
\begin{gather*}
  T    (\Ind \rho)  (\lambda^\natural(\omega))   ( \tilde{\theta} \otimes v )  =
T ((\tilde{\omega} \ast \tilde{\theta})  \otimes v)  = T ((\omega \ast^\natural \theta)^\sim \otimes v)  \\
\phantom{T    (\Ind \rho)  (\lambda^\natural(\omega))   ( \tilde{\theta} \otimes v ) }{}
=   \rho(\lambda^\natural(\omega \ast^\natural \theta)) v  =   \rho(\lambda^\natural(\omega)) \rho(\lambda^\natural( \theta)) v =
  \rho(\lambda^\natural(\omega)) T  ( \tilde{\theta} \otimes v ),
\end{gather*}
so that $T$ intertwines  $ (\Ind \rho)\vert_{\hat{N}_c}  $ restricted to the Hilbert space   $Y$ with $\rho$. Finally, we claim that~$Y$ equals the closure of $(\Ind\rho)(\hat{N}_c) X$. Indeed, for any $\omega \in N_\ast^\sharp$, $\theta \in M_\ast^{\sharp, \beta}$ and $v \in \cHil_\rho$, we see that $(\Ind\rho)(\lambda^\natural(\omega)) ( \theta \otimes v) = (\tilde{\omega} \ast \theta) \otimes v = (\omega \ast^\natural \theta\vert_N)^\sim \otimes v \in Y$. Since $\hat{N}_c \hat{N}_c \subseteq \hat{N}_c$ is dense, it is straightforward to prove that  $Y$ equals the closure of $(\Ind\rho)(\hat{N}_c) X$ in $X$. This concludes the proof by choosing $\pi = \Ind \rho$.
\end{proof}

\begin{cor}\label{CorRep}
For every representation $\rho$ of $\hat{N}_u$ that factors through $\hat{\vartheta}^\natural$, there is a homogeneously cyclic corepresentation $U$ of $M$ such that $\rho$ is equivalent to $\pi_U^\natural$.
\end{cor}

\begin{rmk}\label{RmkGNSCorrespondence}
An essential ingredient for the proof of the quantum version of the Plancherel--Godement theorem is to see to which corepresentation the GNS-map of  $\hat{N}_u$ corresponds.  Recall that the GNS-representation of $\hat{\varphi}_{u}^\natural$ was given by the triple $(P_\gamma P_\beta \cHil, \hat{\Lambda}_u \iota_u, \hat{\pi}_u \iota_u)$, see Proposition~\ref{PropNSFWeightCast}. Since $\hat{\pi}_u \iota_u = \hat{\pi}\vert_{\hat{N}_c} \hat{\vartheta}^\natural$, we can apply Corollary~\ref{CorRep}.  We def\/ine the closed subspace
\begin{gather}\label{EqnSpaceK}
\mathcal{E} = \overline{\{ (\tilde{\omega} \otimes \iota)(W)  P_\beta P_\gamma \cHil  \mid \omega \in N_\ast^\sharp \}} = \overline{\{ ( \hat{\pi}_u \iota_u \lambda_u^\natural(\omega)  P_\beta P_\gamma \cHil  \mid \omega \in N_\ast^\sharp \}} \subseteq \cHil,
\end{gather}
where the closure is with respect to the norm in $\cHil$. It is clear that the representation of~$\hat{N}_u$ that corresponds to the restriction of~$W$ to $\mathcal{E}$ equals  $\hat{\pi}_u \iota_u$.
\end{rmk}

We use the notation $\IR(\hat{M}_u, M_1)$ to denote the irreducible representations $\pi$ of $\hat{M}_u$ such that the representation $\pi \iota_u$ is non-trivial. Under the 1-1 correspondence between $\IR(\hat{M}_u)$ and $\ICM$, see \cite{Kus}, we see from Remark~\ref{RmkAlsoUniv} and the remarks following Def\/inition \ref{DfnRep} that $\IR(\hat{M}_u, M_1)$ corresponds to $\ICMM$. Let $\IR(\hat{M}_c, M_1)$ denote the irreducible representations   of $\hat{M}_c$ such that the restriction to $\hat{N}_c$ is non-trivial. We f\/ind the following diagram of inclusions.
\begin{gather}\label{EqnRepSpaces}
\xymatrix{
&&\IR(\hat{N}_u)  & \hookleftarrow & \IR(\hat{N}_c) \\
\ICMM& \simeq&\IR(\hat{M}_u, M_1) \ar@{^{(}->}[u] & \hookleftarrow & \IR(\hat{M}_c, M_1)  \ar@{^{(}->}[u]^\simeq
}
\end{gather}
Note that the map $\IR(\hat{M}_u, M_1) \hookrightarrow \IR(\hat{N}_u)$ indeed maps into the representations of $\hat{N}_u$ that are irreducible, c.f.\ Remark~\ref{RmkIr}. Hence, also the vertical inclusion on the right hand side of (\ref{EqnRepSpaces}) preserves irreducibility.

The example in Section \ref{SectExample} shows that the inclusion $\IR(\hat{M}_c, M_1) \hookrightarrow \IR(\hat{M}_u, M_1)$ is not surjective. We brief\/ly comment on the fact that also the inclusion $\IR(\hat{M}_u, M_1) \hookrightarrow \IR(\hat{N}_u)$ is generally not surjective. This is a consequence of the fact that there are Lie groups $G$ with compact subgroup $K$ for which the there are non-unitary representations whose restriction to the bi-$K$-invariant functions forms a representation, i.e.\ a homomorphism that preserves the $\ast$-operation, whereas the representation of all $L^1$-functions on $G$ does not preserve the $\ast$. This happens for example for $SL(2, \mathbb{R})$, see \cite[Example~1.1.2 on p.~37 and p.~40]{Vog}.
 This shows that the induction argument contained in the proof of Theorem~\ref{ThmCorrespondenceC} does not work in general on the universal level.

\begin{rmk}\label{RmkIotaUIsNotInjective}
Assume the map $\iota_u: \hat{N}_u \rightarrow \hat{M}_u$ to be injective. Then by general C$^\ast$-algebra theory, it is isometric. With this additional assumption Theorem \ref{ThmCorrespondenceC} holds on the universal level. So for every $\rho \in \Rep(\hat{N}_u)$, there is a representation $\pi \in \Rep(\hat{M}_u)$ on a Hilbert space $\cHil_\pi$, such that $\rho$ is equivalent to the restriction of $\pi \iota_u$ to $(\pi \iota_u(\hat{N}_u))\cHil_\pi$.

The proof is completely analogous to the one of Theorem \ref{ThmCorrespondenceC}, where one takes the universal norm instead of the reduced norm on $N_\ast^\sharp$. The injectivity of $\iota_u$ plays an essential role at two places. First of all, the injectivity of $\iota_u$ can be used to prove positivity of the inner product~(\ref{EqnInnerProduct}), since in this case an element in $\hat{N}_u$ is positive if and only if it is positive in $\hat{M}_u$. Secondly, the universal analogue of~(\ref{EqnAdjAction}) can be recovered from the injectivity of $\iota_u$, since for $\omega \in M_\ast^\sharp$ and $\theta \in M_\ast^{\sharp, \beta}$,
\begin{gather*}
 \Vert \lambda_u^\natural (\theta^\ast \ast \omega^\ast \ast \omega \ast \theta) \Vert_u^\natural
 = \Vert \lambda_u (\theta^\ast \ast \omega^\ast \ast \omega \ast \theta) \Vert_u \\
\phantom{\Vert \lambda_u^\natural (\theta^\ast \ast \omega^\ast \ast \omega \ast \theta) \Vert_u^\natural}{}
 \leq   \Vert  \lambda_u(\omega^\ast \ast \omega) \Vert_u  \Vert \lambda_u (\theta^\ast  \ast \theta) \Vert_u
  = \Vert \lambda_u(\omega) \Vert_u \Vert \lambda_u^\natural (\theta^\ast  \ast \theta) \Vert_u^\natural.
\end{gather*}
The rest of the prove of Theorem \ref{ThmCorrespondenceC} can be copied mutatis mutandis.
\end{rmk}

\section[A quantum group analogue of the Plancherel-Godement theorem]{A quantum group analogue\\ of the Plancherel--Godement theorem}\label{SectPlancherelGodement}

Here we prove a decomposition theorem that may be considered as a locally compact quantum group version of the Plancherel--Godement theorem as can be found in \cite[Th\'eor\`eme IV.2]{Far}. The proof is dif\/ferent from the one given in \cite{Far} and follows the line of the Plancherel theorem as proved by Desmedt \cite{Des}. We show that the C$^\ast$-algebra $\hat{N}_u$ together with the weight $\varphi^\natural_u$ that we introduced and studied so far f\/it into the framework of \cite[Theorem~3.4.5]{Des}. Then we use Theorem~\ref{ThmCorrespondenceC} to translate the results in terms of corepresentations of~$M$ that admit a $M_1$-invariant vector.

Recall that we def\/ined the space $\mathcal{E}$ in (\ref{EqnSpaceK}). Let $\mathcal{L}$ be any Hilbert space and let $\mathcal{L}_0 \subseteq \mathcal{L}$ be a closed subspace. We denote the conjugate  Hilbert space of $\mathcal{L}$ by $\overline{\mathcal{L}}$. Note that the space $\mathcal{L} \otimes \overline{\mathcal{L}_0}$ can canonically be identif\/ied with the Hilbert-Schmidt operators in $B(\mathcal{L}_0, \mathcal{L})$. We denote the latter space by $B_{{\rm 2}}(\mathcal{L}_0, \mathcal{L})$. For results on direct integration, we refer to \cite{DixVNA,Lan} and~\cite{Nus}. In particular, we will use \cite[Theorem 1.10]{Lan} implicitly several times. If $A$ and $B$ are unbounded operators such that $AB$ is closable, we denote~$A \cdot B$ for the closure of~$AB$.

\begin{thm}[Plancherel--Godement]\label{ThmPG}
Let $(M, \Delta)$ be a unimodular locally compact quantum group and let $(M_1, \Delta_1)$ be a compact $($closed$)$ quantum subgroup. Let $\hat{N}$ and $\hat{N}_u$ be the von Neumann algebra and C$^\ast$-algebra as defined earlier in this section. Suppose that $\hat{N}$ is a type~I von Neumann algebra and that $\hat{N}_u$ and $\hat{M}_u$ are separable.

Then, there exists a standard measure $\mu^{M_1}$ on $\ICMM$, a $\mu^{M_1}$-measurable field of Hilbert spaces $(\cHil_U)_{U \in \ICMM}$ of which $(\cHil_U^{M_1})_{U \in \ICMM}$ forms a measurable field of subspaces, a measurable field  of self-adjoint, strictly positive operators $(D_U^{M_1})_{U \in \ICMM}$ acting on $\cHil_U^{M_1}$ and an isomorphism $\mathcal{Q}^{M_1}$ of $\mathcal{E}$ onto $\int^\oplus_{\ICMM} \cHil_U \otimes \overline{ \cHil_U^{M_1}} d\mu^{M_1}(U)$ with properties:
\begin{enumerate}\itemsep=0pt

\item\label{ItemPGOne} For $\omega \in \mathcal{I}_N$ and $\mu^{M_1}$-almost all $U \in \ICMM$, the operator $(\tilde{\omega} \otimes \iota)(U)  (D_U^{M_1})^{-1}$ is bounded and $(\tilde{\omega} \otimes \iota )(U) \cdot (D_U^{M_1})^{-1}$ is in $B_2( \cHil_U^{M_1})$.
\item\label{ItemPGTwoA} For $\omega_1, \omega_2 \in \mathcal{I}_N$, we have
\begin{gather*}
\langle \xi(\tilde{\omega}_1),  \xi(\tilde{\omega}_2) \rangle\\
\qquad{} =
\int_{\ICMM} {\rm Tr} \left( ( ( \tilde{\omega}_2 \otimes \iota)(U) \cdot (D_U^{M_1})^{-1} )^{\ast}( ( \tilde{\omega}_1 \otimes \iota)(U) \cdot (D_U^{M_1})^{-1}) \right) d\mu^{M_1}(U),
\end{gather*}
and we let $\mathcal{Q}_0^{M_1}: P_\beta P_\gamma \cHil \rightarrow \int^\oplus_{\ICMM} \cHil_U^{M_1} \otimes \overline{\cHil_U^{M_1}} d\mu^{M_1}(U)$ be the isometric extension of
\begin{gather*}
\hat{\Lambda}(\lambda(\tilde{\mathcal{I}}_N)) \rightarrow \int_{\ICMM}^\oplus B_2 (\cHil_U^{M_1})d\mu(U): \\ \xi(\tilde{\omega}) \mapsto \int^\oplus_{\ICMM} (\tilde{\omega} \otimes \iota)(U) \cdot (D_U^{M_1})^{-1} d\mu^{M_1}(U).
\end{gather*}
\item\label{ItemPGTwoB}
 $\mathcal{Q}_0^{M_1}$ intertwines $\hat{\pi}_{u} \iota_u$ and $\int_{\ICMM}^\oplus  \pi_U^\natural \otimes 1_{ \overline{ \cHil_U^{M_1}}} d\mu^{M_1}(U)$.

\item\label{ItemPGThree}
$\mathcal{Q}^{M_1}$ intertwines the restriction of $W$ to $\mathcal{E}$ with $\int^\oplus_{\ICMM}  U \otimes 1_{\overline{\cHil_U^{M_1}}}d\mu^{M_1}(U)$. Moreover, the restriction of $\mathcal{Q}^{M_1}$ to $P_\beta P_\gamma \cHil \rightarrow \int^\oplus_{\ICMM}  \cHil_U^{M_1} \otimes \overline{\cHil_U^{M_1}} d\mu^{M_1}(U)$ equals $\mathcal{Q}^{M_1}_0$.

\item\label{ItemPGFour}  Assume moreover that $\hat{M}$ is a type I von Neumann algebra and that $\hat{M}_u$ is separable. Let $\mu$, $D_U$, $\mathcal{Q}$ be defined as in {\rm \cite[{\it Theorem} 3.4.1]{Des}}. Then $\cHil_U^{M_1}$ is an invariant subspace for $\mu$-almost all $D_U$ and $\ICMM$ is a $\mu$-measurable subset of $\ICM$. If one takes:
\begin{itemize}\itemsep=0pt
\item  $\mu^{M_1}$  equal to the restriction of $\mu$ to $\ICMM$;
\item $D_U^{M_1}$ equal to the restriction of $D_U$ to $\cHil_U^{M_1}$;
\item $\mathcal{Q}^{M_1}$ the restriction of $\mathcal{Q}$ to $\mathcal{E}$.
\end{itemize}
Then, $\mu^{M_1}$, $D_U^{M_1}$, $\mathcal{Q}^{M_1}$ satisfy the properties {\rm (\ref{ItemPGOne})}--{\rm (\ref{ItemPGThree})}.
\end{enumerate}
\end{thm}

\begin{proof}
By  Propositions \ref{PropNSFWeightCast} and \ref{PropWAstLift}, $\hat{\varphi}_u^\natural$ is a proper approximate KMS-weight. Therefore, we can apply \cite[Theorem 3.4.5]{Des}, so that we obtain a measure $\mu^{M_1}$ on $\IR(\hat{N}_u)$, a measurable f\/ield of Hilbert spaces $(\mathcal{K}_\sigma^{M_1})_{\sigma \in \IR(\hat{N}_u)}$, a measurable f\/ield of representations $(\pi_\sigma)_{\sigma \in \IR(\hat{N}_u)}$, a measurable f\/ield of self-adjoint, strictly positive operators $(D_\sigma^{M_1})_{\sigma \in \IR(\hat{N}_u)}$ and an isomorphism~$\mathcal{Q}_0^{M_1}$ of $P_\gamma P_\beta \cHil$ onto $\int_{\IR(\hat{N}_u)}^\oplus \mathcal{K}_\sigma^{M_1} \otimes \overline{\mathcal{K}_\sigma^{M_1}} d \mu^{M_1}(\sigma)$ satisfying the properties of this theorem.

Let $\rho \in \IR(\hat{N}_u)$ be in the support of $\mu^{M_1}$.
We claim that $\pi_\rho$ is weakly contained in $\hat{\pi}_u \iota_u$. Suppose that this is not the case, so that there exists $x \in \hat{N}_u$ such that $\pi_\rho(x) \not = 0$ but $\hat{\pi}_u \iota_u(x)  = 0$. Let $X = \{ \sigma \in \IR(\hat{N}_u) \mid \pi_\sigma(x) \not = 0 \}$. Then $X$ is an open neighbourhood of $\rho$. Moreover, it follows form \cite[Theorem 3.4.5]{Des} that $\mathcal{Q}_0^{M_1}$ intertwines $\hat{\pi}_u \iota_u$ with $\int_{\IR(\hat{N}_u)}^\oplus \pi_\sigma \otimes \overline{1_{\mathcal{K}_\sigma}} d\mu^{M_1}(\sigma)$ from which it follows that $\mu^{M_1}(X) = 0$. This contradicts the fact that $\rho$ is in the support of $\mu^{M_1}$, so $\pi_\rho$ is weakly contained in $\hat{\pi}_u \iota_u$.

Since $\hat{\pi}_u \iota_u = \hat{\pi}\vert_{\hat{N}_c} \hat{\vartheta}^\natural$, we see that $\rho$ is in $\IR(\hat{N}_c)$, where $\IR(\hat{N}_c)$ is considered as a subset of~$\IR(\hat{N}_u)$ by the inclusion (\ref{EqnRepSpaces}).
Now we use Theorem \ref{ThmCorrespondenceC} to identify $\IR(\hat{N}_c)$ with $\IR(\hat{M}_c, M_1)$, which we consider as a subspace of~$\ICMM$ by (\ref{EqnRepSpaces}).  We consider $\mu^{M_1}$ as a measure on~$\ICMM$ by def\/ining the complement of~$\IR(\hat{N}_c)$ in~$\ICMM$ to be negligible. Let $U_\sigma \in \ICMM$ denote the corepresentation corresponding to $\sigma \in \IR(\hat{N}_c)$. So,
$
\pi_\sigma = \pi_{U_\sigma}^\natural.
$
We write $D_{U_\sigma}^{M_1}$ for $D_\sigma^{M_1}$ and set $D_U = 0$ for $U \in \ICMM $ not in the support of $\mu^{M_1}$. We denote $\cHil_U$ for the corepresentation Hilbert space of $U \in \ICMM$, and we get $\cHil_{U_\sigma}^{M_1} = \mathcal{K}_\sigma^{M_1}$.  Therefore, since the support of $\mu$ is contained in $\IR(\hat{N}_c)$, we see that $\mathcal{Q}_0^{M_1}$ is a map from $P_\beta P_\gamma \cHil \rightarrow \int^\oplus_{\ICMM} \cHil_U^{M_1} \otimes \overline{\cHil_U^{M_1}} d\mu^{M_1}(U)$.

(\ref{ItemPGOne}). It follows from  Corollary \ref{CorRep}, Remark \ref{RmkInvSubs} and the properties of $D_\sigma^{M_1}$ described in (1) of \cite[Theorem 3.4.5]{Des}, that for $\omega \in \mathcal{I}_N$, the operator $(\tilde{\omega} \otimes \iota)(U) (D_U^{M_1})^{-1}$ is bounded and its closure is Hilbert-Schmidt for $\mu^{M_1}$-almost every $U \in \ICMM$.

(\ref{ItemPGTwoA}) and (\ref{ItemPGTwoB}). We make two observations. First note that by  Proposition \ref{PropNSFWeightCast}, for $\omega \in \mathcal{I}_N$,
$
\hat{\Lambda}_u \iota_u( \lambda_u^\natural(\omega))  = \xi(\tilde{\omega}).
$
Secondly, we have proved that for every $\rho$ in the support of $\mu^{M_1}$, there is a $U \in \ICMM$, such that $\pi_\rho = \pi_U^\natural$. Then (2) of \cite[Theorem 3.4.5]{Des} yields (\ref{ItemPGTwoA}) of the present theorem. The second observation also yields (\ref{ItemPGTwoB}). Note that by Proposition \ref{PropMeas}, we see that $(\cHil_U)_{U \in \ICMM}$ is a measurable f\/ield of Hilbert spaces of which $(\cHil_U^{M_1})_{U \in \ICMM}$ forms a~measurable f\/ield of subspaces. Here we used that $\hat{M}_u$ is separable.

To prove (\ref{ItemPGThree}), we make the following observations. First of all, by Remark \ref{RmkGNSCorrespondence}, we see that $\hat{\pi}_u \iota_u = \pi_{W_{\mathcal{E}}}^\natural$, where $W_{\mathcal{E}}$ denotes the restriction of the multiplicative unitary $W$ to the Hilbert space $\mathcal{E}$. Secondly,
\[
\int_{\ICMM}^\oplus  \pi_U^\natural \otimes 1_{ \overline{ \cHil_U^{M_1}}} d\mu^{M_1}(U) =
\pi_{\int^\oplus_{\ICMM} U \otimes 1_{\overline{\cHil_U^{M_1}}} d\mu(U)}^\natural,
\]
where we use Proposition \ref{PropMeas} to infer that the direct integral on the right hand side exists.
Hence, by Remark \ref{RmkAlsoUniv} we see that $W_\mathcal{E}$ and $\int^\oplus_{\ICMM} U \otimes 1_{\overline{\cHil_U^{M_1}}} d\mu^{M_1}(U)$ are equivalent. Moreover, we def\/ine $\mathcal{Q}^{M_1}$ to be the intertwiner as constructed in the proof of Proposition \ref{PropEqRepCorep}. Then, $\mathcal{Q}^{M_1}$~satisf\/ies~(\ref{ItemPGThree}).

We now prove (\ref{ItemPGFour}). The proof of Proposition \ref{PropTBeta} shows that $\hat{\sigma}_t(P_\gamma) = P_\gamma$.
Write $D = \int_{\ICM} ^\oplus D_U d\mu(U)$. Since $D^{-2}$ is the Radon--Nikodym derivative of $\hat{\varphi}$ with respect to a trace on $\hat{M}$, see the proof of \cite[Theorem 3.4.5]{Des}, we have $D \eta \hat{M}$ and $\hat{\sigma}_t(x) = D^{-2it} x D^{2it}$. Thus,  $P_\gamma$ commutes with $D^{it}$ for all $t \in \mathbb{R}$. Since $P_\gamma \in \hat{M} \simeq \int^\oplus_{\ICM} B(\cHil_U) d\mu(U)$, we have a direct integral decomposition $P_\gamma = \int_{\ICM}^\oplus (P_\gamma)_U d\mu(U)$ and $(P_\gamma)_U$ commutes with $D_U^{it}$ for $\mu$-almost all $U \in \ICM$.

 $P_\gamma$ is the projection of $\cHil$ onto $\cHil^{M_1}$, see Proposition \ref{PropMOneInv}. A vector $v = \int^{\oplus}_{\ICM} v_U d\mu(U)$ is $M_1$-invariant for $W$ if and only if $v_U$ is $M_1$-invariant for $\mu$-almost all $U \in \ICM$, as follows directly from Def\/inition \ref{DfnMOneInv}. Hence $(P_\gamma)_U$ is the projection of $\cHil_U$ onto $\cHil_U^{M_1}$ for $\mu$-almost all $U \in \ICM$. We record three conclusions:
\begin{enumerate}\itemsep=0pt
\item[(P1)] $\cHil_U^{M_1}$ is an invariant subspace of $D_U$ for $\mu$-almost all $U \in \ICM$;
\item[(P2)]  $\ICMM$ is a $\mu$-measurable subset of $\ICM$ by \cite[Proposition II.1.1 (i)]{DixVNA};
\item[(P3)]  The image of $P_\gamma \cHil$ under $\mathcal{Q}$ equals $\int_{\ICMM}^\oplus \cHil_U^{M_1} \otimes \overline{\cHil_U} d\mu(U)$.

\end{enumerate}

For the choice of $\mu^{M_1}$, $D_U^{M_1}$ and $\mathcal{Q}^{M_1}$ made in (\ref{ItemPGFour}), properties (\ref{ItemPGOne}) and (\ref{ItemPGTwoA}) follow from the properties of $\mu$, $D_U$ and $\mathcal{Q}$ as described in \cite[Theorem 3.4.1]{Des} using Proposition \ref{PropDualSqInt}.
 By \cite[Theorem 3.4.1]{Des}, $\mathcal{Q}$ intertwines $W$ with $\int_{\ICM}^\oplus U \otimes 1_{\overline{\cHil_U}} d\mu(U)$. Using (P3) together with \cite[Theorem 3.4.5 (3)]{Des} and the fact that $P_\beta = \hat{J} P_\gamma \hat{J}$, we see that $\mathcal{Q}$ restricts to a unitary map from $P_\beta P_\gamma \cHil$ to  $\int_{\ICMM}^\oplus \cHil_U^{M_1} \otimes \overline{\cHil_U^{M_1}} d\mu(U)$. Hence $\mathcal{Q}$ restricts to a unitary map from $\mathcal{E}$ to~$\int_{\ICMM}^\oplus \cHil_U  \otimes \overline{\cHil_U^{M_1}} d\mu(U)$, which then intertwines the restriction of $W$ to $\mathcal{E}$ with $\int^\oplus_{\ICMM}   U \otimes 1_{\overline{\cHil_U^{M_1}}}d\mu(U)$. This proves (\ref{ItemPGThree}), from which (\ref{ItemPGTwoB}) follows by the construction in Re\-marks~\ref{RmkInvSubs}~and~\ref{RmkGNSCorrespondence}.
\end{proof}

\begin{rmk}
If $\Delta^\natural$ is cocommutative, then $\hat{N}$ is Abelian. By Proposition \ref{PropAbDim} we see that for all $U \in \ICM$, $\dim(\cHil_U^{M_1}) \leq 1$. Hence the operators $D_U^{M_1}$ are scalars. We may assume that $D_U^{M_1} = 1$ by replacing the measure $\mu^{M_1}$ if necessary.  In particular, we see that for a classical Gelfand pair $(G, K)$, the map $\mathcal{Q}_0^{M_1}$ is the spherical Fourier transform.
\end{rmk}

\begin{rmk}
The support of $\mu^{M_1}$ is given by $\IR(\hat{N}_c)$.
Here, $\IR(\hat{N}_c)$ is a subspace of $\ICMM$ as in (\ref{EqnRepSpaces}).  The prove can be done in exactly the same manner as  \cite[Theorem 3.4.8]{Des}, see also \cite[Proposition 8.6.8]{DixCAst}. Note that in the course of the proof of Theorem \ref{ThmPG} we have already proved that the support of $\mu^{M_1}$ is contained in $\IR(\hat{N}_c)$.
\end{rmk}

\begin{rmk}\label{RmkPlancherelAssumption}
As pointed out in remark \cite[Remark 3.3]{CasKoe}, the Theorem \ref{ThmPG} also holds if one assumes that $\hat{N}_c$ and $\hat{M}_c$ are separable, instead of $\hat{N}_u$. Moreover, note that if $\hat{M}_c$ is separable, then so is $\hat{N}_c \subseteq \hat{M}_c$. If $\hat{M}$ is a type I von Neumann algebra, then so is $\hat{N} = P_\gamma \hat{M} P_\gamma$.   In particular, if $\hat{M}$ is type I and $\hat{M}_c$ is separable, then then the result of Theorem \ref{ThmPG} holds for any closed quantum subgroup of $(M, \Delta)$.
\end{rmk}

\section[Example: $SU_q(1,1)_{{\rm ext}}$]{Example: $\boldsymbol{SU_q(1,1)_{{\rm ext}}}$}\label{SectExample}

Let $(M, \Delta)$ be the quantum group analogue of the normaliser of $SU(1,1)$ in $SL(2, \mathbb{C})$ as introduced in the operator algebraic framework in \cite{KoeKus} and further studied in \cite{GroKoeKus}. In this section we  identify the circle as a closed quantum subgroup of $(M,\Delta)$. We show that for this pair the map $\Delta^\natural$ def\/ined in (\ref{EqnShiftedCom})  is not cocommutative. Moreover, the von Neumann subalgebra $N$ of $M$ consisting of bi-invariant elements is not commutative. However, we show how the von Neumann algebras $N$ and $\hat{N}$ as def\/ined in the previous section can be equipped with a $\mathbb{Z}_2$-grading. The grading  allows us to derive similar results as for (quantum) Gelfand pairs. In particular, we make the Fourier transform explicit and show that it preserves the $\mathbb{Z}_2$-grading. Moreover, we derive product formulae for little $q$-Jacobi functions appearing as matrix coef\/f\/icients of corepresentations which admit invariant vectors.

\begin{rmk}
In \cite[Propositions B.2 and B.3]{CasKoe} it is proved that $SU_q(1,1)_{{\rm ext}}$ satisf\/ies the hypotheses of the Plancherel theorem, \cite[Theorem 3.4.1]{Des}. By Remark~\ref{RmkPlancherelAssumption} we can also apply Theorem~\ref{ThmPG}.
\end{rmk}

\begin{nt} In this section we adopt all the notational conventions made in \cite{GroKoeKus}. In particular, in this section we write $\mathcal{K}$ instead of $\mathcal{H}$ to denote the GNS-space of $(M, \Delta)$. For $z \in \mathbb{C}$, we denote $\mu(z) = (z+z^{-1})/2$. For a set $X$ and $x \in X$, we will write $\delta_x$ for the function on $X$ that equals~1 in~$x$ and~0 elsewhere. It should always be clear from the context what the domain of this function is. Recall in particular that $I_q = - q^{\mathbb{N}} \cup q^{\mathbb{Z}}$, where $\mathbb{N}$ denotes the natural numbers excluding 0.
\end{nt}

For the reader's convenience, we summarize the necessary results on the corepresentation theory of $(M, \Delta)$ from~\cite{GroKoeKus}. Let $W$ be the multiplicative unitary of $(M, \Delta)$ and recall \cite{GroKoeKus} that we have a direct integral decomposition
\begin{gather}\label{EqnWDecomposition}
W = \bigoplus_{p \in q^\mathbb{Z}} \left( \int_{[-1,1]}^\oplus
W_{p,x} dx \oplus \bigoplus_{x\in \sigma_d(\Omega_p)} W_{p,x}
\right).
\end{gather}
Here $\sigma_d(\Omega_p)$ is the discrete spectrum of the Casimir
operator \cite[Def\/inition 4.5, Theorem 4.6]{GroKoeKus} restricted to the subspace given in
\cite[Theorem 5.7]{GroKoeKus}. $W_{p,x}$ is a corepresention that is a
direct sum of at most 4 irreducible corepresentations
\cite[Propositions 5.3 and 5.4]{GroKoeKus}.  We will simply write $W = \int^\oplus W_{p,x} d(p,x)$ for the integral decomposition (\ref{EqnWDecomposition}).

The corepresentations in the continuous part of the decomposition are called principal series corepresentations, the corepresentations that appear as a direct summand are called the discrete series corepresentations. In addition the complementary series corepresentations $W_{p,x}$, $x \in (\mu(-q), -1) \cup (1, \mu(q))$, are def\/ined by analytic continuation of matrix coef\/f\/icients, see \cite[Section~10.3]{GroKoeKus}. We mention that it remains unproved that these make up all the corepresentations.

 Using the notation of \cite[Sections 10.2 and 10.3]{GroKoeKus}, an orthonormal basis for the corepresentation
Hilbert space $\mathcal{L}_{p,x}$ of the principal and complementary series $W_{p,x}$ is given by the
vectors
\begin{gather}\label{EqnBasis}
e^{\ep, \eta}_m (p,x), \qquad \ep, \eta \in \{ -, + \}, \quad m \in \mathbb{Z}.
\end{gather}
 For the discrete series corepresentation $W_{p,x}$ a subset of the vectors (\ref{EqnBasis}) gives a basis for the corepresentation space $\mathcal{L}_{p,x}$, see \cite[Proposition~5.2]{GroKoeKus}. For our purposes, it is convenient to use the following notational convention.
\begin{nt}\label{NtZeroVectors}
Let $p \in q^\mathbb{Z}$, $x \in \sigma_d(\Omega_p)$, so that $W_{p,x}$ is a discrete series corepresentation. We denote $e^{\ep, \eta}_m (p,x)$, $\ep, \eta \in \{ -, + \}, m \in \mathbb{Z}$ for the zero vector in case $e^{\ep, \eta}_m(p,x)$ is not in one of the sets def\/ined in cases 1--3 of \cite[Proposition 5.2]{GroKoeKus}. In particular, the non-zero vectors of (\ref{EqnBasis}) form an orthonormal basis of $\mathcal{L}_{p,x}$.
\end{nt}

For any $p \in q^{\mathbb{Z}}$, $x \in (\mu(-q), \mu(q)) \cup \sigma_d(\Omega_p)$, we set $f^{\ep, \eta}_m (p,x) = e^{\ep, \eta}_{m-\frac{1}{2} \chi(p)} (p,x)$, $\ep, \eta \in \{ -, + \}$, where $m \in
 \mathbb{Z}$ if $p \in q^{2 \mathbb{Z}}$ and $m \in \frac{1}{2} + \mathbb{Z}$ if $p \in q^{1+2\mathbb{Z}}$. We remind the reader that the direct integrals over of the vectors $f^{\ep, \eta}_{m}(p,x)$ are vectors in $\mathcal{K}$.
Recall that for $x \in [-1,1]$ the actions of the (unbounded) generators of the dual quantum group, see \cite[equation~(92)]{GroKoeKus}, are given by
\begin{gather}
K f^{\ep, \eta}_m(p,x)  = q^m  f^{\ep, \eta}_m(p,x), \nonumber\\
(q^{-1} - q) E  f^{\ep, \eta}_m(p,x)  = q^{-\frac{1}{2}-m} \vert 1 + \ep \eta q^{2m +1} e^{i \theta} \vert  f^{\ep, \eta}_{m+1}(p,x) , \nonumber\\
U_0^{+ -}  f^{\ep, \eta}_m(p,x)  =  \eta (-1)^{\upsilon(p)}  f^{\ep, -\eta}_{m}(p,x), \nonumber\\
U_0^{-+}  f^{\ep, \eta}_m(p,x)  = \ep \eta^{\chi(p)} (-1)^{m - \frac{1}{2} \chi(p)}  f^{-\ep, \eta}_{m}(p,x),\label{EqnKAction}
\end{gather}
where $\theta$ is such that $x = \mu(e^{i\theta})$.
Similar expressions can be obtained for the discrete and complementary series corepresentations from the expressions in \cite[Lemma~10.1 and Section~10.3]{GroKoeKus}.

\begin{rmk}\label{RmkEquiCorep}
The corepresentations appearing in (\ref{EqnWDecomposition}) are not mutually inequivalent. We give a complete list of equivalences in Proposition \ref{PropEqCorep}.
\end{rmk}

\begin{rmk}\label{RmkExtraCorep}
For every $p \in q^{\mathbb{Z}}$, $x \in \mu(-q^{2\mathbb{Z}+1}p \cup q^{2\mathbb{Z}+1}p)$, one can def\/ine a corepresentation $W_{p,x}$ by def\/ining the action of the generators of $\hat{M}$ by means of \cite[Lemma~10.1]{GroKoeKus}. Since the actions of the generators of $W_{pr,x}$ are equal for any  $r \in q^{4\mathbb{Z}}$,   these corepresentations are all equivalent.  Using \cite[Proposition 5.2]{GroKoeKus} one can check that every such corepresentation is equivalent to at least one corepresentation in the decomposition (\ref{EqnWDecomposition}). Since every discrete series corepresentation is inf\/inite dimensional, it occurs inf\/initely many times in the decomposition (\ref{EqnWDecomposition}) by the Plancherel theorem \cite[Theorem 3.4.1]{Des}, or see Proposition \ref{PropEqCorep} below for a direct proof. Therefore, we see that
\begin{gather*}
W \simeq \bigoplus_{p \in q^\mathbb{Z}} \left( \int_{[-1,1]}^\oplus
W_{p,x} dx \oplus \bigoplus_{x \in \mu(-q^{2\mathbb{Z}+1 }p \cup q^{ 2\mathbb{Z}+1}p)} W_{p,x}
\right).
\end{gather*}
\end{rmk}

We will also need the following expressions for the matrix coef\/f\/icients.
By \cite[Lemma 10.9 and Section 10.3]{GroKoeKus}, for $p \in q^{\mathbb{Z}}$, $x \in (\mu(-q),\mu(q))$,
\begin{gather}\label{EqnWSliceII}
(\iota \otimes \omega_{f^{\ep,\eta}_m(p,x), f^{\ep',\eta'}_{m'}(p,x)})\left(W_{p,x} \right) f_{m_0, p_0, t_0}   \\
= C\big(\ep\eta x;m'-\tfrac{1}{2}\chi(p),\ep',\eta';\ep\ep'\vert p_0\vert q^{-m-m'},p_0,m-m'\big) \delta_{\sgn(p_0),\eta\eta'} f_{m_0 - m + m',\ep\ep' \vert p_0\vert q^{-m-m'} , t_0}. \nonumber
\end{gather}
We refer to \cite{GroKoeKus} for the precise def\/initions of the function  $C(\cdot)$ in terms of basic hypergeometric series. (\ref{EqnWSliceII}) also holds for the discrete series corepresentations.

\subsection*{The diagonal subgroup}

Let $\mathbb{T}$ denote the circle group and let $(L^\infty (\mathbb{T}), \Delta_{\mathbb{T}})$ be the usual locally compact quantum group associated with $\mathbb{T}$ with dual quantum group $(L^\infty(\mathbb{Z}), \Delta_{\mathbb{Z}})$. We identify $(L^\infty (\mathbb{T}), \Delta_{\mathbb{T}})$ as a closed quantum subgroup of $(M, \Delta)$. Recall from \cite[Def\/inition~4.3]{GroKoeKus} that the spectrum of the opera\-tor~$K$  equals $0 \cup q^{\frac{1}{2}\mathbb{Z}}$.
\begin{dfn}\label{DfnPiHat}
 We def\/ine a normal, injective $\ast$-homomorphism $\hat{\pi}: L^\infty(\mathbb{Z}) \rightarrow \hat{M}$ by setting $\hat{\pi}(\delta_k) = \delta_{q^{\frac{k}{2}}}(K)$.
\end{dfn}

Note that $\hat{\pi}$ preserves the comultiplication, since
\begin{gather*}
  (\hat{\pi} \otimes \hat{\pi}) \Delta_{\mathbb{Z}}(\delta_k)  =  (\hat{\pi} \otimes \hat{\pi}) \left( \sum_{l \in \mathbb{Z}} \delta_l \otimes \delta_{k - l} \right) \\
\phantom{(\hat{\pi} \otimes \hat{\pi}) \Delta_{\mathbb{Z}}(\delta_k)}{}
 =   \sum_{l \in \mathbb{Z}} \delta_{q^{\frac{l}{2}}}(K) \otimes \delta_{q^{\frac{k-l}{2}}}(K) =  \delta_{q^{\frac{k}{2}}}(K \otimes K) =  \delta_{q^{\frac{k}{2}}}(\hat{\Delta}(K)).
\end{gather*}
Therefore, $\hat{\pi}$ identif\/ies $(\mathbb{T}, \Delta_{\mathbb{T}})$ as a closed quantum subgroup of $(M, \Delta)$.
Furthermore, $\hat{\pi}$ induces a morphism $\pi$ between the universal quantum groups $(M_u, \Delta_u)$ and $(C(\mathbb{T}), \Delta_{\mathbb{T}})$, where  here with slight abuse of notation $\Delta_{\mathbb{T}}$ is restricted to a map $C(\mathbb{T}) \rightarrow C(\mathbb{T} \times \mathbb{T})$.

\subsection*{Spherical corepresentations}

We compute the actions $\gamma$ and $\beta$ of left and right translation and determine which of the corepresentations found in \cite{GroKoeKus} admit a $L^\infty(\mathbb{T})$-invariant vector.

\begin{prop}\label{PropActions}
For all $p \in q^{\mathbb{Z}}$, $x \in (\mu(-q),\mu(q)) \cup \sigma_d(\Omega_p)$,
 \begin{gather}
\beta\big(\iota \otimes \omega_{ f_m^{\ep,\eta}(p,x), f_{m'}^{\ep',\eta'} (p,x)}\big)(W_{p,x})   =
\big(\iota \otimes \omega_{ f_m^{\ep,\eta}(p,x), f_{m'}^{\ep',\eta'} (p,x)}\big)(W_{p,x}) \otimes \zeta^{2m}, \label{EqnBetaAction} \\
\gamma \big(\iota \otimes \omega_{ f_m^{\ep,\eta}(p,x), f_{m'}^{\ep',\eta'} (p,x)}\big)(W_{p,x})   =
\big(\iota \otimes \omega_{ f_m^{\ep,\eta}(p,x), f_{m'}^{\ep',\eta'} (p,x)}\big)(W_{p,x}) \otimes \zeta^{2m'}. \label{EqnGammaAction}
\end{gather}
Here, $\zeta$ is the identity function on the complex unit circle $\mathbb{T}$.
\end{prop}

\begin{proof}
Note that $\hat{\pi}$ has a direct integral decomposition $\hat{\pi} =  \int^\oplus
\hat{\pi}_{p,x} d(p,x)$. Here $\hat{\pi}_{p,x}: L^\infty(\mathbb{Z}) \rightarrow B(\mathcal{L}_{p,x})$ is determined by $\hat{\pi}_{p,x}(\delta_k) = \delta_{q^{\frac{k}{2}}}(K)$, where the action of $K$ on the representation space is given by (\ref{EqnKAction}) and similarly for the discrete series corepresentations.

By the def\/inition of $\beta$, see (\ref{EqnBeta}), and the decomposition (\ref{EqnWDecomposition}), we see that
\begin{gather*}
 \int^\oplus (\beta \otimes \iota)(W_{p,x}) d(p,x) =
 \int^\oplus (W_{p,x})_{13}   (\iota \otimes \hat{\pi}_{p,x}) (W_{\mathbb{T}})_{23}   d(p,x).
\end{gather*}
By (\ref{EqnKAction}) and the def\/inition of $\hat{\pi}$,  for almost all pairs $(p,x)$ in the decomposition (\ref{EqnWDecomposition}),
\begin{gather}\label{EqnActionDec}
(\beta \otimes \iota)(W_{p,x}) = (W_{p,x})_{13} (\iota \otimes \hat{\pi}_{p,x})(W_{\mathbb{T}})_{23}.
\end{gather}
This proves (\ref{EqnBetaAction}) for the discrete series corepresentations.
It follows from (\ref{EqnWSliceII}) and the fact that the function $C$ given there is analytic on a neighbourhood of $(\mu(-q), \mu(q))$, see \cite[Section~10.3]{GroKoeKus}, that for every $p \in q^{\mathbb{Z}}$,
\begin{gather}\label{EqnLHS}
(\mu(-q), \mu(q)) \rightarrow M \otimes L^\infty(\mathbb{T}): \ x \mapsto \beta \left( \big(\iota \otimes \omega_{ f_m^{\ep,\eta}(p,x) , f_{m'}^{\ep',\eta'} (p,x) }\big)(W_{p,x}) \right) ,
\end{gather}
extends to an analytic function on a neighbourhood of  $(\mu(-q), \mu(q))$. Here, we use the fact that~$\beta$ is normal and the fact that a function is $\sigma$-weak analytic if and only if it is analytic with respect to the norm \cite[Result 1.2]{KusAut}. Similarly, it follows that for $p \in q^{\mathbb{Z}}$,
\begin{gather}\label{EqnRHS}
(\mu(-q), \mu(q))\rightarrow M \otimes L^\infty(\mathbb{T}): \ x \mapsto  \big(\iota \otimes \omega_{ f_m^{\ep,\eta}(p,x), f_{m'}^{\ep',\eta'}(p,x) }\big)(W_{p,x}) \otimes \zeta^{2m},
\end{gather}
extends to an analytic function on a neighbourhood of  $(\mu(-q), \mu(q))$. Note that for all $p \in q^{\mathbb{Z}}$, $x \in (\mu(-q), \mu(q))$,
\begin{gather*}
 \big(\iota \otimes \omega_{ f_m^{\ep,\eta}(p,x), f_{m'}^{\ep',\eta'}(p,x) }\big)(W_{p,x}) \otimes \zeta^{2m}\\
  \qquad {}= \big(\iota \otimes \iota \otimes \omega_{ f_m^{\ep,\eta}(p,x) , f_{m'}^{\ep',\eta'} (p,x) }\big) \left( (W_{p,x})_{13} (\iota \otimes \hat{\pi}_{p,x})(W_{\mathbb{T}})_{23} \right).
\end{gather*}
Now, (\ref{EqnActionDec}) yields that (\ref{EqnLHS}) and (\ref{EqnRHS}) agree on a dense subset of $[-1,1]$. Since (\ref{EqnLHS}) and (\ref{EqnRHS}) have analytic extensions, they are equal for every $x \in (\mu(-q), \mu(q))$.
The proof of (\ref{EqnGammaAction}) is similar.
\end{proof}

\begin{rmk}
Note that in the preceding proof, we cannot apply \cite[Theorem~IV.8.25]{TakI} directly to the map $\hat{\pi}$ on the von Neumann algebraic level since $L^\infty(\mathbb{T})$ is not separable. Therefore, we have def\/ined the representations $\hat{\pi}_{p,x}$ explicitly.
\end{rmk}

Recall that we def\/ined $L^\infty(\mathbb{T})$-invariant vectors in Def\/inition~\ref{DfnMOneInv}.

\begin{cor}\label{CorHom}
Let $p \in q^{\mathbb{Z}}$,  $x \in ( \mu(-q), \mu(q))\cup \sigma_d(\Omega_p)$.
\begin{enumerate}\itemsep=0pt
\item\label{ItemCorZero} If $p \in q^{2 \mathbb{Z}}$, then the space of $L^\infty(\mathbb{T})$-invariant vectors of  $W_{p,x}$ is spanned by $f_0^{\ep, \eta}(p,x)$, with $\ep, \eta \in \{ +, -\}$.
\item\label{ItemCorI} If $p \in q^{2 \mathbb{Z}}$, $x \in \sigma_d(\Omega_p)$.  The space of $L^\infty(\mathbb{T})$-invariant vectors of $W_{p,x}$ is one dimensional.
\item\label{ItemCorII}  If $p \in q^{1 + 2 \mathbb{Z}}$, then $W_{p,x}$ has no $L^\infty(\mathbb{T})$-invariant vectors.
\item\label{ItemHomIV} $P_\gamma = \delta_1(K) = \hat{\pi}(\delta_0)$ and hence the von Neumann algebras $N$ and $\hat{N}$ that are constructed in Section~{\rm \ref{SectHomogeneousSpaces}} are given by
\begin{gather*}
N  = \overline{\left\{
(\iota \otimes \omega_{\delta_1(K) v, \delta_1(K) w})(W) \mid v,w \in   \mathcal{K}
\right\} }^{\sigma\text{\rm -strong-}\ast} ,\qquad
\hat{N}  =   \delta_1(K) \hat{M} \delta_1(K).
\end{gather*}
\end{enumerate}
\end{cor}

\begin{proof}
From the considerations in Sections \ref{SectHomogeneousSpaces} and \ref{SectSphericalCoreps}, (\ref{ItemCorZero}), (\ref{ItemCorII}) and (\ref{ItemHomIV}) follow. For  (\ref{ItemCorI}) consider $\lambda \in -q^{2 \mathbb{Z} +1} \cup q^{2 \mathbb{Z} +1}$ be such that $x = \mu(\lambda)$ and $\vert \lambda \vert \geq 1$. Put $j', l' \in \mathbb{Z}$ by setting $\vert \lambda \vert = q^{1-2j'} =   q^{1+2l '}$. In particular, $l' < j'$. For case (i) of \cite[Proposition 5.2]{GroKoeKus}, $f^{-,+}_0(p,x)$ is zero if and only if $0 > l'$ if and only if $\vert \lambda \vert \geq 1$, which is true by assumption. Similarly, $f^{+,-}_0(p,x) = 0$. Hence the only $L^\infty(\mathbb{T})$-invariant vector is $f^{+,+}_0(p,x)$. Cases (ii) and (iii) of \cite[Proposition 5.2]{GroKoeKus} follow in exactly the same manner.
\end{proof}

\begin{rmk}
Note that the space of $L^\infty(\mathbb{T})$-invariant vectors for the irreducible components of $W_{p,x}$ is not necessarily of dimension $\leq 1$. For example $W_{p,x}$, $p \in q^{2\mathbb{Z}}$, $x \in [-1,1] \backslash \{0\}$ splits as a~sum of 2 irreducible corepresentations of which the $L^\infty(\mathbb{T})$-invariant vectors form a 2-dimensional vector space, see \cite[Section 10.2]{GroKoeKus}. This implies that $\Delta^\natural$ is not cocommutative and we can not use (\ref{EqnProdForm}) directly to obtain product formulae. However, we will def\/ine gradings on the spaces $N$ and $\hat{N}$ which still allows us to f\/ind such formulae.
\end{rmk}

In Corollary~\ref{CorHom} we determined the corepresentations appearing in the decomposition that admit a $L^\infty(\mathbb{T})$-invariant vector. These corepresentations are not mutually inequivalent. Here, we give a list of the equivalences. We only consider the spherical corepresentations and consider the principal, discrete as well as the complementary series.  Only the principal and discrete series are important to determine the spaces $N$, $\hat{N}$ and the spherical Fourier transform, see Theo\-rem~\ref{ThmPG}. Nevertheless, the complementary series is still important, since  the product formulae we derive later still hold for the spherical matrix elements of the complementary series.

Motivated by \cite[Lemma 10.11]{GroKoeKus}, we introduce the following basis vectors. For $ p \in q^{2\mathbb{Z}}$, $x \in (\mu(-q),\mu(q))  $, set
\begin{gather*}
g_m^{1,+}(p,x) = \frac{1}{2}\sqrt{2}( f^{+,+}_m(p,x)+i^{\chi(p)}f^{-,-}_m(p,x)) ,\nonumber\\
g_m^{1,-}(p,x) = \frac{1}{2}\sqrt{2}(f^{+,-}_m(p,x)-i^{\chi(p)}f^{-,+}_m(p,x)),\nonumber\\
g_m^{2,+}(p,x) = \frac{1}{2}\sqrt{2}(f^{+,+}_m(p,x)-i^{\chi(p)}f^{-,-}_m(p,x)),\nonumber\\
g_m^{2,-}(p,x) = \frac{1}{2}\sqrt{2}( f^{+,-}_m(p,x)+i^{\chi(p)}f^{-,+}_m(p,x)).
\end{gather*}
For every $j \in \{1,2\}$, $p \in q^{2\mathbb{Z}}$, $x \in (\mu(-q),\mu(q))$, the vectors $g_m^{j, \sigma}(p,x)$, $\sigma \in \{ +,-\}$, $m \in \mathbb{Z}$ form an orthonormal basis for $\mathcal{L}_{p,x}^j$, the corepresentation space of one of the  summands of $W_{p,x}$, see \cite[equation~(95)]{GroKoeKus}.
For $  p \in q^{2\mathbb{Z}}$, $x \in \sigma_p(\Omega_d) $, set
\begin{gather}
g_m^{1,+}(p,x) =   f^{+,+}_m(p,x)+i^{\chi(p)}f^{-,-}_m(p,x) ,\nonumber\\
g_m^{1,-}(p,x) = f^{+,-}_m(p,x)-i^{\chi(p)}f^{-,+}_m(p,x),\nonumber\\
g_m^{2,+}(p,x) = f^{+,+}_m(p,x)-i^{\chi(p)}f^{-,-}_m(p,x),\nonumber\\
g_m^{2,-}(p,x) = f^{+,-}_m(p,x)+i^{\chi(p)}f^{-,+}_m(p,x).\label{EqnGVectorsDisc}
\end{gather}
Recall that we made the convention that $f^{\ep,\eta}_m(p,x) = 0$ in case $f^{\ep,\eta}_m(p,x)$ is not in the basis given in \cite[Proposition~5.2]{GroKoeKus}. Hence, for $x \in \sigma_d(\Omega_p)$, we see from \cite[Proposition~5.2]{GroKoeKus} that the vectors def\/ined in (\ref{EqnGVectorsDisc}) are dependent and any of them is equal to a vector of the form $f^{\ep, \eta}_m(p,x)$ for some $\ep, \eta \in \{ +,-\}$ modulo a phase factor.

 Now, we determine which of the discrete, principal and complementary series corepresentations are equivalent. By considering the action of the Casimir operator as in the proof of \cite[Proposition~B.1]{CasKoe}, it is clear that any two corepresentations that fall within a dif\/ferent series are inequivalent. We restrict ourselves to the spherical corepresentations.
\begin{prop}\label{PropEqCorep}\qquad
\begin{enumerate}\itemsep=0pt
 \item\label{ItemEqCorep} Let $x, x' \in (\mu(-q),\mu(q)) \backslash \{ 0 \}$, $p,p'\in q^{2\mathbb{Z}}$, $j, j' \in \{1,2\}$. $W_{p,x}^j \simeq W_{p',x'}^{j'}$ if and only if either $x = \pm x'$, $j=j'$, $p/p' \in q^{4 \mathbb{Z}}$ or $x = \pm x'$, $j \not= j'$, $p/p' \in q^{2+4 \mathbb{Z}}$.
\item  Let $p,p'\in q^{2\mathbb{Z}}$, $j, j', k, k' \in \{1,2\}$. $W_{p,0}^{j,k} \simeq W_{p',0}^{j',k'}$ if and only if either $j=j'$, $k=k'$, $p/p'\in q^{4\mathbb{Z}}$ or $j\not = j'$, $k=k'$, $p/p'\in q^{2+4\mathbb{Z}}$.
\item Let $x, x' \in \mu(-q^{\mathbb{Z}} \cup q^{\mathbb{Z}})$, $p,p'\in q^{2\mathbb{Z}}$. $W_{p,x} \simeq W_{p', x'}$ if and only if $\vert x \vert = \vert x' \vert$.
\end{enumerate}
\end{prop}
\begin{proof} The proposition follows from a careful comparison of the action of the generators, see~(\ref{EqnKAction}) for the principal series, \cite[Proposition 5.2 and Lemma 10.1]{GroKoeKus} for the discrete series and \cite[Section 10.3]{GroKoeKus} for the complementary series. We prove (\ref{ItemEqCorep}). By considering the Casimir operator \cite[Def\/inition 4.5]{GroKoeKus} one sees that if an irreducible component of $W_{p,x}$ is equivalent to an irreducible component of $W_{p',x'}$, then $\vert x \vert = \vert x' \vert$.

In case $x = x'$, an intertwiner must send $g_m^{j,\pm}(p,x)$ to a non-zero scalar multiple of $g_m^{j',\pm}(p',x)$ as follows by considering the actions of $K$ and $E$. Writing out the actions of $U^{+-}_0$ and $U^{-+}_0$ one sees that there exists such an intertwiner only in the following two cases:
\begin{enumerate}\itemsep=0pt
\item[(i)] $ p/p' \in q^{4 \mathbb{Z}}$, $j =j',$ for which it sends $g_m^{j,\pm}(p,x)$ to  $g_m^{j,\pm}(p',x)$;
\item[(ii)] $ p/p' \in q^{2+ 4 \mathbb{Z}}$, $j \not = j'$, in which case it sends $g_m^{j,\pm}(p,x)$ to  $\pm g_m^{j',\pm}(p',x)$.
\end{enumerate}

 Similarly, in case $x = -x'$, an intertwiner must send $g_m^{j,\pm}(p,x)$ to  $g_m^{j',\mp}(p',-x)$ as follows from the actions of $K$ and $E$. The actions of $U^{+-}_0$ and $U^{-+}_0$ show that this is only possible if
\begin{enumerate}\itemsep=0pt
\item[(i)] $ p/p' \in q^{4 \mathbb{Z}}$, $j =j'$ for which it sends $g_m^{j,\pm}(p,x)$ to  $\pm g_m^{j,\mp}(p',-x)$;
\item[(ii)]  $ p/p' \in q^{2+ 4 \mathbb{Z}}$, $j \not = j'$, in which case it sends $g_m^{j,\pm}(p,x)$ to  $g_m^{j',\pm}(p',-x)$.
\end{enumerate}
This proves (\ref{ItemEqCorep}), the other cases follows similarly.
\end{proof}

\begin{figure}[h!]
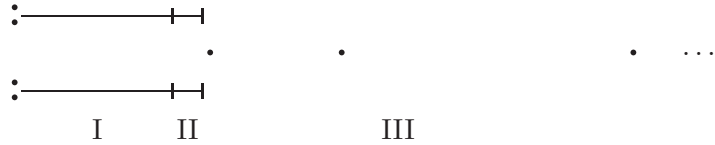
\centering
$$
\xy
  (0,5)*{}="a0";
  (0,-5)*{}="b0";
  (20,5)*{}="a1";
  (20,-5)*{}="b1";
  (24,5)*{}="a2";
  (24,-5)*{}="b2";
  "a0"; "a1" **\dir{-};
  "b0"; "b1" **\dir{-};
  "b1"; "b2" **\dir{-};
  "a1"; "a2" **\dir{-};
  (20,-6)*{};(20,-4)*{} **\dir{-};
  (20,4)*{};(20,6)*{} **\dir{-};
  (24,-6)*{};(24,-4)*{} **\dir{-};
  (24,4)*{};(24,6)*{} **\dir{-};
  (10,-10)*{{\rm I }};
  (22,-10)*{{\rm II }};
  (50,-10)*{{\rm III }};
  "a0"+(-1,1)*{ \mh{\cdot}};
  "a0"+(-1,-1)*{ \mh{\cdot}};
  "b0"+(-1,1)*{\mh{\cdot}};
  "b0"+(-1,-1)*{ \mh{\cdot}};
  (25,0)*{\mh{\cdot}};
  (42.5,0)*{\mh{\cdot}};
  (81.25,0)*{\mh{\cdot}};
  (90,0)*{   \ldots };
\endxy
$$
\caption{Known part of $\ICMM$ for $SU_q(1,1)_{\rm{ext}}$ with the diagonal subgroup, c.f.~(\ref{EqnICMM}). Part I: Principal spherical series.  Part II: Complementary spherical series. Part III: Discrete spherical series.  }\label{PicSpectrum}
\end{figure}

Summarizing Corollary \ref{CorHom} and Proposition \ref{PropEqCorep}, we f\/ind that~$\ICMM$, the space of equiva\-lence classes of irreducible spherical corepresentations is partly given by
\begin{gather}\label{EqnICMM}
\begin{array}{@{}r@{\,}l@{\,}ll}
\ICMM &\supseteq  &(0,1,1) \cup (0,1,2) \cup (0,2,1) \cup (0,2,2) \cup ( (0,1] \times \{ 1,2 \}) & \textrm{(principal)} \\
  & &\cup  \quad \mu(q^{2\mathbb{N}+1}) & {\rm (discrete)} \\
  & & \cup \quad (1, \mu(q)) \times \{ 1,2 \} & {\rm (complementary)}
\end{array}\hspace*{-10mm}
\end{gather}
 Here, we identify the points  $(0,1,1)$, $(0,1,2)$, $(0,2,1)$, $(0,2,2)$ with the respective irreducible corepresentations $W_{1,0}^{1,1}$, $W_{1,0}^{1,2}$, $W_{1,0}^{2,1}$, $W_{1,0}^{2,2}$, see \cite[Proposition~10.14~(ii)]{GroKoeKus}. We let a point $(x,j) \in (0,1] \times \{ 1,2 \}$ correspond to $W_{1,x}^j$, see \cite[Proposition~10.12]{GroKoeKus}. The points $x \in  \mu(q^{2\mathbb{N}+1})$ corresponds to $W_{1,x}$, see \cite[Proposition 5.2]{GroKoeKus}. We emphasize that it is not known if the corepresentations described in~\cite{GroKoeKus} are all the corepresentations, therefore we do not know if this description completely describes~$\ICMM$. For the von Neumann algebras $N$ and $\hat{N}$ as well as the spherical Fourier transform, only the principal and discrete (spherical) series matter. These are fully identif\/ied within (\ref{EqnICMM}).
 For completeness, we illustrate the right hand side of~(\ref{EqnICMM}) by means of Fig.~\ref{PicSpectrum}.

\subsection*{The von Neumann algebra $\boldsymbol{N}$}

The next step is to make $N$ and $\hat{N}$ more explicit and meanwhile def\/ine a grading on these spaces. Therefore, we will f\/irst f\/ind an alternative formula for the mapping $T_\beta T_\gamma$ which is  convenient for computations. This is also going to play a role when we   derive product formulae. Recall that the spectrum of $K$ is given by $q^{\frac{1}{2}\mathbb{Z}}$. For $k, l \in \mathbb{Z}$, set
\[
M_{k,l}   = \overline{\left\{
(\iota \otimes \omega_{\delta_{q^{\frac{1}{2}(k-l)}}(K) v, \delta_{q^{\frac{1}{2}(k+l)}}(K) w})(W) \mid v,w \in   \mathcal{K}
\right\} }^{\sigma\text{-strong-}\ast}.
\]
By (\ref{EqnKAction})   and (\ref{EqnWSliceII}), the spaces $M_{k,l}$ have mutually trivial intersection. Moreover, $M$ equals the $\sigma$-strong-$\ast$ closure of the direct sum of vector spaces $\oplus_{k,l \in \mathbb{Z}} M_{k,l}$. Note that $N = M_{0,0}$.

\begin{dfn}
For $m \in \mathbb{Z}$, $p, t \in -q^{\mathbb{Z}} \cup q^{\mathbb{Z}}$, let $P_{m,p,t}$ be the orthogonal rank one projection of $\mathcal{K}$ onto the space spanned by the vector $f_{m,p,t}$. Here we def\/ine $f_{m,p,t}$ to be the zero vector if either $p \not \in I_q$ or $t \not \in I_q$.  Def\/ine mappings $T^+$ and $T^-$ def\/ined by
\[
T^{\pm}: \ M \rightarrow M: \  x \mapsto \sum_{m \in \mathbb{Z}, p,t \in I_q} P_{m,\pm p,t} x P_{m,p,t},
\]
where the sum converges in the strong topology.
\end{dfn}
\begin{prop}\label{PropTpTm}
$T^+ + T^- = T_\beta T_\gamma$.
\end{prop}
\begin{proof}
 For a normal functional $\omega = \sum_{i \in I} \omega_{\xi_i, \eta_i} \in M_\ast$, with $\sum_{i \in I} \Vert \xi_i \Vert^2, \sum_{i \in I} \Vert \eta_i \Vert^2 < \infty$, the linear map
\[
M \rightarrow \mathbb{C}: \  x \mapsto \omega T^\pm(x) = \sum_{i \in I} \sum_{m \in \mathbb{Z}, p, t \in I_q} \langle x P_{m,p,t} \xi_i, P_{m,\pm p,t} \eta_i \rangle,
\]
is normal. Hence, the maps $T^\pm$ are normal. Moreover, for $x \in M_{k,l}$, using (\ref{EqnWSliceII}),
\begin{gather}\label{EqnTpTm}
T^+(x)  + T^-(x) = \left\{
\begin{array}{ll}
x, & \textrm{if } k=l=0, \\
0, & \textrm{otherwise}.
\end{array}
\right.
\end{gather}
For $x \in M_{k,l}$ we see by using (\ref{EqnBetaAction}) and (\ref{EqnGammaAction}) that $T_\beta T_\gamma (x)$ also equals the right hand side of~(\ref{EqnTpTm}).
Since also $T_\beta T_\gamma$ is normal, this proves that $T^+ + T^- = T_\beta T_\gamma$.
\end{proof}

Let $u_0\in B( L^2(I_q)) $ be the partial isometry determined by $u_0: \delta_p \mapsto \delta_{-p}$, and recall that $u$ was def\/ined as $u = 1 \otimes u_0 \in B( L^2(\mathbb{Z}) \otimes L^2(I_q))$.
By \cite[Lemma~2.4]{KoeKus} $M \simeq L^\infty(\mathbb{T}) \otimes B(L^2(I_q))$. Proposition~\ref{PropTpTm} yields,
\[
N = T_\beta T_\gamma (M) = (T^+ + T^-)(M),
\]
which is isomorphic to the von Neumann subalgebra of $M$ generated by $1 \otimes L^\infty(I_q)$ and~$u$. Therefore, introduce the following identif\/ications.

\begin{dfn}
We identify $N$ with the von Neumann algebra acting on $L^2(I_q)$ being generated by $L^\infty(I_q)$ and $u_0$.
 We split $N$ as a direct sum of vector spaces
\begin{gather*}
N = N_+ \oplus N_-,\quad \textrm{where} \quad N_+ = L^\infty(I_q), \quad \textrm{and}\quad  N_- =  L^\infty(I_q)u_0 =  L^\infty(I_q \cap (-1,1))u_0.
\end{gather*}
This turns $N$ into a $\mathbb{Z}_2$-graded algebra.
\end{dfn}

  We f\/ind that $\varphi^\natural$, i.e.\ the restriction of the Haar weight $\varphi$ \cite[Def\/inition~4.1]{KoeKus} to~$N$, equals the measure given by:
\[
\varphi^\natural(f) = \sum_{p_0 \in I_q} f(p_0) p_0^2, \qquad f \in N^+ = L^\infty(I_q)^+.
\]
For $f \in \mathfrak{m}_{\varphi^\natural}^+$ we f\/ind that $u_0f \in \mathfrak{m}_{\varphi^\natural}$ and it follows by \cite[Def\/inition~4.1]{KoeKus} that $\varphi^\natural(u_0 f) = 0$, so that $\Lambda(\nphi \cap N_+)$ and $\Lambda(\nphi \cap N_-)$  are orthogonal spaces. The discussion so far allows us to make the following identif\/ications.

\begin{dfn}
Identify the closure of $\Lambda(\nphi \cap N_+)$ with $L^2(I_q)$. We will write $L^2(N_+)$ for~$L^2(I_q)$ to indicate explicitly that $L^2(I_q)$ is considered as part of the GNS-space of $\varphi^\natural$.
 Similarly, we identify the closure of $\Lambda(\nphi \cap N_-)$ with $L^2(I_q \cap (-1,1))$, by identifying $f u_0 \in N_-$, where  $f \in  L^\infty(I_q \cap (-1,1)) \cap L^2(I_q \cap (-1,1))$ with $f \in L^2(I_q \cap (-1,1))$.  We write short hand $L^2(N_-)$ for  $L^2(I_q \cap (-1,1))$. We write $L^2(N) = L^2(N_+) \oplus L^2(N_-)$. We emphasize that here the spaces~$L^2(I_q)$ (respectively  $L^2(I_q \cap (-1,1))$) should be understood with respect to the integral given by the weighted sum $\sum_{p_0 \in I_q }( \:\: \cdot \:\: )   p_0^2$ (respectively $\sum_{p_0 \in I_q \cap (-1,1)}( \:\: \cdot \:\: )   p_0^2$).
\end{dfn}

\subsection*{The von Neumann algebra $\boldsymbol{\hat{N}}$}
We now turn our attention to the von Neumann algebra $\hat{N}$ as def\/ined in Section~\ref{SectHomogeneousSpaces}. Considered as a subalgebra of $\hat{M}$, it inherits the  $\mathbb{Z}_2$-grading def\/ined in \cite[Def\/inition~4.7]{GroKoeKus}.

\begin{dfn}
We set
\[
\hat{N}_+ = \hat{N} \cap \hat{M}_+, \qquad \hat{N}_- = \hat{N} \cap \hat{M}_-,
\]
see \cite[Def\/inition 4.7]{GroKoeKus} for $\hat{M}_+$ and $\hat{M}_-$. This turns $\hat{N}$ into a $\mathbb{Z}_2$-graded algebra.
\end{dfn}

\begin{prop}
We have an isomorphism of von Neumann algebras:
\begin{gather}\label{EqnIsoVNA}
 \hat{N} \simeq \int_{[0,1]}^\oplus M_2(\mathbb{C}) \oplus M_2(\mathbb{C}) dx \oplus \bigoplus_{x \in \mu(q^{2\mathbb{N}+1})} \mathbb{C}.
\end{gather}
The isomorphism is determined by the map
\begin{gather}\label{EqnVNAIso}
(\tilde{\omega} \otimes \iota)(W) \mapsto \int_{[0,1]}^\oplus (\tilde{\omega} \otimes \iota)(W_{1,x}^1) \oplus (\tilde{\omega} \otimes \iota)(W_{1,x}^2)  dx \oplus \bigoplus_{x \in \mu(q^{2\mathbb{N}+1})}  (\tilde{\omega} \otimes \iota)(W_{1,x}) .
\end{gather}
Under this isomorphism $\hat{N}_+$ corresponds to the direct integrals over matrices whose entries vanish off the diagonal. $\hat{N}_-$ corresponds to the direct integrals over matrices with values vanishing on the diagonal.
\end{prop}
\begin{proof}
 It follows from Corollary \ref{CorHom} that the map def\/ined in (\ref{EqnVNAIso}) indeed maps into the proper matrix algebras. Let $P \in \hat{M} \cap \hat{M}'$ be the projection as in the proof of \cite[Proposition~B.2]{CasKoe}. By considering the Casimir operator, one sees that $P$ is the spectral projection on the interval $[-1,1]$ of the Casimir operator.  In \cite[Proposition~B.2]{CasKoe}, it is shown that
\[
 P \hat{M}P \simeq \int^{\oplus}_{x \in [0,1]} \hat{M}_x dx,
\]
where $\hat{M}_x$ is generated by $\{ (\omega \otimes \iota)(W_x) \mid \omega \in M_\ast \}$ and $W_x = (\oplus_{p \in q^{\mathbb{Z}}} W_{p,x}) \oplus  (\oplus_{p \in q^{\mathbb{Z}}} W_{p,-x})$.  By similar techniques as in \cite[Proposition~B.2]{CasKoe}, one can show that
\begin{gather}\label{EqnDirectInt}
  \hat{M}  \simeq \int^{\oplus}_{x \in [0,1]} \hat{M}_x dx \oplus  \bigoplus_{x \in  \sigma_d(\Omega) \cap (1, \infty) } \hat{M}_x,
\end{gather}
where  $\hat{M}_x, x \in \sigma_d(\Omega ) \cap (1, \infty)$ is generated by $\{ (\omega \otimes \iota)(W_x) \mid \omega \in M_\ast \}$ and
\[
W_x = (\bigoplus_{p \in q^{\mathbb{Z}}, {\rm s.t. } x \in \sigma_d(\Omega_p)} W_{p,x}) \oplus  (\bigoplus_{p \in q^{\mathbb{Z}}, {\rm s.t. } -x \in \sigma_d(\Omega_p)} W_{p,-x}) .
\]
 Since $P_\gamma \in \hat{M}$, it has a direct integral decomposition as in (\ref{EqnDirectInt}), i.e.
 \[
  P_\gamma \hat{M} P_\gamma      \simeq \int^{\oplus}_{x \in [0,1]} (P_\gamma)_x \hat{M}_x (P_\gamma)_x dx  \oplus  \bigoplus_{x \in  \sigma_d(\Omega) \cap (1, \infty) }(P_\gamma)_x  \hat{M}_x (P_\gamma)_x .
\]
Similarly, $K$ decomposes with respect to the direct integral decomposition (\ref{EqnDirectInt}) as a direct integral over a f\/ield $(K_x)_{x \in [0,1] \cup (\sigma_d(\Omega) \cap (1, \infty))}$ and by Proposition \ref{PropActions}, we see that $(T_\beta T_\gamma \otimes \iota)(W_x) = (1 \otimes \delta_1(K_x)) W_{x} (1 \otimes \delta_1(K_x))  = (1 \otimes (P_\gamma)_x) W_x (1 \otimes (P_\gamma)_x)$. Hence, $(P_\gamma)_x \hat{M}_x (P_\gamma)_x$ is generated by $\{ (\tilde{\omega} \otimes \iota)(W_x) \mid \omega \in M_\ast \}$. Due to Corollary \ref{CorHom} and Proposition \ref{PropEqCorep}, the latter von Neumann algebra is isomorphic to $M_2(\mathbb{C})\oplus M_2(\mathbb{C})$ in case $x \in (0,1]$ and to $\mathbb{C}$ in case $x \in \sigma_d(\Omega) \cap (1, \infty)$. This proves (\ref{EqnIsoVNA}). That (\ref{EqnVNAIso}) gives the isomorphism follows directly from this proof. The claim on the gradings follows from \cite[Proposition 4.8]{GroKoeKus}.
\end{proof}

\begin{rmk}
The von Neumann algebra $\hat{M}_0$ in the previous proof is isomorphic to 4 copies of~$\mathbb{C}$ as follows from \cite[Proposition 10.13]{GroKoeKus}. Since this does not matter for the integral decomposition (\ref{EqnIsoVNA}), and to avoid some redundant extra notation we have not treated the point $x =0$ separately.
\end{rmk}

We claim that $\hat{\varphi}^\natural$, i.e.\ the restriction of the dual left Haar weight to $\hat{N}$, is a trace. Indeed, it follows from \cite[Section 5]{CasKoe} and \cite[Proposition 3.5.5]{Des} that for any $p \in q^{\mathbb{Z}}$, $x \in [-1,1] \cup \sigma_d(\Omega_p)$, there is a constant $c(p,x)$ such that
\begin{gather*}
D_{p,x} f_{m}^{\epsilon, \eta}(p,x) = p^{1/2} q^m c(p,x) f_m^{\ep, \eta}(p,x).
\end{gather*}
Here $D = \int^\oplus D_{p,x} d(p,x)$, with short hand notation $D_{p,x} = D_{W_{p,x}}$, is the Duf\/lo--Moore operator arising from the Plancherel theorem, see \cite[Theorem 3.4.1]{Des} and also Theorem \ref{ThmPG}.
Since the eigenvalues of $D_{p,x}$ are independent of the signs $\ep, \eta$, we see that $\int^\oplus D_{p,x}^{L^\infty(\mathbb{T})} d(p,x) = P_\gamma \int^\oplus D_{p,x}  d(p,x) P_\gamma$ is central in $\hat{N}$.
Since $\int^\oplus (D_{p,x}^{L^\infty(\mathbb{T})} )^{-2} d(p,x)$ is the Radon--Nikodym derivative of $\hat{\varphi}^\natural$ with respect to a trace, see the proof of \cite[Theorem 3.4.5]{Des}, we see that $\hat{\varphi}^\natural$ is a trace. Therefore, under the identif\/ication (\ref{EqnIsoVNA}), we see that
\begin{gather}\label{EqnDualWeight}
 \hat{\varphi}^\natural = \int^\oplus_{[0,1]} \left( {\rm Tr}_{M_2(\mathbb{C})} \oplus  {\rm Tr}_{M_2(\mathbb{C})}  \right) d(x)  dx \oplus  \bigoplus_{x \in \mu(q^{2\mathbb{N}+1})} d(x) {\rm Tr}_{\mathbb{C}},
\end{gather}
where $d$ is a scalar valued function, which we leave undetermined.
\begin{rmk}
It follows from \cite[Theorem 3.4.1 (6)]{Des} that the function $c(1,x)$ depends on the Plancherel measure. In case the Plancherel measure is choosen as in (\ref{EqnWDecomposition}), one has $d(x) = c(1,x)^{-2}$. For the discrete series these constants are computed in \cite[Section 3.5]{Des}. For the principal series, the exact values cannot be found in the literature. We will not derive them here, since this is beyond the scope of our example.
\end{rmk}
 We identify the GNS-space of $\varphi^\natural$ using the isomorphism (\ref{EqnIsoVNA}). That is, the GNS-space is given by
\begin{gather*}
L^2(\hat{N}) = \int_{[0,1]}^\oplus M_2(\mathbb{C}) \oplus M_2(\mathbb{C})  dx \oplus \bigoplus_{x \in \mu(q^{2\mathbb{N}+1})} \mathbb{C},
\end{gather*}
where the direct integral and direct sums are taken as Hilbert spaces (as opposed to the direct integrals and sums of von Neumann algebras given in (\ref{EqnIsoVNA})), where the inner product comes from the traces $d(x) ( {\rm Tr}_{M_2(\mathbb{C})} \oplus  {\rm Tr}_{M_2(\mathbb{C})}   ) $ for the integral part and $d(x) {\rm Tr}_{\mathbb{C}}$ for the direct summands.

It follows from (\ref{EqnIsoVNA}) and (\ref{EqnDualWeight}) that $\hat{\Lambda}(\mathfrak{n}_{\hat{\varphi}} \cap \hat{N}_+)$ and $\hat{\Lambda}(\mathfrak{n}_{\hat{\varphi}} \cap \hat{N}_-)$  are orthogonal spaces, we denote their closures interpreted within $L^2(\hat{N})$ by $L^2(\hat{N}_+)$ and $L^2(\hat{N}_-)$. They consist of direct integrals of diagonal matrices and of\/f-diagonal matrices, respectively.

\subsection*{The spherical Fourier transform}
Here, we determine the spherical Fourier transform, i.e.\ we determine the map $\mathcal{Q}_0^{L^\infty(\mathbb{T})}$ for Theorem \ref{ThmPG}.
 In particular, we are interested in the kernels of the integral transformations that appear in this transform, since we will need them later. Using the identif\/ications of the GNS-spaces for $\varphi^\natural$ and $\hat{\varphi}^\natural$ with $L^2(N)$ and $L^2(\hat{N})$, we may consider $\mathcal{Q}_0^{L^\infty(\mathbb{T})}$ as a map from~$L^2(N)$ to~$L^2(\hat{N})$ and use the short hand notation $\mathcal{F}_2: L^2(N) \rightarrow L^2(\hat{N})$ for this map.

If $f \in L^1(I_q) \cap L^\infty(I_q) \subseteq N_+$, then there is a   functional $f\cdot \varphi^\natural \in N_\ast$ given by  $(f\cdot \varphi^\natural )(x) = \varphi^\natural(x f)$, $x \in N$. Then, $\xi(f\cdot \varphi^\natural ) = \Lambda(f)$ by def\/inition of $\xi$, see Section \ref{SectPreliminaries}. Thus, for such $f$, we f\/ind by Theorem \ref{ThmPG} that,
\begin{gather}\label{EqnFourierTr}
\mathcal{F}_2: \Lambda(f) \mapsto    ((f \cdot \varphi^\natural)^\sim  \otimes \iota) \left( \int^\oplus_{[0,1]} W_{1,x}  dx  \oplus \bigoplus_{x \in \mu(q^{2\mathbb{N}+1})}      W_{1,x}  \right).
\end{gather}
Here, the right hand side indeed is an element of $L^2(\hat{N})$ by Corollary \ref{CorHom}. Note that we choose the inner product on $L^2(\hat{N})$ to be the direct integral over the traces  $d(x) ( {\rm Tr}_{M_2(\mathbb{C})} \oplus  {\rm Tr}_{M_2(\mathbb{C})}   ) $ for the integral part and $d(x) {\rm Tr}_{\mathbb{C}}$ for the direct summands. Hence, the Duf\/lo--Moore operators which are direct integrals of scalar multiples of the identity are contained in the inner product by means of the function $d$. Hence, they do not appear in (\ref{EqnFourierTr}).
Similarly, for $f \in L^1(I_q) \cap L^\infty(I_q)$, so that $f u_0 \in N_-$,
\begin{gather}\label{EqnFourierTrII}
\mathcal{F}_2: \Lambda(f u_0) \mapsto    ((f u_0 \cdot \varphi^\natural)^\sim  \otimes \iota) \left( \int^\oplus_{[0,1]} W_{1,x}  dx  \oplus \bigoplus_{x \in \mu(q^{2\mathbb{N}+1})}      W_{1,x}  \right).
\end{gather}

We see from (\ref{EqnFourierTr}) that the spherical Fourier transformation   is in fact a combination of integral transformations with the spherical matrix elements of the corepresentations $W_{1,x}$ in its kernel. The next step is to make these kernels explicit.

For any $m_0 \in \mathbb{Z}$, $t_0 \in I_q$, $x \in [-1,1] \cup \mu(-q^{2\mathbb{N}+1} \cup q^{2\mathbb{N}+1})$, $j \in \{1,2\}$, $\sigma, \tau \in \{ +, - \}$, put
\[
 K_j^{\sigma, \tau}(p_0; x)  = \langle (\iota \otimes \omega_{g_0^{j,\sigma}(1,x), g_0^{j,\tau}(1,x) })(W_{1,x})f_{m_0, p_0, t_0},f_{m_0, \sigma \tau p_0, t_0} \rangle.
 \]
We emphasize that $ K_j^{\sigma, \tau}(p_0; x)$ for $x$ not in $\sigma_d(\Omega_1)$ is def\/ined by Remark \ref{RmkExtraCorep}.
This expression is independent of $m_0$ and $t_0$. Due to Proposition \ref{PropEqCorep}, we have the following symmetry
\begin{gather}\label{EqnSymmLqJ}
 K_j^{\sigma, \tau}(p_0; x) = \sigma \tau K_j^{-\sigma, -\tau}(p_0; -x).
\end{gather}
From (\ref{EqnWSliceII}), (\ref{EqnFourierTr}) and (\ref{EqnFourierTrII}), we see that we get a graded Fourier transform $\mathcal{F}_{2,+} \oplus \mathcal{F}_{2,-}: L^2(N_+) \oplus L^2(N_-) \rightarrow L^2(\hat{N}_+) \oplus L^2(\hat{N}_-)$, which is def\/ined by sending $f \oplus g u_0 \in (L^2(I_q) \cap L^1(I_q)) \oplus (L^2(I_q \cap (-1,1)) \cap L^1(I_q \cap (-1,1))) u_0 \subseteq L^2(N_+) \oplus L^2(N_-)$ to the matrix valued function on ${\rm IC}(M, L^\infty(\mathbb{T}))$ determined by sending $ (x,j)  \in  [0, 1]   \times \{0,1\} $ to
\begin{gather}\label{EqnFT}
\sum_{p_0 \in I_q}
\left(
\begin{array}{ll}
 K_j^{+, +}(p_0; x) f(p_0) &  K_j^{-, +}(p_0; x) g(p_0) \\
 K_j^{+, -}(p_0; x) g(p_0) &  K_j^{-, -}(p_0; x) f(p_0)
\end{array}
\right)
p_0^2,
\end{gather}
and $x \in \mu(q^{2 \mathbb{Z}+1})$ to
\[
\sum_{p_0 \in I_q}  K_0^{+, +}(p_0; x) f(p_0) p_0^2 = \sum_{p_0 \in I_q}  K_1^{+, +}(p_0; x) f(p_0) p_0^2.
\]
Note that for $x =0$, the matrix appearing in (\ref{EqnFT}) is actually a direct sum of two matrix blocks after a basis transformation, since $W_{1,0}$ splits as a direct sum of four irreducible corepresentations, see \cite[Proposition 10.13]{GroKoeKus}. By Theorem \ref{ThmPG}, this map is unitary.

For completeness, we give the analogous result for the inverse Fourier transform $\mathcal{F}_2^{-1}:   L^2(\hat{N}_+) \oplus L^2(\hat{N}_-)\rightarrow L^2(N_+) \oplus L^2(N_-)$. Since $\mathcal{Q}_0^{L^\infty(\mathbb{T})}$ is a restriction of the Plancherel transformation $\mathcal{Q}$, see Theorem \ref{ThmPG}, we see that  $\mathcal{F}_2^{-1}$ can be considered as the restriction of $\mathcal{Q}^{-1}$.
The transform  $\mathcal{Q}^{-1}$ is described in \cite{Cas}  on the operator algebraic level, see in particular \cite[Example~3.13]{Cas}. See also \cite{DaeFou} for the algebraic counterpart. More explicitly, the spherical inverse Fourier transform is determined by
\[
\mathcal{F}_2^{-1}: \ f \mapsto (\iota \otimes (f \cdot \hat{\varphi}))(W^\ast),
\]
where $f \in \hat{N} \cap \mathfrak{n}_{\hat{\varphi}}$ is such that there is a normal functional on $\hat{M}$, denoted by $(f \cdot \hat{\varphi})$, which is determined by $(f \cdot \hat{\varphi})(x) = \hat{\varphi}(x f)$, $x \in \mathfrak{n}_{\hat{\varphi}}^\ast$. By the decomposition of $W$ (\ref{EqnWDecomposition}), we f\/ind the following theorem.

\begin{thm}
For   $\sigma, \tau \in \{+, -\}$, let
$
f_1^{\sigma, \tau},   f_2^{\sigma, \tau}  \in L^1([0,1]) \cap L^2([0,1]),
$
where $L^1([0,1])$ and $L^2([0,1])$ should be understood with respect to the integral $\int_{[0,1]} \: d(x)dx$. Let
$
f_d \in L^1(\mu(q^{2\mathbb{N}+1}) ) \cap  L^2(\mu(q^{2\mathbb{N}+1}) ),
$
where  $L^1(\mu(q^{2\mathbb{N}+1}) )$ and $L^2(\mu(q^{2\mathbb{N}+1}) )$ should be understood with respect to the integral $\sum_{x \in \mu(q^{2\mathbb{N}+1}) }\:  d(x)$.
 Define the function $f \in L^2(\hat{N})$ by
\[
f(x) = \left\{
\begin{array}{ll}
\left(
\begin{array}{ll}
f_1^{+,+}(x) & f_1^{-,+}(x) \\
f_1^{+,-}(x) & f_1^{-,-}(x)
\end{array}
\right)
\oplus
\left(
\begin{array}{ll}
f_2^{+,+}(x) & f_2^{-,+}(x) \\
f_2^{+,-}(x) & f_2^{-,-}(x)
\end{array}
\right), \quad
&
x \in [0,1], \\
f_d(x), & x \in \mu(q^{2\mathbb{N}+1}).
\end{array}
\right.
\]
Then,
\begin{gather}
  \int_{[0,1]} \left( f_1^{+,+}(x) \overline{K_1^{+,+}(p_0;x)}  + f_1^{-,-}(x) \overline{K_1^{-,-}(p_0;x)} \right) d(x) dx
  \nonumber\\
 \qquad{}+   \int_{[0,1]}  \left( f_2^{+,+}(x) \overline{K_2^{+,+}(p_0;x)}  +  f_2^{-,-}(x) \overline{K_2^{-,-}(p_0;x)} \right) d(x) dx\nonumber\\
  \qquad{}
+   \sum_{x \in \mu(q^{2\mathbb{N}+1})} \left( f_d(x) \overline{ K_1^{+,+}(p_0; x)} d(x) \right)\nonumber \\
\qquad{} \bigoplus  \left( \int_{[0,1]} \left( f_1^{-,+}(x) \overline{K_1^{+,-}(p_0;x)}  + f_1^{+,-}(x) \overline{K_1^{-,+}(p_0;x)} \right) d(x) dx \right.
\nonumber\\ \left.\qquad{}+    \int_{[0,1]} \left( f_2^{-,+}(x) \overline{K_2^{+,-}(p_0;x)}  +  f_2^{+,-}(x) \overline{K_2^{-,+}(p_0;x)} \right)  d(x) dx  \right)  u_0\label{EqnInverseTransform}
\end{gather}
exists for every $p_0 \in I_q$. Moreover, \eqref{EqnInverseTransform} considered as a direct sum of functions in $p_0$ forms an element of $L^2(N_+) \oplus L^2(N_-)$. This mapping extends to a unitary map $L^2(\hat{N}_+) \oplus L^2(\hat{N}_-) \rightarrow L^2(N_+) \oplus L^2(N_-)$, which is inverse to $\mathcal{F}_2$.
\end{thm}

We explicitly state the formulae for the kernels $K_j^{\sigma, \tau}(p_0;x)$ for $x \in [-1,1] \cup \mu(-q^{2\mathbb{N}+1} \cup q^{2\mathbb{N}+1})$, which can be expressed in terms of little $q$-Jacobi functions.  Using the notation of \cite[Section~9]{GroKoeKus} for $S(\cdot)$ and $A(\cdot)$, we f\/ind
\begin{gather*}
K_1^{+, +}(p_0; x)  = S(-\lambda, p_0, p_0, 0) \times \left\{
\begin{array}{ll}
1, & p_0 < 0, \vspace{1mm}\\
\dfrac{A(\lambda, 1, 0, +, +)}{A(\lambda, 1, 0, -, -)}, & p_0 > 0,
\end{array}
\right. \\
K_2^{+, +}(p_0; x)   = S(-\lambda, p_0, p_0, 0) \times \left\{
\begin{array}{ll}
1, & p_0 < 0, \vspace{1mm}\\
-\dfrac{A(\lambda, 1, 0, +, +)}{A(\lambda, 1, 0, -, -)} ,& p_0 > 0,
\end{array}
\right. \\
K_1^{+, -}(p_0; x)   = -S(\lambda, -p_0, p_0, 0) \times \left\{
\begin{array}{ll}
\dfrac{A(\lambda, 1, 0, +, +)}{A(-\lambda, 1, 0, +, -)}, & p_0 < 0, \vspace{1mm}\\
-\dfrac{A(\lambda, 1, 0, +, +)}{A(-\lambda, 1, 0, -, +)}, & p_0 > 0,
\end{array}
\right. \\
K_2^{+, -}(p_0; x)    = -S(\lambda, -p_0, p_0, 0) \times \left\{
\begin{array}{ll}
\dfrac{A(\lambda, 1, 0, +, +)}{A(-\lambda, 1, 0, +, -)}, & p_0 < 0, \vspace{1mm}\\
\dfrac{A(\lambda, 1, 0, +, +)}{A(-\lambda, 1, 0, -, +)}, & p_0 > 0.
\end{array}
\right. \\
\end{gather*}
The fractions of the functions $A(\cdot)$ are phase factors.
Here, $\lambda \in \mathbb{T}$ such that $\mu(\lambda) = x$. And, for $\lambda \in \mathbb{C} \backslash \{0\}$, simplifying \cite[Lemma~9.1~(ii)]{GroKoeKus},
\begin{gather*}
S(\pm \lambda, \mp p_0, p_0, 0) =   \vert p_0 \vert^2 \nu(p_0)^2 c_q^2 \sqrt{(\pm \kappa(p_0), - \kappa(p_0); q^2)_\infty } (\mp q^2; q^2)_\infty    \\
\qquad{} \times\frac{(q^2, -q^2/\kappa(p_0), \lambda q^2, 1/\lambda, -q/\lambda ; q^2)_\infty}{(\sgn(p_0) 1/\lambda, \sgn(p_0) \lambda q^2, \pm q/\lambda; q^2)_\infty } \rFs{2}{1}{-q/\lambda, -\lambda q}{-q^2}{q^2, -q^2/\kappa(p_0)}.
\end{gather*}
The other kernels occurring in (\ref{EqnFT}) can be expressed explicitly by means of the symmetry relation~(\ref{EqnSymmLqJ}). The matrix coef\/f\/icients $K^{\sigma, \tau}_j(p_0, x)$ are special types of little $q$-Jacobi functions, see also \cite[Appendix B.5]{GroKoeKus} and references given there.

\subsection*{Product formulae for little $\boldsymbol{q}$-Jacobi functions}

Let $N_{\ast, +}$ be the space of normal functionals in $N_\ast$ which are zero
on $N_-$. Let $N_{\ast, -}$ be the space of normal functionals in $N_\ast$ which
are zero on $N_+$. Note that $\varphi^\natural$ is a trace on $N$.
Therefore, every functional in $N_{\ast, +}$ is given by $f \cdot
\varphi^\natural$, where $f$ is a function on $I_q$ so that $\sum_{p_0 \in I_q} f(p_0) p_0^2 < \infty$.
Every
functional in $N_{\ast, -}$ is given by $f u_0 \cdot
\varphi^\natural$, where $f$ is a function on $I_q \cap (-1,1)$ so that $\sum_{p_0 \in I_q \cap (-1,1)} f(p_0) p_0^2 < \infty$.  So $N_\ast = N_{\ast, +} \oplus N_{\ast, -}$ as
vector spaces.

 For $p \in I_q$, recall that $\delta_p$ denotes the function on $I_q$ whose value is 1 in $p$ and 0
elsewhere. We write $\delta_{p,+}$ for   $\delta_p \cdot \varphi^\natural \in  N_{\ast,+}$. For
$p \in I_q $, we write $\delta_{p,-}$ for the
functional $\delta_p u_0 \cdot \varphi^\natural \in N_{\ast,-}$. Note that only for $p \in I_q \cap (-1,1)$, this functional is non-zero, but it is convenient to keep this notation.

\begin{rmk}\label{RmkFunctionals}
Note that if one identif\/ies $N$ with the subalgebra $(1 \otimes L^\infty(I_q) \otimes 1 )\cup (1 \otimes L^\infty(I_q)u_0 \otimes 1)$   of $M$ acting on the GNS-space $\mathcal{K}$, then $\delta_{p,+} = \omega_{f_{m_0, p, t_0}, f_{m_0, p, t_0}}$ and $\delta_{p,-} = \omega_{f_{m_0, -p, t_0}, f_{m_0, p, t_0}}$, where $m_0 \in \mathbb{Z}$, $t_0 \in I_q$. This functional is independent of $m_0$ and $t_0$ if consi\-de\-red as a functional on $N$.
\end{rmk}
Using the fact that the Fourier transform preserves the $\mathbb{Z}_2$-gradings on $N$ and $\hat{N}$, we see that for $p_1, p_2 \in I_q$:
\begin{alignat*}{3}
& (\delta_{p_1,+} \otimes \delta_{p_2,+}) \Delta^\natural  \in N_{\ast,+},  \qquad &&
(\delta_{p_1,-} \otimes \delta_{p_2,-}) \Delta^\natural  \in N_{\ast,+}, & \\
& (\delta_{p_1,+} \otimes \delta_{p_2,-}) \Delta^\natural  \in N_{\ast,-}, \qquad &&
(\delta_{p_1,-} \otimes \delta_{p_2,+}) \Delta^\natural  \in N_{\ast,-}.&
\end{alignat*}
 Hence, we see that for $f, g \in L^\infty(I_q)$, there exist constants $A_{p_0}(p_1, p_2)$, $B_{p_0}(p_1, p_2)$,  $p_0, p_1, p_2 \in I_q$ and $C_{p_0}(p_1, p_2)$, $D_{p_0}(p_1, p_2)$,  $p_0 \in I_q \cap (-1,1)$, $p_1, p_2 \in I_q$  such that for any $p_1, p_2 \in I_q$ the four sums
\begin{alignat*}{3}
& \sum_{p_0 \in I_q} \vert A_{p_0}(p_1, p_2)\vert p_0^2, \qquad && \sum_{p_0 \in I_q} \vert B_{p_0}(p_1, p_2)\vert p_0^2 , & \\
& \sum_{p_0 \in I_q\cap (-1,1)} \vert C_{p_0}(p_1, p_2)\vert p_0^2,\qquad && \sum_{p_0 \in I_q\cap (-1,1)} \vert D_{p_0}(p_1, p_2)\vert p_0^2,&
\end{alignat*}
are f\/inite and
\begin{gather}
(\delta_{p_1,+} \otimes \delta_{p_2,+}) \Delta^\natural (f)   =   \sum_{p_0
\in I_q} A_{p_0}(p_1, p_2) f(p_0) p_0^2,\label{EqnPrI}  \\
(\delta_{p_1,-} \otimes \delta_{p_2,-}) \Delta^\natural (f)   =   \sum_{p_0
\in I_q} B_{p_0}(p_1, p_2) f(p_0) p_0^2, \\
(\delta_{p_1,+} \otimes \delta_{p_2,-}) \Delta^\natural (f u_0)   =   \sum_{p_0
\in I_q \cap (-1,1)} C_{p_0}(p_1, p_2) f(p_0) p_0^2, \\
(\delta_{p_1,-} \otimes \delta_{p_2,+}) \Delta^\natural (f u_0)   =   \sum_{p_0
\in I_q \cap (-1,1)} D_{p_0}(p_1, p_2) f(p_0) p_0^2.\label{EqnPrIV}
\end{gather}
Putting   $f = K_j^{\sigma, \tau}(p_0 ;x)$ for any $x \in [-1,0) \cup (0,1] \cup \mu(-q^{2\mathbb{N}+1} \cup q^{2\mathbb{N}+1} )$, $j \in \{1,2\}$, this
yields a~product formula for $_2 \varphi_1$-series. For $p_1, p_2 \in I_q$ and $p_3, p_4 \in I_q \cap (-1,1)$,
\begin{gather}
K_j^{+, +}(p_1 ;x) K_j^{+, +}(p_2 ;x)   =   \sum_{p_0 \in
I_q} A_{p_0}(p_1, p_2)K_j^{+, +}(p_0 ;x) p_0^2, \label{EqnProdI} \\
K_j^{+, -}(p_3 ;x) K_j^{-, +}(p_4 ;x)    =   \sum_{p_0 \in
I_q} B_{p_0}(p_3, p_4) K_j^{-, -}(p_0 ;x) p_0^2, \label{EqnProdII} \\
K_j^{+, +}(p_1 ;x) K_j^{-, +}(p_3 ;x)   =   \sum_{p_0 \in
I_q \cap (-1,1)} C_{p_0}(p_1, p_3) K^{-,+}_{j}(p_0; x) p_0^2,  \\
K_j^{-, +}(p_3 ;x) K_j^{-, -}(p_1 ;x)    =   \sum_{p_0 \in
I_q \cap (-1,1)} D_{p_0}(p_3, p_1) K_j^{-, +}(p_0 ;x) p_0^2. \label{EqnProdIV}
\end{gather}
\begin{rmk}
Note that the gradings on $N$ and $\hat{N}$ make that the left hand sides of (\ref{EqnProdI})--(\ref{EqnProdIV}) consists of a single product of two $_2 \varphi_1$-functions.
\end{rmk}

\begin{rmk}
If $x \in \mu(-q^{2\mathbb{N}+1} \cup q^{2\mathbb{N}+1})$, so that $W_{1,x}$ is a discrete series corepresentation, then equations (\ref{EqnProdII})--(\ref{EqnProdIV}) are trivial, i.e.\ they equate~0 to~0.
 If in addition $x <0$, then~(\ref{EqnProdI}) is also trivial.
\end{rmk}
\begin{rmk}
From (\ref{EqnPrI})--(\ref{EqnPrIV}), we can also get formulae for the products
\begin{alignat*}{3}
& K_j^{-, -}(p_1 ;x) K_j^{-, -}(p_2 ;x), \qquad && K_j^{+, -}(p_3 ;x) K_j^{-, +}(p_4 ;x),&\\
& K_j^{+, +}(p_1 ;x) K_j^{-, +}(p_3 ;x),\qquad && K_j^{-, +}(p_3 ;x) K_j^{-, -}(p_1 ;x),&
\end{alignat*}
where  $j \in \{1, 2\}$, $p_1, p_2 \in I_q$ and $p_3, p_4 \in I_q \cap (-1,1)$. However, using the symmetry (\ref{EqnSymmLqJ}), these formulae  are already contained in (\ref{EqnProdI})--(\ref{EqnProdIV}).
\end{rmk}

In the remainder of this section we determine the coef\/f\/icient functions $A$, $B$, $C$, and $D$. We explicitly show how to f\/ind $A$, the other coef\/f\/icients can be found by the same method. Let $f = \delta_{p_0}$, so that the right hand side of (\ref{EqnPrI}) is equal to $A_{p_0}(p_1, p_2) p_0^2$. To determine the left hand side, note that $\delta_{p_0} = \delta_{\sgn(p_0)}(e) \delta_{p_0^{-2}} (\gamma^\ast \gamma)$. Recall that $\Delta(e) = e \otimes e \in N \otimes N$, so that $ \Delta( \delta_{\sgn(p_0)}(e) ) = \delta_{\sgn(p_0)}(e \otimes e) \in N \otimes N$. Hence,
\[
\Delta^\natural (\delta_{\sgn(p_0)}(e)) = \Delta( \delta_{\sgn(p_0)}(e) ) =\delta_1(e) \otimes \delta_{\sgn(p_0) }(e) + \delta_{-1}(e) \otimes \delta_{-\sgn(p_0) }(e) \in N \otimes N.
\]
Note that by the relation $\Delta(x) = W^\ast (1 \otimes x) W$, $x \in M$, we f\/ind that $\Delta( \delta_{p_0^{-2}} (\gamma^\ast \gamma) )$ equals the projection onto the closure of
\[
\textrm{span} \left\{ W^\ast f_{m', p' , t'} \otimes f_{m, p_0, t} \mid m, m' \in \mathbb{Z}, p', t, t' \in I_q \right\}.
\]
This projection is given by the formula
\[
\mathcal{K} \otimes \mathcal{K}  \rightarrow \mathcal{K} \otimes \mathcal{K}:  v \mapsto \sum_{m, m' \in \mathbb{Z}, p', t, t'\in I_q} \langle v , W^\ast f_{m', p', t'} \otimes f_{m, p_0, t} \rangle \: W^\ast f_{m', p', t'} \otimes f_{m, p_0, t},
\]
where the sum is norm convergent. Note that for $x \in N$, $\Delta(x) \in M^\gamma \otimes M^\beta$, so that $\Delta^\natural(x) = (\iota \otimes T_\beta T_\gamma) \Delta(x) = (\iota \otimes (T^+ +T^-)) \Delta(x)$. Since $\Delta^\natural ( \delta_{p_0^{-2}} (\gamma^\ast \gamma) ) \in N \otimes N$, we f\/ind that for any $m_1, m_2 \in \mathbb{Z}$ and $t_1, t_2 \in I_q$,
\begin{gather}
 (\delta_{p_1} \otimes \delta_{p_2}) \Delta^\natural ( \delta_{p_0^{-2}} (\gamma^\ast \gamma) )  \nonumber\\
 \qquad{} =
 \langle (\iota \otimes (T^+ + T^-))( \Delta  ( \delta_{p_0^{-2}} (\gamma^\ast \gamma) )) f_{m_1, p_1, t_1} \otimes f_{m_2, p_2, t_2}, f_{m_1, p_1, t_1} \otimes f_{m_2, p_2, t_2} \rangle \nonumber\\
 \qquad{} =
 \langle   \Delta  ( \delta_{p_0^{-2}} (\gamma^\ast \gamma) ) f_{m_1, p_1, t_1} \otimes f_{m_2, p_2, t_2}, f_{m_1, p_1, t_1} \otimes f_{m_2, p_2, t_2} \rangle
 \nonumber\\
 \qquad{} =
\sum_{m', m \in \mathbb{Z}, p', t , t' \in I_q} \vert \langle f_{m_1, p_1, t_1} \otimes f_{m_2, p_2, t_2} , W^\ast f_{m', p', t'} \otimes f_{m, p_0, t} \rangle \vert^2 ,\label{EqnGammaComp}
\end{gather}
where the equations follow from Remark \ref{RmkFunctionals} and Proposition \ref{PropTpTm}, the def\/inition of~$T^+$ and~$T^-$ and the discussion above.
Recall the functions $a_p(x,y)$, $p, x, y \in I_q$ from \cite[Def\/inition~6.2]{GroKoeKus}. From \cite[Propositions~4.5 and~4.10]{KoeKus}, we see that (\ref{EqnGammaComp}) equals
\begin{gather*}
  \sum_{t \in I_q}  \left(\frac{t}{t_2}\right)^2 a_{t}( \sgn(p_0 p_2 t_2) p_1 q^{m_2} t, t_2)^2 a_{p_0}(p_1, p_2)^2 \nonumber\\
  \qquad{} =
  \sum_{t \in I_q}  a_{t_2}( \sgn(t_2) \vert p_1 \vert q^{m_2} t, t)^2 a_{p_0}(p_1, p_2)^2 =
a_{p_0}(p_1, p_2)^2.\label{EqnACoefficient}
\end{gather*}
Here, the f\/irst equality follows from \cite[equation~(24)]{GroKoeKus} and the fact that by def\/inition $a_{p_0}(p_1, p_2)^2$ $ = 0$ if $\sgn(p_0p_2) \not = \sgn(p_1)$. The last equality follows from \cite[Proposition~6.3]{GroKoeKus}. Hence, using \cite[Section~IX.4, equation~(4)]{TakII} in the second equality,
\begin{gather*}
  A_{p_0}(p_1, p_2) p_0^2 =   (\delta_{p_1} \otimes \delta_{p_2}) \Delta^\natural \big(\delta_{\sgn(p_0)}(e)  \delta_{p_0^{-2}} (\gamma^\ast \gamma) \big)   \\
  \phantom{A_{p_0}(p_1, p_2) p_0^2}{} =
(\delta_{p_1} \otimes \delta_{p_2})(\iota \otimes T_\gamma)\big( \delta_{\sgn(p_0)}( \Delta(e))  \delta_{p_0^{-2}} \Delta(\gamma^\ast \gamma)\big) \\
\phantom{A_{p_0}(p_1, p_2) p_0^2}{}
=
(\delta_{p_1} \otimes \delta_{p_2})\delta_{\sgn(p_0)}( \Delta(e)) (\iota \otimes T_\gamma)\big(  \delta_{p_0^{-2}} \Delta(\gamma^\ast \gamma)\big)  \\
\phantom{A_{p_0}(p_1, p_2) p_0^2}{}
=
 \delta_{\sgn(p_0), \sgn(p_1p_2)} a_{p_0}(p_1, p_2)^2 = a_{p_0}(p_1, p_2)^2.
\end{gather*}
Recall \cite[Def\/inition 3.1]{KoeKus} that by def\/inition $a_{p_0}(p_1, p_2) = 0$ if $\delta_{\sgn(p_0), \sgn(p_1p_2)} = 0$, so indeed the last equality follows.
By a similar computation we get,
\begin{gather*}
 B_{p_0}(p_1, p_2) p_0^2   =  a_{p_0}(p_1, p_2) a_{p_0}(-p_1, -p_2),  \\
 C_{p_0}(p_1, p_2) p_0^2   =  a_{p_0}(p_1, -p_2) a_{-p_0}(p_1, p_2),\\
 D_{p_0}(p_1, p_2) p_0^2   =  a_{p_0}(-p_1, p_2) a_{-p_0}(p_1, p_2).
\end{gather*}

\begin{rmk}
In particular, the sums in (\ref{EqnProdI})--(\ref{EqnProdIV}) run only through either the positive or the negative numbers, since $a_z(x,y) = 0$ if $\sgn(xyz) = -1$.
\end{rmk}

\begin{rmk}
It is also possible to obtain (\ref{EqnPrI})--(\ref{EqnPrIV}) from \cite[Proposition~4.10]{GroKoeKus}, using the pairing between~$M$ and $\left\{ \lambda(\omega) \mid \omega \in M_\ast \right\} \subseteq \hat{M}$, def\/ined by $\langle x, \lambda(\omega) \rangle = \omega(x)$.
 \end{rmk}

\begin{rmk}
In case of the group $SU(1,1)$ we know \cite[Chapter~6]{VilKli} that there exists an addition formula corresponding to the product formula, i.e.\ the product formula corresponds to the constant term in the addition formula. It would be of interest to obtain addition formulae corresponding to (\ref{EqnProdI})--(\ref{EqnProdIV}).
\end{rmk}

\subsection*{Acknowledgement}
The author likes to thank Erik Koelink for the useful discussions and Noud Aldenhoven for providing Fig.~\ref{PicSpectrum}. Also, the author benef\/its from a detailed referee report.

\newpage

\pdfbookmark[1]{References}{ref}
\LastPageEnding

\end{document}